\documentstyle[12pt]{article}
\newcommand{\sect}[1]{\section{#1}\setcounter{equation}{0}}

\font\mbn=msbm10 scaled \magstep1
\font\mbs=msbm7 scaled \magstep1
\font\mbss=msbm5 scaled \magstep1
\newfam\mbff
\textfont\mbff=\mbn
\scriptfont\mbff=\mbs
\scriptscriptfont\mbff=\mbss\def\mbf{\fam\mbff}
\def\Re{{\mbf R}}

\def\Z{{\mbf Z}}
\def\Co{{\mbf C}}

\def\N{{\mbf N}}
\def\H{{\mbf H}}
\newtheorem{Th}{Theorem}[section]
\newtheorem{Lm}[Th]{Lemma}
\newtheorem{C}[Th]{Corollary}
\newtheorem{Con}[Th]{Conjecture}
\newtheorem{D}[Th]{Definition}
\newtheorem{Proposition}[Th]{Proposition}
\newtheorem{R}[Th]{Remark}

\author{Alexander Brudnyi\thanks{Research supported in part by NSERC.
}\\
Department of Mathematics and Statistics\\
University of Calgary, Calgary\\
Canada\\
\\
Yuri Brudnyi\\
Department of Mathematics\\
Technion, Haifa\\
Israel}
\title{Metric Spaces with Linear Extensions Preserving Lipschitz Condition}
\date{}
\begin{document}
\maketitle
\begin{abstract}
{We study a new bi-Lipschitz invariant $\lambda(M)$ of a metric
space $M$; its finiteness means that Lipschitz functions on an
arbitrary subset of $M$ can be linearly extended to functions on
$M$ whose Lipschitz constants are expanded by a factor controlled
by $\lambda(M)$. We prove that $\lambda(M)$ is finite for several
important classes of metric spaces. These include metric trees of
arbitrary cardinality, groups of polynomial growth, Gromov-hyperbolic
groups, certain classes of Riemannian manifolds of
bounded geometry and the finite direct sums of arbitrary combinations
of these objects. On the other hand we construct an example of a
two-dimensional Riemannian manifold $M$ of bounded geometry for
which $\lambda(M)=\infty$.}
\end{abstract}
\tableofcontents
\sect{Introduction} {\bf A.} The concept of a metric space is
arguably one of the oldest and important in mathematics, whereas
analysis on metric spaces has been developed only within last few
decades. The important part of this development is devoted to the
selection and study of classes of metric spaces ''supporting''
certain basic analytic facts and theories known for $\Re^{n}$.
Some results and problems appearing in this area are presented, in
particular, in the surveys [CW], [Gr1] Appendix B (by Semmes) and
[HK]. The main topic of our paper also belongs to this direction
of research and is devoted to study of classes of metric spaces
possessing the following property.
\begin{D}\label{d0}
A metric space $(M,d)$ has the Lipschitz condition preserving linear
extension property (abbreviated ${\cal LE}$), if for each of its subspaces 
$S$ there is a linear continuous extension operator acting from $Lip(S)$ into
$Lip(M)$.
\end{D}
Here $Lip(S)$ is the space of real-valued functions on $S$ equipped with
the seminorm
\begin{equation}\label{e11}
|f|_{Lip(S)}:=\sup_{m'\neq m''}\frac{|f(m')-f(m'')|}{d(m',m'')}\ ;
\end{equation}
hence the linear operator $E:Lip(S)\to Lip(M)$ of this definition meets
the following conditions.
\begin{itemize}
\item[(a)] The restriction of $Ef$ to $S$ satisfies
\begin{equation}\label{e12}
Ef|_{S}=f\ ,\ \ \ \ \ f\in Lip(S)\ .
\end{equation}
\item[(b)]
The norm of $E$ given by
\begin{equation}\label{e12'}
||E||:=\sup\{|Ef|_{Lip(M)}\ :\ |f|_{Lip(S)}\leq 1\}
\end{equation}
is finite.
\end{itemize}
In the sequel the linear space of these operators will be denoted by
$Ext(S,M)$, that is to say,
\begin{equation}\label{e13}
Ext(S,M):=\{E\in {\cal L}(Lip(S),Lip(M))\ :\ E\ {\rm satisfies}\ (a)\
{\rm and}\ (b)\}\ .
\end{equation}

The adjective ``linear'' in Definition \ref{d0} drastically
changes the situation compared to that for nonlinear extensions
of Lipschitz functions. In fact, McShane [Mc] gave two
simple nonlinear formulas for the extension of functions $f\in Lip(S)$
preserving their Lipschitz constants. At the same year Kirszbraun
[Ki] proved existence of a Lipschitz constant preserving
(nonlinear) extension for maps between two Euclidean spaces; then
Valentine [V] remarked that this result remains true for the case
of general Hilbert spaces. Another generalization of Kirszbraun's
theorem was due to Lang and Schroeder [LSch] who proved such a
result for Lipschitz maps between metric path spaces with upper
and lower curvature bounds. Most of these results either fail to be
true or are unknown for the linear extension case (even for scalar
functions). For instance, Theorem \ref{te110} below presents an
example of a Riemannian two-dimensional manifold $\Sigma$ of
bounded geometry and a (metric) subspace $S$ such that
$Ext(S,\Sigma)=\emptyset$. (In the forthcoming paper [BB] we prove a similar 
result for infinite-dimensional Hilbert spaces.)

In the present paper we will study the following quantitative characteristic
of spaces with ${\cal LE}$.
\begin{D}\label{d2}
Given a metric space $(M,d)$ one defines its Lipschitz condition preserving
linear extension constant $\lambda(M)$ by
\begin{equation}\label{e14}
\lambda(M):=\sup_{S\subset M}\ \inf\{||E||\ :\ E\in Ext(S,M)\}\ .
\end{equation}
\end{D}

It is proved for a wide class of metric spaces, see Theorem
\ref{te14} below, that finiteness of (\ref{e14}) is equivalent to
the ${\cal LE}$ of $M$. In particular, $\Re^{n}$ equipped with an
arbitrary norm and the hyperbolic space $\H ^{n}$ with the inner
path (geodesic) metric belong to this class. In the former case,
the following estimate can be derived from the
classical Whitney extension theorem [W1]
\begin{equation}\label{w1}
\lambda(l_{1}^{n})\leq c^{n}
\end{equation}
where $c$ is some absolute (numerical) constant. The proof of (\ref{w1}) is
based on the Whitney covering lemma that is not true for a relatively simple
metric spaces, e.g., for $\H ^{n}$. Using a new approach based on a
quasi-isometric representation of $\H ^{n}$ as a space of balls in
$\Re^{n-1}$ with a corresponding path metric,  it was proved in
[BSh2, Proposition 5.33] that
\begin{equation}\label{w2}
\lambda(\H ^{n})\leq c^{n}
\end{equation}
with an absolute constant $c$.

In the present paper, we essentially enlarge the class of metric spaces with
finite $\lambda(M)$ and, in some cases, give even relatively sharp
estimates of this constant. For example, we show that for some $0<c<1$
$$
c\sqrt{n}\leq\lambda(\Z^{n})\leq 24\ \! n\ ,
$$
for $\Z^{n}$ regarded as an abelian group with the word metric,
and that the same upper estimate holds for an arbitrary Carnot
group of homogeneous dimension $n$.

One of the main tools of our approach is the {\em finiteness property}
of the characteristic (\ref{e14}), see Theorem \ref{te11} and
Corollaries \ref{c12} and \ref{c13} below, asserting, in particular, that
$$
\lambda(M)=\sup_{S}\lambda(S)
$$
where $S$ runs through all {\em finite point} subspaces of $M$
(with the induced metric). As the first consequence of this fact
we prove, see Theorem \ref{te14a} below, that for the direct sum
of arbitrary nontrivial metric trees ${\cal T}_{i}$ with path metrics
$d_{{\cal T}_{i}}$, $1\leq i\leq n$,
$$
C_{1}\sqrt{n}\leq\lambda(\oplus_{i=1}^{n}{\cal T}_{i})\leq C_{2}n
$$
with absolute constants $C_{i}>0$. This implies a similar estimate
for the Cayley graph of the direct product of free groups of
arbitrary cardinality.

The next important result, Theorem \ref{te17}, gives sufficient conditions
for finiteness of $\lambda(M)$ in the case of {\em locally doubling} metric
spaces $M$ (see corresponding definitions of this and other notions used here
in the next section). In particular, this implies the corresponding results
for a metric space of {\em bounded geometry} and for one framed by a group of
its isometries acting freely, properly and cocompactly, see Corollaries
\ref{c1.9} and \ref{c111b}.

The problem of nontrivial lower estimates of $\lambda(M)$ is unsolved
even for the ``relatively simple'' case of the $n$-dimensional Euclidean
space.\footnote{In the forthcoming paper [BB] we show that in this case
$\lambda(M)\geq c\sqrt[8]{n}$ for a numerical constant $c>0$.}
However,
we do prove such a result for the space $l_{p}^{n}$ with $p\neq 2$ based
on a sharp in order estimate of $\lambda_{conv}(l_{p}^{n})$, see
Theorem \ref{te1.11} below. Here $\lambda_{conv}$ is defined for metric
spaces $M$ with {\em convex structure}, e.g., for normed spaces, by
the same formula (\ref{e14}) but with $S$ running through {\em convex}
subsets of $M$.

Finally, Theorem \ref{te114} gives sufficient conditions for
finiteness of $\lambda(M)$ for a wide class of metric spaces that
includes, in particular, fractals, Carnot groups, groups of
polynomial growth, Gromov-hyperbolic groups and certain
Riemannian manifolds with curvature bounds, see section 3 below.\\
{\bf B.} The linear extension problem for spaces of continuous
functions was first studied by Borsuk [Bor] in 1933. Important
results in this area were obtained by Kakutani, Dugundji,
Lindenstrauss, Pe\l czynski and many other mathematicians, see [P]
and [BL, Chapters 2 and 3] and references therein. For the case of
uniformly continuous functions a negative result was proved by
Pe\l czynski [P, Remarks to $ \S 2$]. His argument, going back to
the paper [L] by Lindenstrauss, can be modified (cf. the proof of
Proposition \ref{p51} below) to establish that there is no 
linear bounded extension operator from $Lip(Y)$ into $Lip(X)$, if
$Y$ is a reflexive subspace of a Banach space $X$. It was noted in
the cited book [P] that ``... our knowledge of existence of linear
extension operators for uniformly continuous or Lipschitz
functions is rather unsatisfactory''. Not much research has,
however, been done in this area. Several important linear extension
results were however proved (but linearity there was not
formulated explicitly) in another area of research on Lipschitz
extensions. The main point there is to estimate (nonlinear)
Lipschitz extension constants for mappings from finite metric
spaces into Banach spaces, see [MP], [JL], [JLS]. In particular,
in the paper [JLS] by Johnson, Lindenstrauss and Schechtman the
proposed extension operator is linear and therefore their results
give for an $n$-dimensional Banach space $B$ the estimate
$\lambda(B)\leq Cn$, and for a finite metric space $M$ the
estimate $\lambda(M)\leq C\log(card M)$, where $C$ is an absolute
constant. Another important result was proved by Matou\v{s}ek
[Ma]; for scalar valued functions this gives an estimate of
$\lambda({\cal T})$ for an arbitrary metric tree ${\cal T}$ by
some universal constant.

For differentiable functions on $\Re^{n}$ a method of linear
extension was discovered by Whitney [W1] in 1934. It has been then
used in variety of problems of Analysis. To discuss the few
results in this field we recall that $C_{b}^{k}(\Re^{n})$ and
$C_{u}^{k}(\Re^{n})$ are spaces of $k$-times continuously
differentiable functions on $\Re^{n}$ whose higher derivatives
are, respectively, bounded or uniformly continuous. We also
introduce the space $C^{k,\omega}(\Re^{n})\subset C^{k}(\Re^{n})$
defined by the seminorm
\begin{equation}\label{wh1}
|f|_{C^{k,\omega}}:=\max_{|\alpha|=k}\sup_{x,y\in\Re^{n}}
\frac{|D^{\alpha}f(x)-D^{\alpha}f(y)|}{\omega(|x-y|)}\ .
\end{equation}
Here $\omega:\Re_{+}\to\Re_{+}$ is nondecreasing, equal to 0 at 0 and
concave; we will write $C^{k,s}(\Re^{n})$ for $\omega(t):=t^{s}$,
$0<s\leq 1$.\\
Finally, $\Lambda^{\omega}(\Re^{n})$ stands for the Zygmund space defined
by the seminorm
\begin{equation}\label{wh2}
|f|_{\Lambda^{\omega}}:=\sup_{x\neq y}\frac{|f(x)-2f(\frac{x+y}{2})+f(y)|}
{\omega(|x-y|)}\ ;
\end{equation}
here $\omega:\Re_{+}\to\Re_{+}$ is as in (\ref{wh1}), but we assume now that
$\omega(\sqrt{t})$ is concave.

Let now $S\subset\Re^{n}$ be an arbitrary closed subset and $X$ be one of
the above introduced function spaces. 
Then $X|_{S}$ denotes the linear space of
traces of functions from $X$ to $S$ endowed with the seminorm
\begin{equation}\label{wh3}
|f|_{X}^{S}:=\inf\{|g|_{X}\ :\ g|_{S}=f\}\ .
\end{equation}

Then {\em the linear extension problem} can be formulated as follows.

{\em Does there exist a linear continuous extension operator from
$X|_{S}$ into $X$?}\\
One can also consider the {\em restricted linear extension problem} with
$S$ belonging to a fixed class of closed subspaces of $\Re^{n}$.

Whitney's paper [W2] is devoted to a criterium for a function
$f\in C(S)$ with $S\subset\Re$ to belong to the trace space 
$C_{b}^{k}(\Re)|_{S}$ and
gives, in fact, a positive solution to the linear extension problem for
$C_{b}^{k}(\Re)$. It was noted in [BSh2] that Whitney's method gives
the same result for the spaces $C^{k,\omega}(\Re)$ and $C_{u}^{k}(\Re)$.

The situation for the multidimensional case is much more complicated.
The restricted problem, for the class of {\em compact} subsets of
$\Re^{n}$ was solved positively by \penalty-10000 
G. Glaeser [Gl] for the space
$C_{b}^{1}(\Re^{n})$ using a special construction of the geometry
of subsets in $\Re^{n}$. However, for the space
$C_{u}^{1}(\Re^{n})$, $n\geq 2$, the linear extension problem
fails to be true, see [BSh2, Theorem 2.5]. In [BSh2] (see also
[BSh1]) the linear extension problem was solved positively for the
spaces $C^{1,\omega}(\Re^{n})$ and $\Lambda^{\omega}(\Re^{n})$. A
recent breakthrough due to Ch.$\!$ Fefferman [F1] in the problem
of a constructive characterization of the trace space
$C^{k,1}(\Re^{n})|_{S}$, allowed him to solve the linear extension
problem for the space $C^{k,\omega}(\Re^{n})$, see [F2], [F3] and [F4].\\
{\bf C.} The paper is organized as follows.

Section 2 introduces basic
classes of metric spaces involved in our considerations and formulates the
main results, Theorems \ref{te11}, \ref{te14a}, \ref{te14},
\ref{te110}, \ref{te1.11} and their corollaries.

The next section presents
some important examples of metric spaces possessing ${\cal LE}$, while in
section 4 we discuss several open problems.

All the remaining sections are
devoted to proofs of the aforementioned main theorems and corollaries.

Finally, the Appendix presents an alternative proof of a
Kantorovich-Rubinshtein duality
type theorem used in the proof of Theorem \ref{te11}. \\
\\
{\bf Notations.} Throughout the paper we often suppress the symbol
$d$ in the notation $(M,d)$ and simply refer to $M$ as a metric
space. The same simplification will be used for all notation related
to $M$, e.g., we will write $Lip(M)$, $\lambda(M)$, see
(\ref{e11}) and (\ref{e14}), and use similar notations $Lip(S)$
and $\lambda(S)$ for $S\subset M$ regarded as a metric subspace
of $M$ (with the induced metric). Points of $M$ are denoted by
$m,m',m''$ etc, and $B_{r}(m)$ stands for the open ball of $M$
centered at $m$ and of radius $r$. We will write $M\in {\cal LE}$
for $M$ satisfying Definition \ref{d0}. Recall also that
$Ext(S,M)$ has already introduced by (\ref{e13}). Finally, we
define the {\em Lipschitz constant} of a map $\phi:(M,d)\to
(M_{1},d_{1})$ by
\begin{equation}\label{not}
|\phi|_{Lip(M,M_{1})}:=\sup_{m'\neq m''}
\frac{d_{1}(\phi(m'),\phi(m''))}{d(m',m'')}\ .
\end{equation}
Let us recall that $\phi$ is a {\em quasi-isometry} (or {\em bi-Lipschitz
equivalence}), if $\phi$ is a bijection and Lipschitz constants
for $\phi$ and $\phi^{-1}$ are finite. If, in addition,
$$
\max\{|\phi|_{Lip(M,M_{1})}\ , \ |\phi^{-1}|_{Lip(M_{1},M)}\}\leq C
$$
for some $C>0$, then $\phi$ is called a {\em $C$-isometry}
(and {\em isometry}, if $C=1$).
\sect{Formulation of the Main Results}
To formulate our first result we require the notion of a {\em dilation}.
This is a quasi-isometrty $\delta$ of $M$ such that the operator
$\Delta: Lip(M)\rightarrow Lip(M)$ given by
$$
(\Delta f)(m):=f(\delta(m))\ ,\ \ \ m\in M\ ,
$$
satisfies
\begin{equation}\label{dil}
||\Delta||\cdot ||\Delta^{-1}||=1\ .
\end{equation}
\begin{Th}\label{te11}
Assume that $S$ is a subspace of $M$ such that for some dilation
\penalty-10000 $\delta:M\to M$ we have
\begin{itemize}
\item[{\rm (a)}]
$S\subset\delta(S)$;
\item[{\rm (b)}]
$\cup_{j=0}^{\infty}\delta^{j}(S)$ is dense in $M$.
\end{itemize}
Then
$$
\lambda(M)=\sup_{F\subset S}\lambda(F)
$$
where $F$ runs through all finite point subspaces of $S$.
\end{Th}

Choosing $S=M$ and $\delta$ equal to the identity map we get from this
\begin{C}\label{c12}
\begin{equation}\label{e15}
\lambda(M)=\sup_{F}\lambda(F)
\end{equation}
where $F$ runs through all finite point subspaces of $M$.
\end{C}

Together with Theorem \ref{te11} this immediately implies
\begin{C}\label{c13}
Let $S\subset M$ satisfy the assumptions of Theorem \ref{te11}. Then
$$
\lambda(M)=\lambda(S)\ .
$$
\end{C}

The results presented above will be used in almost all subsequent
proofs. As the first application we give a rather sharp estimate
of $\lambda(M)$ for $M$ being the direct sum of metric trees. To
formulate the result let us recall the corresponding notions.

A {\em tree} ${\cal T}$ is a connected graph with no cycles, see,
e.g. [R, Ch.9] for more details. We turn ${\cal T}$ into a
{\em path metric space} by identifying each edge $e$ with a
bounded interval of $\Re$ of length $l(e)$ and then determining
the distance between two points of the 1-dimensional CW-complex
formed by these edges to be the infimum of the lengths of the paths
joining them. Since every two vertices of a tree can be joined by a
unique path, the metric space $({\cal T}, d_{{\cal T}})$ obtained
in this way is, in fact, a geodesic space, see e.g. [BH, pp. 8-9].

Let now $(M_{i},d_{i})$, $1\leq i\leq n$, be metric spaces. Their
direct $p$-sum $\oplus_{p}\{(M_{i},d_{i})\}_{1\leq i\leq n}$ is a metric 
space with the
underlying set $\displaystyle \prod_{i=1}^{n} M_{i}$ and a metric
$d$ given by
\begin{equation}\label{2.2'}
d(m,m'):=\left(\sum_{i=1}^{n}d_{i}(m_{i},m_{i}')^{p}\right)^{1/p};
\end{equation}
here $m=(m_{1},\dots,m_{n})$, $m'=(m_{1}',\dots,m_{n}')$.

\begin{Th}\label{te14a}
Let ${\cal T}_{i}$ be a nontrivial metric tree, $1\leq i\leq n$. Then for
$p=1,\infty$
$$
c_{0}\sqrt{n}\leq\lambda(\oplus_{p}\{{\cal T}_{i}\}_{1\leq i\leq n})\leq cn
$$
where $c_{0}, c$ are absolute constants.
\end{Th}

The basic fact of independent interest used along with Corollary
\ref{c12} in the proof of this theorem asserts that every infinite
metric tree with uniformly bounded vertex degrees 
admits a quasi-isometric embedding into the hyperbolic plane with
distortion\footnote{recall that {\em distortion} of a
quasi-isometry $\phi: M_{1}\to M_{2}$ of metric spaces
is defined by $|\phi|_{Lip(M_{1},M_{2})}\cdot 
|\phi^{-1}|_{Lip(M_{2},M_{1})}$.} bounded by a numerical constant.
It seems to be strange to use here the
hyperbolic plane instead of a Euclidean space of some dimension.
Strikingly, by a result of Bourgain [Bou]  this cannot be done
even if we use an infinite dimensional Hilbert space.

For $n=1$ the above result was proved by Matou\v{s}ek [Ma] by
another method. It is worth noting that an important class of spaces,
Gromov-hyperbolic spaces of bounded geometry, have metric structure close to 
that of metric trees. This implies the corresponding Lipschitz
extension result for spaces of this class, see Corollary \ref{2.12''} 
below.

Our next result relates the ${\cal LE}$ of $M$ to the finiteness
of $\lambda(M)$. For its formulation we introduce the following
two classes of metric spaces.
\begin{D}\label{def1}
(a) A metric space $M$ is said to be proper (or boundedly compact), if every
closed ball in $M$ is compact.\\
(b) A metric space $M$ has the weak transition property (WTP), if
for some $C\geq 1$ and every finite set $F$ and open ball $B$ in
$M$ there is a $C$-isometry $\sigma:M\rightarrow M$ such that
$$
B\cap\sigma(F)=\emptyset\ .
$$
\end{D}
\begin{Th}\label{te14}
Assume that $M$ is either proper or has the WTP. Then the ${\cal LE}$ of $M$ 
is equivalent to the finiteness of $\lambda(M)$.
\end{Th}

We now shall discuss some general conditions under which a metric space
possesses the required extension property. For this purpose we use a
modification of the well-known {\em doubling condition}.
\begin{D}\label{def2}
A metric space $M$ 
is  locally doubling, if for some $R>0$ and integer $N$ each
ball of radius $r<R$ in $M$ can be covered by at most $N$ balls of radius
$r/2$.
\end{D}

The class of such spaces will be denoted by ${\cal D}(R,N)$. The class
of {\em doubling metric spaces} is then
$\bigcup_{N}(\ \!\! \bigcap_{R>0}{\cal D}(R,N))$.
We will write $M\in {\cal D}(N)$, if the assumption of Definition \ref{def2}
holds for all $r<\infty$.

The second notion that will be used is introduced by
\begin{D}\label{def3}
A set $\Gamma\subset M$ is said to be an $R$-lattice, if the family of
open balls $\{B_{R/2}(\gamma)\ :\ \gamma\in\Gamma\}$ forms a cover of $M$,
while the balls $B_{cR}(\gamma)$, $\gamma\in\Gamma$, are pairwise
disjoint for some $c=c_{\Gamma}\in (0,1/4]$.
\end{D}

The existence of $R$-lattices follows easily from Zorn's lemma.
\begin{Th}\label{te17}
Assume that a metric space $M\in {\cal D}(R,N)$ and $\Gamma\subset M$ is
an $R$-lattice. Assume also that the constants $\lambda(\Gamma)$ and
\begin{equation}\label{e16}
\lambda_{R}:=\sup\{\lambda(B_{R}(m))\ :\ m\in M\}
\end{equation}
are finite.

Then $\lambda(M)$ is bounded by a constant
depending only on $\lambda(\Gamma)$, $\lambda_{R}$, $c_{\Gamma}$,
$R$ and $N$.
\end{Th}

In order to formulate a corollary of this result we introduce a subclass
of the class $\cup_{N,R}{\cal D}(R,N)$ consisting of metric spaces of
{\em bounded geometry}, cf. the corresponding definition in [CG] for
the case of Riemannian manifolds.
\begin{D}\label{def1.8}
A metric space $M$ is of bounded geometry with parameters $n\in\N$, $R$,
$C>0$ (written $M\in {\cal G}_{n}(R,C)$), if each open ball of radius
$R$ in $M$ is $C$-isometric to a subset of $\Re^{n}$.
\end{D}

Let us note that if $B_{R}(m)$ is $C$-isometric to a subset $S$ of
$\Re^{n}$, then
$$
C^{-2}\cdot\lambda(S)\leq\lambda(B_{R}(m))\leq C^{2}\cdot\lambda(S)\ ,
$$
and by the classical Whitney extension theorem, see, e.g., [St, Ch.6],
$\lambda(S)\leq\lambda(\Re^{n})<\infty$. So the previous theorem
leads to
\begin{C}\label{c1.9}
Let $M\in {\cal G}_{n}(R,C)$. Then $\lambda(M)$ is finite if and only if
for some $R$-lattice $\Gamma$ we have
$$
\lambda(\Gamma)<\infty\ .
$$
\end{C}

To formulate the second corollary we recall the definition of Gromov
hyperbolicity [Gr3]. We choose a definition  attributed (by Gromov) to
Rips (see equivalent formulations in [BH] Chapter 3).

Let $(M,d)$ be a {\em geodesic metric space}; this means that every two 
points $m, n\in M$ can be joined by a geodesic segment, the image of a map
$\gamma:[0,a]\to M$ such that $\gamma(0)=m$, $\gamma(a)=n$ and
$d(\gamma(t),\gamma(s))=|t-s|$ for all $t$, $s$ in $[0,a]$.
\begin{D}\label{2.12'}
A geodesic metric space $M$ is said to be $\delta$-hyperbolic, if every
geodesic triangle in $M$ is $\delta$-slim, meaning that each of its sides
is contained in the $\delta$-neighbourhood of the union of the remaining
sides.
\end{D}
\begin{C}\label{2.12''}
Let $M$ be the finite direct $p$-sum of hyperbolic metric spaces of
bounded geometry, $1\leq p\leq\infty$. Then $\lambda(M)$ is finite.
\end{C}

The next consequence concerns a {\em path metric space}\footnote{i.e., the
distance between every pair of points equals the infimum of the lengths of
curves joining the points.} with a {\em group action}. For its formulation
we need
\begin{D}\label{d111a} (see, e.g., [BH, p.131]).
A subgroup $G$ of the group of isometries of a metric space $M$ acts
properly, freely  and cocompactly on $M$, if
\begin{itemize}
\item[(a)]
for every compact set $K\subset M$ the set
$\{g\in G\ :\ g(K)\cap K\neq\emptyset\}$ is finite;
\item[(b)]
for every point $m\in M$ the identity $g(m)=m$ implies that $g=1$;
\item[(c)]
there is a compact set $K_{0}\subset M$ such that
\begin{equation}\label{e17a}
M=G(K_{0})\ .
\end{equation}
\end{itemize}
\end{D}

By the \v{S}varc-Milnor lemma, see, e.g., [BH, p.140], the group
$G$ of this definition is {\em finitely generated} whenever $M$ is
a path space. If $A$ is a (finite) generating set for $G$, then
$d_{A}$ stands for the {\em word metric} on $G$ determined by $A$,
see, e.g., [Gr1,p.89]. Replacing $A$ by another (finite)
generating set one obtains a corresponding word metric
bi-Lipschitz equivalent to $d_{A}$. Therefore the ${\cal LE}$ of
$G$ regarded as a metric space in this way does not depend on the
choice of $A$.
\begin{C}\label{c111b}
Let $M$ be a path space  framed by a group $G$ acting on $M$ by
isometries. Assume that
\begin{itemize}
\item[(a)]
$M$ is a metric space of bounded geometry;
\item[(b)]
$G$ acts on $M$ properly, freely and cocompactly.
\end{itemize}
Then $\lambda(M)$ is finite if and only if $\lambda(G,d_{A})$ is.
\end{C}

In view of Corollary \ref{c1.9}
it would be natural to conjecture that $\lambda(M)$
is finite for every $M$ of bounded geometry. The following
counterexample disproves this assertion.
\begin{Th}\label{te110}
There exists a connected two-dimensional metric space $M_{0}$ of bounded
geometry such that
$$
Ext(S,M_{0})=\emptyset
$$
for some subset $S\subset M_{0}$.
\end{Th}

The basic step used in our construction of $M_{0}$ is the following result
of independent interest. For its formulation we introduce the functional
$$
\lambda_{conv}(M):=\sup _{C}\lambda(C,M)
$$
where $C$ runs through all {\em convex subsets} of a {\em normed linear
space} $M$. Here we set
$$
\lambda(S,M):=\inf\{||E||\ :\ E\in Ext(S,M)\}\ .
$$
\begin{Th}\label{te1.11}
There exists an absolute constant $c_{0}>0$ such that for all $n$ and
$1\leq p\leq\infty$
$$
c_{0}\leq n^{-\left|\frac{1}{p}-\frac{1}{2}\right|}
\cdot\lambda_{conv}(l_{p}^{n})
\leq 1\ .
$$
\end{Th}

Let us now return to metric spaces of bounded geometry. Corollary
\ref{c1.9} tells us that the existence of the desired extension
property is reduced to that for lattices. The example of Theorem
\ref{te110} makes the following conjecture to be rather plausible.
\begin{Con}\label{con1}
A lattice $\Gamma\subset M$ has the ${\cal LE}$, if it is uniform.

The latter means
that for some increasing function $\phi_{\Gamma}:\Re_{+}\rightarrow
\Re_{+}$ and constant $0<c\leq 1$ the number of points of
$\Gamma\cap B_{R}(m)$ for every $R>0$ and $m\in\Gamma$ satisfies
\begin{equation}\label{e17}
c\phi_{\Gamma}(R)\leq |\Gamma\cap B_{R}(m)|\leq\phi_{\Gamma}(R)\ .
\end{equation}
\end{Con}

We confirm this conjecture for lattices of {\em polynomial
growth}, i.e., for \penalty-10000 $\phi_{\Gamma}(R)=aR^{n}$ for
some $a, n\geq 0$ and for some other lattices including even those
of {\em exponential growth}. These will follow from an extension
result presented below. In its introduction we use the notion of a measure
{\em doubling at a point} $m$. 
This is a nonnegative Borel measure $\mu$ on $M$
such that every open ball centered at $m$ is of finite strictly positive
$\mu$-measure and the {\em doubling constant}
$$
D_{m}(\mu):=\sup_{R>0}\frac{\mu(B_{2R}(m))}{\mu(B_{R}(m))}
$$
is finite. If, in addition,
$$
D(\mu):=\sup_{m\in M}D_{m}(\mu)<\infty,
$$
the $\mu$ is said to be a {\em doubling measure}.

A metric space endowed with a doubling measure is said to be of 
{\em homogeneous type} [CW]. 

Our basic class of metric spaces is presented by
\begin{D}\label{d218}
A metric space $(M,d)$ is said to be of pointwise homogeneous type if
there is a fixed family $\{\mu_{m}\}_{m\in M}$ of Borel measures on 
$M$ satisfying the following properties.
\begin{itemize}
\item[(i)]
\underline{Uniform doubling condition:}

$\mu_{m}$ is doubling at $m$ and
\begin{equation}\label{eq27}
D:=\sup_{m\in M}D_{m}(\mu_{m})<\infty.
\end{equation}
\item[(ii)]
\underline{Consistency with the metric:}

For some constant $C>0$ and all $m_{1},m_{2}\in M$ and $R>0$
\begin{equation}\label{eq28}
|\mu_{m_{1}}-\mu_{m_{2}}|(B_{R}(m))\leq\frac{C\mu_{m}(B_{R}(m))}{R}
d(m_{1},m_{2})
\end{equation}
where $m=m_{1}$ or $m_{2}$.
\end{itemize}
\end{D}
\begin{R}\label{r219}
{\rm The conditions (\ref{eq27}), (\ref{eq28}) hold trivially 
for $M$ equipped with a doubling measure $\mu$ (i.e., in this case
$\mu_{m}=\mu$ for all $m$). So metric spaces of homogeneous type belong to
the class introduced by this definition.}
\end{R}
\begin{Th}\label{te114}
If $M$ is of pointwise homogeneous type 
with the optimal constants $C$ and $D$, then 
the following inequality
\begin{equation}\label{e2.9'}
\lambda(M)\leq k_{0}(C+1)(\log_{2}D+1)
\end{equation}
holds with some numerical constant $k_{0}$.
\end{Th}

In particular, for a metric space $M$ of homogeneous type 
(when $C=0$) $D$ is greater than, say, $\sqrt[4]{2}$, see, e.g., [CW],
and (\ref{e2.9'}) gives the inequality
$$
\lambda(M)\leq 2k_{0}\log_{2}D.
$$
\begin{R}\label{2.22'}
{\rm The last inequality can be easily derived from Theorem 1.4 of the paper
[LN] by Lee and Naor on ``absolute'' Lipschitz extendability of doubling
metric spaces.\footnote{A metric space is {\em doubling}, if for every $R>0$
there is an integer $N>1$ such that any closed ball of radius $2R$ can be
covered by $N$ balls of radius $R$.} Their proof is based on a probabilistic
argument. Using a modification of the proof of Theorem \ref{te114} one can
give a constructive proof of the Lee-Naor result.}
\end{R}

We formulate several consequences of Theorem \ref{te114}.
\begin{D}\label{de221}
A metric space $(M,d)$ is said to be of pointwise
$(a,n)$-homogeneous type, $n\geq 0$, $a\geq 1$, with respect to a
family of Borel measures $\{\mu_{m}\}_{m\in M}$ on $M$, if it
satisfies condition (\ref{eq28}) and the condition
\begin{equation}\label{eq29}
\frac{\mu_{m}(B_{lR}(m))}{\mu_{m}(B_{R}(m))}\leq al^{n}
\end{equation}
for arbitrary $l\geq 1$, $m\in M$ and $R>0$.
\end{D}
\begin{C}\label{co222}
If $M$ is of pointwise $(a,n)$-homogeneous type, then
\begin{equation}\label{eq210}
\lambda(M)\leq K_{0}(C+1)a^{2}(n+1)
\end{equation}
where $K_{0}$ is a numerical constant $(<225)$ and $C$ is the
constant in (\ref{eq28}).
\end{C}

Let us note that for a metric space of homogeneous type condition
(\ref{eq28}) trivially holds, and we can take $C=0$ in
(\ref{eq210}).

We now single out a special case of the above result with a better
estimate of $\lambda(M)$. Specifically, suppose now that for all balls
in $M$
\begin{equation}\label{eq211}
\mu_{m}(B_{R}(m))=\gamma R^{n}\ ,\ \ \ \gamma, n>0\ .
\end{equation}
Under this assumption the following holds.
\begin{C}\label{c116}
$$
\lambda(M)\leq 24(n+C)\ .
$$
\end{C}

Finally, we establish the finiteness of $\lambda(M)$ for a metric
space being the direct sum of spaces of pointwise homogeneous type. Unlike
the situation for spaces of homogeneous type this case requires an
additional restriction on the related families of Borel
measures introduced as follows.
\begin{D}\label{de224}
A family of measures $\{\mu_{m}\}$ on a metric space $M$ is said
to be $K$-uniform $(K\geq 1)$, if for all $m_{1},m_{2}$ and $R>0$
$$
\mu_{m_{1}}(B_{R}(m_{1}))\leq K\mu_{m_{2}}(B_{R}(m_{2}))\ .
$$
\end{D}
\begin{Th}\label{c118}
Let $M_{i}$ be of pointwise homogeneous type with respect to a $K_{i}$-uniform
family of Borel measures $\{\mu_{m}^{i}\}_{m\in M_{i}}$ satisfying
conditions of Definition \ref{d218} with the optimal
constants $D_{i}$, $C_{i}$, $1\leq i\leq N$. Then the following inequality
\begin{equation}\label{e2.13'}
\lambda(\oplus_{p}\{(M_{i},d_{i})\}_{1\leq i\leq N})\leq c_{0}
(\widetilde C_{p}+1)
(\log_{2}D+1) 
\end{equation}
is true with
$$
D:=\prod_{i=1}^{N}D_{i},\ \ \ \widetilde C_{p}:=\left(\sum_{i=1}^{N}
C_{i}^{q}\right)^{1/q}\prod_{i=1}^{N}K_{i}.
$$
Here $c_{0}$ is a numerical constant and $q$ is the exponent conjugate to 
$p$, i.e., $\frac{1}{p}+\frac{1}{q}=1$.

If, in particular, (\ref{eq211}) holds for $M_{i}$ with $n=n_{i}$
and $\gamma=\gamma_{i}$, $1\leq i\leq N$, then
\begin{equation}\label{eq212}
\lambda(\oplus_{\infty}\{(M_{i},d_{i})\}_{1\leq i\leq N})\leq 
24\sum_{i=1}^{N}(n_{i}+C_{i}).
\end{equation}
\end{Th}

The extension results of this section are true for Banach-valued 
Lipschitz functions, if the Banach space is complemented in its second
dual space (e.g., dual Banach spaces possess this property [Di]).
This can be derived straightforwardly from the scalar results.
However, the Banach-valued version of Theorem \ref{te114} is true without any
restriction. It can be established by an appropriate modification of
the proof presented here.  This and other results in that direction
will be presented in a forthcoming paper.
\sect{Examples}
{\bf 3.1. Groups with a metric space structure.}\\
{\bf A.} {\em Carnot Groups} (see [FS] and [He] for basic facts).

A {\em Carnot group} (also known as a {\em homogeneous group})
is a connected and simply connected real Lie group $G$
whose Lie algebra $g$ admits a stratification
\begin{equation}\label{21}
g=\bigoplus_{i=1}^{m}V_{i}\ \ \ {\rm with}\ \ \ [V_{i},V_{i}]=V_{i+1}\ ;
\end{equation}
here $V_{m+1}=\{0\}$ and $V_{m}\neq\{0\}$.\\
Being nilpotent, $G$ is diffeomorphic to $\Re^{n}$ with $n:=dim\ \! G$.
Together with the topological dimension $n$ an important role is played by 
the
homogeneous dimension of $G$ given by
\begin{equation} \label{22}
dim_{h}\ \!G:=\sum_{j=1}^{m}j\ \!dim\ \! V_{j}\ .
\end{equation}
The group $G$ can be equipped with a left-invariant
({\em Carnot-Caratheodory}) metric $d$ for which the ball
$B_{r}(x):=\{y\in G\ :\ d(x,y)<r\}$ satisfies
$$
|B_{r}(x)|=r^{Q}\ ,\ \ \ x\in G\ ,\ r>0\ .
$$
Here $|\cdot |$ is the (normed) left-invariant Haar measure on $G$ and
$Q:=dim_{h}\ \!G$. Therefore Corollary \ref{c116} immediately implies that
\begin{equation}\label{23}
\lambda(G,d)\leq 24\ \!dim_{h}\ \!G\ .
\end{equation}

The simplest example of a Carnot group is $\Re^{n}$. In this case,
$dim_{h}\ \!G=dim\ \!G=n$, and an arbitrary Banach norm on
$\Re^{n}$ defines a Carnot-Caratheodory metric. This gives the
aforementioned extension result of [JLS] with a better constant.
In particular, (\ref{23}) and Theorem \ref{te1.11} yield
\begin{equation}\label{24}
c_{0}\
\!n^{\left|\frac{1}{p}-\frac{1}{2}\right|}\leq\lambda(l_{p}^{n})\leq
24 n
\end{equation}
with $c_{0}>0$ independent of $n$ and $p$.

Another interesting example of a Carnot group is the Heisenberg
group $H_{n}(\Re)$ that, as a set, is equal to
$\Re^{n}\times\Re^{n}\times\Re$. The group operation is defined by
\begin{equation}\label{25}
(x,y,t)\cdot (x',y',t')=(x+x',y+y',t+t'+<x',y>)
\end{equation}
where the first two coordinates are vectors in $\Re^{n}$ and
$<\cdot,\cdot>$ is the standard scalar product. The topological
dimension of $H_{n}(\Re)$ is clearly $2n+1$ while its homogeneous
dimension equals $2n+2$. Finally, a Carnot-Caratheodory metric $d$
is given by
\begin{equation}\label{26}
d((x,y,t),(x',y',t')):=|(x,y,t)^{-1}\cdot (x',y',t')|
\end{equation}
where $|(x,y,t)|:=(<x,x>^{2}+<y,y>^{2}+t^{2})^{1/4}$.

For the metric space $(H_{n}(\Re),d)$ inequality (\ref{23}) gives
an upper bound $\lambda(H_{n}(\Re))\leq 48(n+1)$. It is
interesting to note that the Whitney extension method does not
work even to prove that $\lambda(H_{n}(\Re))<\infty$. In fact, its
basic geometric ingredient, Whitney's covering lemma cannot be
proved in this setting in a form allowing the required Lipschitz
partition of unity.

Finally, consider the discrete subgroup $H_{n}(\Z)$ of $H_{n}(\Re)$
consisting of elements of the set $\Z^{n}\times\Z^{n}\times\Z$. It is easily
seen that the map $\delta:(x,y,t)\mapsto\frac{1}{2}(x,y,t)$ is a dilation
in the sense used in Theorem \ref{te11}. Since
$\delta(H_{n}(\Z))\supset H_{n}(\Z)$ and $\delta^{j}H_{n}(\Z)$ is dense
in $H_{n}(\Re)$, Theorem \ref{te11} implies that
\begin{equation}\label{27}
\lambda(H_{n}(\Z))=\lambda(H_{n}(\Re))\ .
\end{equation}
{\bf B.} {\em Groups of Polynomial Growth.}

Let $G$ be a finitely generated group with the word metric $d_{A}$
associated with a set of generators $A$. This group is said to be of
{\em polynomial growth}, if for every $R>0$ the number of elements in a
ball of radius $R$ is bounded by $cR^{n}$ with fixed constants
$c,n>0$. By the Gromov result [Gr2] such a group is virtually nilpotent and
therefore by [B] for every ball $B_{R}(g_{0})$ the inequality
\begin{equation}\label{28}
c_{1}R^{Q}\leq |B_{R}(g_{0})|\leq c_{2}R^{Q}
\end{equation}
is true. Here $|\cdot|$ is the counting measure, $Q$ is the homogeneous
degree (\ref{22}) of the Zariski closure of the maximal torsion free
nilpotent subgroup of $G$ and $c_{1},c_{2}>0$ depend only on $G$.
Then Theorem \ref{te114} implies for this group the inequality
\begin{equation}\label{29}
\lambda(G)\leq c
\end{equation}
with $c$ depending on $G$.
In case $G$ being torsion free nilpotent,
the constants in (\ref{28}) depend only on
$Q$, see [B], and $c$ in (\ref{29}) does, as well.

For the special case of the abelian group $\Z^{n}$ this result can be 
sharpened.
In this case we use another representation of the metric space $(G,d_{A})$
related to the {\em Cayley graph} ${\cal C}_{A}(G)$. The latter is a metric
graph
whose vertices are in a one-to-one correspondence with elements of $G$ and
which has edges $e_{a}$ of length one joining each $g\in G$ with
$ga$, $a\in A\cup A^{-1}$. The metric subspace $G\subset {\cal C}_{A}(G)$ is
then isometric to $(G,d_{A})$, see, e.g. [BH, p.8]. For $\Z^{n}$ with
the set of generators $A:=\{a_{1},\dots,a_{n}\}$ being the standard basis of
$\Re^{n}$ the Cayley graph ${\cal C}_{A}(\Z^{n})$ is the 1-dimensional
CW-complex
determined by $\Z^{n}$ with the metric induced from $l_{1}^{n}$. Therefore
$(\Z^{n},d_{A})$ coincides with $\Z_{1}^{n}\subset l_{1}^{n}$ and the
application of Theorem \ref{te11} with the dilation
$\delta:x\mapsto\frac{1}{2}x$, $x\in\Re^{n}$, yields
$$
\lambda(\Z^{n},d_{A})=\lambda(l_{1}^{n})\ .
$$
This and Corollary \ref{c116} immediately imply that
\begin{equation}\label{210}
c\sqrt{n}\leq\lambda(\Z^{n},d_{A})\leq 24 n
\end{equation}
with $c$ independent of $n$.\\
{\bf C.} {\em Gromov-Hyperbolic Groups.}

A finitely generated group $G$ is called Gromov-hyperbolic if 
its Cayley graph with respect to some finite generating set is a Gromov
hyperbolic metric space.
Every such group is finitely presented, and conversely,
in a certain statistical sense, almost every finitely presented group is
hyperbolic [Gr3]. It is also known that any infinite, non-virtually cyclic 
hyperbolic group is of exponential growth. 
On the other hand, Corollary \ref{2.12''} asserts that $\lambda(G)$ is 
finite for these groups $G$. We now present several examples of hyperbolic 
groups.
\\
(1) Any finite group is hyperbolic.\\
(2) Any free group of finite rank is $0$-hyperbolic.\\
(3) The fundamental group of a compact Riemann surface is hyperbolic. \\
(4) A discrete cocompact group $G$ of isometries of the hyperbolic 
$n$-space $\H^{n}$ is hyperbolic, see [Bo] and [GrP] for
existence of such arithmetic and nonarithmetic groups $G$.\\
{\bf D.} {\em Free Groups.}

Let ${\cal F}(A)$ be a free group with the set of generators $A$
of arbitrary cardinality. It is easily seen that the Cayley graph
of ${\cal F}(A)$ is a metric tree rooted at the unit of ${\cal
F}(A)$ (the empty word). Hence ${\cal F}(A)$ equipped with the
word metric $d_{A}$ is an infinite rooted metric tree with all
edges of length one. In turn, the direct product $\prod_{i=1}^{n}{\cal
F}(A_{i})$ with the word metric generating by the generating set
$\prod_{i=1}^{n}A_{i}$ is isometric to 
$\oplus_{1}\{({\cal F}(A_{i}), d_{A_{i}})\}_{1\leq i\leq n}$.
Therefore Theorem \ref{te14a} immediately implies that
$$
c_{1}\sqrt{n}\leq\lambda(\prod_{i=1}^{n}{\cal F}(A_{i}))\leq
c_{2}n
$$
with $0<c_{1}<c_{2}$ independent of $n$.

We conjecture that this quantity is actually 
equivalent to $n$ as $n\to\infty$.\\
{\bf 3.2. Riemannian Manifolds.}\\
{\bf A.} {\em Nilpotent Lie Groups.}

Let $G$ be a simply connected real nilpotent Lie group of
dimension $N$ equipped with a left-invariant Riemannian metric.
Unlike Carnot groups the unit ball of $G$ is far from being a
Euclidean ellipsoid, see Figure 1 in [K] plotting the ball of a
big radius for the Heisenberg group $H^{3}(\Re)$ equipped with a
Riemannian metric. Nevertheless, it was shown in that paper that
the volume of the ball $B_{R}(g)$ of $G$ (with respect to the
geodesic metric) satisfies the inequality
$$
aR^{Q}\leq v(B_{R}(g))\leq bR^{Q}\ ,\ \ \ 1<R<\infty\ ;
$$
here $a,b>0$ are independent of $R$ and $g$, and $Q$ is the
homogeneous dimension of $G$, see (\ref{22}). A similar inequality
with $N$ instead of $Q$ holds for $0\leq R\leq 1$ (this follows
from the definition of the Riemannian structure on $G$). Hence the
dilation function for the measure $v$ is bounded by $(b/a)l^{Q}$
for $l>1$, and the metric space $G$ is of homogeneous type with respect to
the measure $v$. Then the conditions of Corollary \ref{co222} hold
with $C=0$ and $n=Q$. This implies the inequality
$$
\lambda(G)\leq K_{0}(b/a)^{2}(Q+1)
$$
with $K_{0}<225$.\\
{\bf B.} {\em Regular Riemannian Coverings.}

Let $p:M\to C$ be a regular Riemannian covering of a compact Riemannian
manifold $C$.
By definition, the deck transformation group of this covering $G_{p}$
acts on $M$ by isometries properly, freely and cocompactly. Hence by
Corollary \ref{c111b} $\lambda(M)$ is finite, if $\lambda(G_{p})$ is. For
instance, finiteness of the latter is true for $G_{p}$ being one of the
discrete groups presented in section 3.1. In particular, if $C$ is a compact
$n$-dimensional Riemannian manifold whose mean curvature tensor
is everywhere positive definite the group $G_{p}$ is of a polynomial growth
[M]. Hence $\lambda(M)<\infty$ in this case.
Suppose now that $C$ is a compact
Riemannian manifold whose sectional curvature is bounded above by a negative
number. Then its fundamental group $\pi_{1}(C)$ is Gromov-hyperbolic,
see, e.g., [BH, p.448]. Moreover, this group acts properly, freely and
cocompactly on the universal covering $M$ of $C$. Hence $\lambda(M)<\infty$
by Corollary \ref{c111b}.\\
{\bf C.} {\em Riemannian Manifolds of Nonnegative Ricci
Curvature.}

Let $M_{n}$ be a complete noncompact $n$-dimensional Riemannian
manifold regarded as a metric space with the geodesic metric.
Assume that the Ricci curvature of $M_{n}$ is nonnegative. Then by
the Laplacian comparison theorem, see, e.g., [Gr, p.283], the
volume of its balls satisfies
$$
\frac{v(B_{R_{2}}(m))}{v(B_{R_{1}}(m))}\leq
\left(\frac{R_{2}}{R_{1}}\right)^{n}\ ,\ \ \ 0<R_{1}\leq R_{2}\ .
$$
Hence $M_{n}$ equipped with the measure $v$ is of homogeneous type and the
dilation function $D(l)$ of $v$ is bounded by  $l^{n}$. Applying
Corollary \ref{co222} with $C=0$ and $a=1$, see (\ref{eq210}), we
get
$$
\lambda(M_{n})\leq K_{0}(n+1)
$$
with $K_{0}<225$.\\
{\bf D.} {\em Riemannian Manifolds of Pinched Negative Sectional Curvature.}

Let $M$ be a complete, simply-connected Riemannian manifold whose sectional 
curvature $\kappa$ satisfies $-b^{2}<\kappa<-a^{2}<0$ for some $a,b\in\Re$.
Then $M$ is a Gromov-hyperbolic metric space, see, e.g., [BH].
Rauch's comparison theorem implies also that $M$ is of bounded geometry,
see, e.g., [CE]. Now, application of Corollary \ref{2.12''} immediately yields
the finiteness of $\lambda(M)$.\\
{\bf E.} {\em Other Riemannian Manifolds.}

Let $H_{\rho}^{n+1}$ be a complete Riemannian manifold with the underlying
set \penalty-10000
$\Re_{+}^{n+1}:=\{(x,t)\in\Re^{n}\times\Re\ :\ t>0\}$ and the
Riemannian metric
$$
ds^{2}:=\rho(t)^{-2}(dx_{1}^{2}+\dots +dx_{n}^{2}+dt^{2})\ .
$$
We assume that $\rho$ is continuous and nondecreasing and
$\rho(0)=0$, while $\rho(t)/t^{2}$ is nonincreasing. We endow
$H_{\rho}^{n+1}$ with the inner (geodesic) metric and show that
the metric space obtained possesses the ${\cal LE}$. For this goal
we use the collection of measures $\{\mu_{m}\ \! :\ \!m\in
H_{\rho}^{n+1}\}$ introduced in [BSh2]. They are given by the
formula
$$
\mu_{m}(U):=\int_{U}\chi(t-s)\ \!\widehat\rho(s)\prod_{i=1}^{n}
\widehat\rho(t+|y_{i}-x_{i}|)\ \!dy_{1}\cdots dy_{n}\ \!ds
$$
provided that $m=(x,t)\in\Re_{+}^{n+1}$ and
$U\subset\Re_{+}^{n+1}$. Here $\widehat\rho:=1/\rho$
and $\chi$ is the Heaviside function, indicator of $[0,\infty)$.
It was proved in [BSh2, pp.537-540] that these measures satisfy
the conditions of Definition \ref{d218}. Hence the application of
Theorem \ref{te114} yields the inequality
$$
\lambda(H_{\rho}^{n+1})<\infty\ .
$$
\sect{Open Problems}
(a) Does there exist a metric space $M\in {\cal LE}$ but with $\lambda(M)=\infty$?
\\
(b) It follows from Theorem \ref{te1.11} that
$$
\lambda(l_{p})=\infty\ ,\ \ \ {\rm if}\ \ \ p\neq 2\ .
$$
Is the same true for infinite dimensional Hilbert spaces?$\!$
\footnote{
The answer is positive. In the forthcoming paper [BB] we show that
$\lambda(X)=\infty$ for an arbitrary infinite dimensional Banach space $X$.}\\
(c) It was established within the proof of Theorem \ref{te1.11}, that
there is a surjection $\phi$ of $l_{p}^{n}$ onto its arbitrary convex
subset $C$ such that
\begin{equation}\label{open1}
|\phi|_{Lip(l_{p}^{n},C)}\leq n^{\left|\frac{1}{p}-\frac{1}{2}\right|}\ ,\ \ \
1\leq p\leq\infty\ .
\end{equation}
Is it true that the metric projection of $l_{p}^{n}$ onto $C$ satisfies
a similar inequality with a constant $c(p)$ depending only on $p$? The
result is known for $p=2$ with $c(p)=1$, but in general $c(p)$ should be more
than $1$. E.g., the sharp Lipschitz constant for the metric projection in
$l_{\infty}^{2}$ is $2$, not $\sqrt{2}$ (V. Dol'nikov, a personal
communication).\\
(c) Is it true that $\lambda(\Gamma)$ is finite for any uniform lattice
$\Gamma$ of a metric space $M$ (Conjecture \ref{con1})?\\
(d) Is it true that for a complete simply connected length space $M$ of a
non-positive curvature in the A. D. Alexandrov sense $\lambda(M)<\infty$?\\
(e) Assume that $(M_{i},d_{i})$ satisfy $\lambda(M_{i})<\infty$, $i=1,2$.
Is the same true for $M_{1}\times M_{2}$ endowed with the metric
$\max(d_{1},d_{2})$?
\sect{Proof of Theorem \ref{te11}: the Case of $S=M$} {\bf 5.1.}
By a technical reason we prefer to deal with a {\em punctured}
metric space $(M,d,m^{*})$ with a designed point $m^{*}\in M$. A
{\em subspace} $S$ of this space has to contain $m^{*}$ while
$Lip_{0}(S)$ stands for a closed subspace of $Lip(S)$ given by the
condition
\begin{equation}\label{e21}
f(m^{*})=0\ .
\end{equation}
Since $f\mapsto f-f(m^{*})$ is a projection of norm $1$ from
$Lip(S)$ onto $Lip_{0}(S)$, the constant $\lambda(S)$ is unchanged
after replacing $Lip(S)$ by $Lip_{0}(S)$.

In the sequel we will exploit the following fact.\\
{\bf Theorem A}\ {\em There is a Banach space $K(S)$ predual to
$Lip_{0}(S)$ and such that all evaluation functionals
$\delta_{m}:f\mapsto f(m)$, $m\in S$, belong to $K(S)$.}
\ \ \ \ \ $\Box$

$\!\!\!$For compact metric spaces $S$ this can be derived from the
Kantorovich-Rubinshtein theorem [KR] and for separable $S$ from
the result of Dudley [Du] who generalized the
Kantorovich-Rubinshtein construction to this case. For bounded
metric spaces a predual space with the required in Theorem A
property can be found in [We], see also [GK]. For the convenience
of the reader we present below a simple alternative proof of
Theorem A, see Appendix.
\\
{\bf 5.2.} We begin with the next result constituting the first part of
the proof (for $S=M$, i.e., for Corollary \ref{c12}).
\begin{Proposition}\label{p21}
Let $F$ be a given finite point subset of $(M,m^{*})$. Assume that for every
finite $G\supset F$ there is an extension operator $E_{G}\in Ext(F,G)$
from $F$ to $G$ and
\begin{equation}\label{e22}
A:=\sup_{G}||E_{G}||<\infty\ .
\end{equation}
Then there is $E\in Ext(F,M)$ such that
\begin{equation}\label{e23}
||E||\leq A\ .
\end{equation}
\end{Proposition}
{\bf Proof.} Introduce a map $J:M\rightarrow K(M)$ by
\begin{equation}\label{e24}
J(m):=\delta_{m}\ .
\end{equation}
\begin{Lm}\label{l22}
\begin{itemize}
\item[{\rm (a)}]
$J$ is an isometric embedding;
\item[{\rm (b)}]
$K(M)$ is the minimal closed subspace containing $J(M)$.
\end{itemize}
\end{Lm}
{\bf Proof.} Let $f\in K(M)^{*}=Lip_{0}(M)$. Then for $m',m''\in M$
$$
|<f,\ \delta_{m'}-\delta_{m''}>|:=|f(m')-f(m'')|\leq ||f||_{Lip(M)}
d(m',m'')\ .
$$
Taking supremum over $f$ from the unit ball of $Lip_{0}(M)$ to have
$$
||J(m')-J(m'')||_{K(M)}\leq d(m',m'')\ .
$$
To prove the converse, one defines a function $g\in Lip_{0}(M)$ by
$$
g(m):=d(m,m')-d(m^{*},m')\ .
$$
Then
$$
|g(m')-g(m'')|=d(m',m'')\ \ \ {\rm and}\ \ \
|g(m_{1})-g(m_{2})|\leq d(m_{1},m_{2})\ ,\ \ \ m_{1},m_{2}\in M\ .
$$
In particular,
$$
||g||_{Lip(M)}=1\ \ \ {\rm and}\ \ \
||J(m')-J(m'')||_{K(M)}\geq |g(m')-g(m'')|=d(m',m'')\ ,
$$
and the first assertion is done.

Let now $X:=\overline{lin J(M)}\neq K(M)$. Then there is a nonzero linear
functional $f\in K(M)^{*} (=Lip_{0}(M))$ which is zero on $X$. By the
definition of $X$
$$
f(m)=<f,\ \delta_{m}-\delta_{m^{*}}>=<f,\ J(m)>=0\
$$
for every $m\in M$ and so $f=0$, a contradiction.\ \ \ \ \ $\Box$

Let now
\begin{equation}\label{e25}
\kappa_{S}:K(S)\rightarrow K(S)^{**}
\end{equation}
be the canonical isometric embedding. Since $K(F)$ is finite-dimensional,
$\kappa_{F}$ is an isomorphism onto $K(F)^{**}$.

For a finite $G\supset F$ one introduces a linear operator
$\rho_{G}:K(G)\rightarrow K(F)^{**}$ by
\begin{equation}\label{e26}
\rho_{G}:=E_{G}^{*}\kappa_{G}
\end{equation}
where the star designates a conjugate operator. Then define a
vector-valued function $\phi_{G}:G\rightarrow K(F)^{**}$ by
\begin{equation}\label{e27}
\phi_{G}(m):=\rho_{G}J(m)\ ,\ \ m\in G\ .
\end{equation}
\begin{Lm}\label{l23}
\begin{itemize}
\item[{\rm (a)}]
$\phi_{G}\in Lip(G,K(F)^{**})$ and its norm satisfies
\begin{equation}\label{e28}
||\phi_{G}||\leq A\ .
\end{equation}
\item[{\rm (b)}]
For $m\in F$
\begin{equation}\label{e29}
\phi_{G}(m)=\kappa_{G}J(m)
\end{equation}
In particular, $\phi_{G}(m^{*})=0$.
\end{itemize}
\end{Lm}
{\bf Proof.} (a) Let $m\in G$ and $h\in Lip_{0}(F) (=K(F)^{*})$. By
(\ref{e26}) and (\ref{e27})
$$
<\phi_{G}(m),\ h>=<\kappa_{G}J(m),\ E_{G}h>=<E_{G}h,\ J(m)>=(E_{G}h)(m)\ .
$$
This immediately implies that
$$
|<\phi_{G}(m')-\phi_{G}(m''),\ h>|\leq ||E_{G}||
d(m',m'')||h||_{Lip_{0}(F)}\ ,\ \ \ m',m''\in G\ .
$$
This, in turn, gives (\ref{e28}).\\
(b) Since $(E_{G}h)(m)=h(m)$, $m\in F$, and $h(m^{*})=0$, the previous
identity implies (\ref{e29}).\ \ \ \ \ $\Box$

The family of the functions $\{\phi_{G}\}$ is indexed by the
elements $G\supset F$ forming a {\em net}. We now introduce a
topology on the set $\Phi$ of functions $\psi:M\rightarrow
K(F)^{**}$ satisfying the inequality
\begin{equation}\label{e210}
||\psi(m)||_{K(F)^{**}}\leq Ad(m,m^{*})\ ,\ m\in M\ ,
\end{equation}
that allows to find a limit point of the family $\{\phi_{G}\}$. Denote by
$B_{m}$, $m\in M$, the closed ball in $K(F)^{**}$ centered at $0$ and of
radius $Ad(m,m^{*})$. Then (\ref{e210}) means that for every $\psi\in\Phi$
$$
\psi(m)\in B_{m}\ ,\ m\in M\ .
$$
Let $Y:=\prod_{m\in M}B_{m}$ equipped with the product topology. Since
$K(F)^{**}$ is finite-dimensional, $B_{m}$ is compact and therefore
$Y$ is compact, as well. Let $\tau:\Phi\rightarrow Y$ be the natural
bijection given by
$$
\psi\mapsto (\psi(m))_{m\in M}\ .
$$
Identifying $\Phi$ with $Y$ one equips $\Phi$ with the topology of $Y$. Then
$\Phi$ is compact.

Let now $\widehat\phi_{G}:M\rightarrow K(F)^{**}$ be the extension
of $\phi_{G}$ from $G$ by zero. By Lemma \ref{l23} $\widehat\phi_{G}$ meets
condition (\ref{e210}), i.e., $\{\widehat\phi_{G}\}\subset\Phi$.
By compactness of $\Phi$ there is a subnet $N=\{\widehat\phi_{G_{\alpha}}\}$
of the net $\{\widehat\phi_{G}\ : G\supset F\}$ such that
\begin{equation}\label{e211}
\lim\widehat\phi_{G_{\alpha}}=\phi
\end{equation}
for some $\phi\in\Phi$, see e.g. [Ke, Chapter 5, Theorem 2].
By the definition of
the product topology one also has
\begin{equation}\label{e212}
\lim\widehat\phi_{G_{\alpha}}(m)=\phi(m)\ ,\ m\in M\ ,
\end{equation}
(convergence in $K(F)^{**}$).

Show that $\phi$ is Lipschitz. Let $m',m''\in M$ be given, and
$\widetilde N$ be a subnet of $N$ containing those of
$\widehat\phi_{G_{\alpha}}$ for which $m',m''\subset G_{\alpha}$.
Then by (\ref{e26}), (\ref{e27}) and (\ref{e212}) one has for
$h\in Lip_{0}(F)$
$$
\begin{array}{c}
\displaystyle
<\phi(m')-\phi(m''),\ h>=\lim_{\widetilde N}<\phi_{G_{\alpha}}(m')-
\phi_{G_{\alpha}}(m''),\ h>=
\\
\\
\displaystyle
\lim_{\widetilde N}<E_{G_{\alpha}}h,\ J(m')-J(m'')>=
\lim_{\widetilde N}[(E_{G_{\alpha}}h)(m')-(E_{G_{\alpha}}h)(m'')]\ .
\end{array}
$$
Together with (\ref{e22}) this leads to the inequality
$$
|<\phi(m')-\phi(m''),\ h>|\leq A||h||_{Lip_{0}(F)}d(m',m'')\ ,
$$
that is to say,
\begin{equation}\label{e213}
||\phi||_{Lip(M,K(F)^{**})}\leq A\ .
\end{equation}
Using (\ref{e29}) we also similarly prove that for $m\in F$
\begin{equation}\label{e214}
\phi(m)=\kappa_{F}J(m)\ ,\ \ \ {\rm and}\ \ \ \phi(m^{*})=0\ .
\end{equation}

Utilizing the function $\phi$ we, at last, define the required extension
operator \penalty-10000 $E:Lip_{0}(F)\rightarrow Lip_{0}(M)$ as follows. Let
\begin{equation}\label{e215}
\widetilde \kappa_{F}:K(F)^{*}\rightarrow K(F)^{***}
\end{equation}
be the canonical embedding (an isometry in this case). For
$h\in Lip_{0}(F)=K(F)^{*}$ we define $Eh$ by
\begin{equation}\label{e216}
(Eh)(m):=<\widetilde \kappa_{F}h,\ \phi(m)>\ ,\ \ m\in M\ .
\end{equation}
Then by (\ref{e213})
$$
|(Eh)(m')-(Eh)(m'')|=|<\phi(m')-\phi(m''),\ h>|\leq
A||h||_{Lip_{0}(F)}d(m',m'')\ ,
$$
and (\ref{e23}) is proved.

Now by (\ref{e214}) we have for $m\in F$
$$
(Eh)(m)=<\phi(m),\ h>=<\kappa_{F}J(m),\ h>=h(m)\ ;
$$
in particular, $(Eh)(m^{*})=0$.

The proof of the proposition is done.\ \ \ \ \ $\Box$

We are now ready to prove Theorem \ref{te11} for the case of
$S=M$.

Since the inequality
\begin{equation}\label{e217}
A:=\sup_{F}\lambda(F)\leq\lambda(M)
\end{equation}
with $F$ running through finite subspaces of $(M,m^{*})$ is trivial,
we have to establish the converse. In other words, we have to prove
that for every $S\subset (M,m^{*})$ and $\epsilon>0$ there is
$E\in Ext(S,M)$ with
\begin{equation}\label{e218}
||E||\leq A+\epsilon\ .
\end{equation}
By the definition of $A$, for each pair $F\subset G$ of finite subspaces
of $S$ there is $E_{G}\in Ext(F,G)$ with $||E_{G}||\leq A+\epsilon$.
Applying to family $\{E_{G}\}$ Proposition \ref{p21} we find for every
finite $F$ an operator $\widetilde E_{F}\in Ext(F,M)$ with
\begin{equation}\label{e219}
||\widetilde E_{F}||\leq A+\epsilon\ .
\end{equation}
To proceed with the proof we need the following fact.
\begin{Lm}\label{l24}
There is a linear isometric embedding $I_{S}:K(S)\rightarrow\kappa_{M}(K(M))$.
\end{Lm}
{\bf Proof.}
Let $R_{S}:f\mapsto f|_{S}$ be the restriction to $S\subset (M,m^{*})$.
Clearly, $R_{S}$ is a linear mapping from $Lip_{0}(M)$ onto
$Lip_{0}(S)$ with norm bounded by $1$. Since each $f\in Lip_{0}(S)$
has a preserving norm extension to $Lip_{0}(M)$, see e.g. [Mc],
\begin{equation}\label{e220}
||R_{S}||=1\ .
\end{equation}
Introduce now the required linear operator $I_{S}$ by
\begin{equation}\label{e221}
I_{S}:=R_{S}^{*}\kappa_{S}\ .
\end{equation}
Since $R_{S}^{*}$ maps the space $Lip_{0}(S)^{*}=K(S)^{**}$ into
$Lip_{0}(M)^{*}=K(M)^{**}$ and $\kappa_{S}:K(S)\rightarrow
K(S)^{**}$ is the canonical embedding, the operator $I_{S}$ maps
$K(S)$ into $K(M)^{**}$. Show that it, in fact, sends $K(S)$ into
$\kappa_{M}(K(M))$. Let $m\in S$ and $g\in Lip_{0}(M)$. Then by
(\ref{e221}) and the definition of $R_{S}$ we get
$$
<I_{S}J(m),\ g>=<\kappa_{S}J(m),\ g|_{S}>=<g,\ J(m)>:=g(m)=
<\kappa_{M}J(m),\ g>\
$$
which implies the embedding
$$
I_{S}(J(S))\subset \kappa_{M}(J(M))\ .
$$
According to Lemma \ref{l22} this, in turn, implies the required embedding
of $I_{S}(K(S))$ into $\kappa_{M}(K(M))$. At last, by (\ref{e220}) and
(\ref{e221})
$$
||I_{S}||=||R_{S}^{*}||=||R_{S}||=1\ ,
$$
and the result is done.\ \ \ \ \ \ $\Box$

Let now $\widetilde E_{F}$ be the operator from (\ref{e219}). Let us
define an operator \penalty-10000 $P_{F}: K(M)\rightarrow K(F)^{**}$ by
\begin{equation}\label{e222}
P_{F}:=\widetilde E_{F}^{*}\kappa_{M}\ .
\end{equation}
Using the isometric embedding $I_{F}:K(F)\rightarrow \kappa_{S}(K(S))$ of
Lemma \ref{e24} we then introduce an operator
$Q_{F}:K(M)\rightarrow \kappa_{S}(K(S))\subset K(S)^{**}$
by
\begin{equation}\label{e223}
Q_{F}:=I_{F}(\kappa_{F})^{-1}P_{F}\ .
\end{equation}
Since $dim\ \! K(F)<\infty$, this is well-defined. Introduce, at last,
a vector-valued function $\phi_{F}:M\rightarrow K(S)^{**}$ by
\begin{equation}\label{e224}
\phi_{F}(m):=Q_{F}J(m)\ ,\ \ m\in M\ .
\end{equation}
Arguing as in Lemma \ref{l23} and using the estimate
$||Q_{F}||\leq ||P_{F}||\leq ||\widetilde E_{F}||\leq A+\epsilon$, see
(\ref{e222}), (\ref{e223}) and (\ref{e219}), we obtain the inequality
\begin{equation}\label{e225}
||\phi_{F}||_{Lip(M,K(S)^{**})}\leq A+\epsilon\ .
\end{equation}
Moreover, for each $m\in F$ and $h\in Lip_{0}(M)$
$$
\begin{array}{c}
<\phi_{F}(m),\ h>=<R_{F}^{*}\kappa_{F}(\kappa_{F})^{-1}P_{F}J(m),\ h>=
<h|_{F},\ P_{F}J(m)>=\\
\\
<E_{F}h,\ J(m)>=h(m)\ ,
\end{array}
$$
see (\ref{e221}) and (\ref{e223})-(\ref{e225}). Hence for $m\in F$
\begin{equation}\label{e226}
\phi_{F}(m)=\kappa_{M}J(m)\ ;
\end{equation}
in particular, $\phi_{F}(m^{*})=0$.

From here and (\ref{e225}) we derive that the set $\{\phi_{F}(m)\}$ with
$F$ running through the net of all finite subspaces of $S$ is a subset of the
closed ball $B_{m}\subset K(S)^{**}$ centered at 0 and of radius
$(A+\epsilon)d(m,m^{*})$. In the weak$^{*}$ topology $B_{m}$ is compact.
From this point our proof repeats word for word that of Proposition \ref{p21}.
Namely, consider the set $\Phi$ of functions $\psi:M\rightarrow K(S)^{**}$
satisfying
$$
||\psi(m)||_{K(S)^{**}}\leq (A+\epsilon)d(m,m^{*})\ ,\ \ m\in M\ .
$$
Equip $B_{m}$ with the weak$^{*}$ topology and introduce the set
$Y:=\prod_{m\in M}B_{m}$ equipped with the product topology. Then $Y$ is
compact and, so, $\Phi$ is too in the topology induced by the bijection
$\Phi\ni\psi\mapsto (\psi(m))_{m\in M}\in Y$. Then there is a subnet $N$ of
the net $\{\phi_{F}\ :\ (F,m^{*})\subset (S,m^{*})\ ,\ \# F<\infty\}$ such
that
$$
\lim_{N}\phi_{F}=\phi
$$
for some $\phi\in\Phi$. By the definition of the product topology
$$
\lim_{N}\phi_{F}(m)=\phi(m)\ ,\ \ m\in M\
$$
(convergence in the weak$^{*}$ topology of $K(S)^{**}$). Arguing as in the
proof of Proposition \ref{p21}, see (\ref{e213}), we
derive from (\ref{e225}) that
\begin{equation}\label{e227}
||\phi||_{Lip(M,K(S)^{**})}\leq A+\epsilon
\end{equation}
and, moreover, for $m\in S$
\begin{equation}\label{e228}
\phi(m)=\lim_{N'}\phi_{F}(m)=\kappa_{M}J(m)\ .
\end{equation}
Here $N':=\{\phi_{F}\in N\ :\ m\in F\subset S\}$ is a subnet of $N$,
and we have used (\ref{e226}).

Using now the canonical embedding
$\widetilde\kappa_{S}:K(S)^{*}=Lip_{0}(S)\rightarrow K(S)^{***}$ we introduce
the required extension operator $E\in Ext(S,M)$ by
$$
(Eh)(m):=<\widetilde\kappa_{S}h,\ \phi(m)>\ ,\ \ m\in M\ ,\ h\in Lip_{0}(S)\ .
$$
Since $\phi(m)\in K(S)^{**}$, this is well-defined. Then for
$m',m''\in M$ we get from (\ref{e227})
$$
\begin{array}{c}
|(Eh)(m')-(Eh)(m'')|\leq ||h||_{Lip_{0}(S)}||\phi(m')-\phi(m'')||_{K(S)^{**}}
\leq \\
\\
(A+\epsilon)||h||_{Lip_{0}(S)}d(m',m'')\ .
\end{array}
$$
Moreover, by (\ref{e228}) we have for $m\in S$
$$
(Eh)(m)=<\widetilde\kappa_{S}h,\ \kappa_{M}J(m)>=<h,\ J(m)>=h(m)\ .
$$
Hence $E\in Ext(S,M)$ and $||E||\leq A+\epsilon$. This implies the
converse to (\ref{e217}) inequality
$$
\lambda(M)\leq\sup_{F}\lambda(F) (=A)\ .
$$
Proof of Theorem \ref{te11} for $S=M$ is complete.\ \ \ \ \ \ $\Box$
\sect{Proof of Theorem \ref{te11}: the Final Part}
We begin with a slight modification of Proposition \ref{p21} using an
increasing sequence of subspaces $\{S_{j}\}_{j\geq 0}$ of the space
$(M,m^{*})$ instead of the net $\{G\}$ of its finite point subspaces
containing a given $F$. Repeating line-to-line the proof of Proposition
\ref{p21} for
this setting we obtain as a result an extension operator
$E\in Ext(F,S_{\infty})$ where $S_{\infty}:=\cup_{j}S_{j}$ with the
corresponding bound for $||E||$. If, in addition, $S_{\infty}$ is dense in
$M$, then there is a canonical isometry $Lip_{0}(S)\leftrightarrow Lip_{0}(M)$
generated by continuous extensions of functions from $Lip_{0}(S_{\infty})$.
This leads to the following assertion.
\begin{Proposition}\label{p31}
Assume that $\{S_{j}\}_{j\geq 0}$ is an increasing sequence of subspaces
of $(M,m^{*})$ whose union $S_{\infty}$ is dense in $M$, and $F$ is a finite
point subspace in $\cap_{j\geq 0}S_{j}$. Suppose that
for every $j$ there exists $E_{j}\in Ext(F,S_{j})$ and
\begin{equation}\label{e31}
A:=\sup_{j}||E_{j}||<\infty\ .
\end{equation}
Then there is an operator $E\in Ext(F,M)$ with $||E||\leq A$.\ \ \ \ \
$\Box$
\end{Proposition}

Let now $\{m_{k}^{j}\}_{j\in\N}$ be a sequence in $M$ convergent to $m_{k}$,
$1\leq k\leq n$. Set $F:=\{m^{*},m_{1},\dots, m_{n}\}$ and
$F_{j}:=\{m^{*},m_{1}^{j},\dots, m_{n}^{j}\}$.
\begin{Proposition}\label{p32}
Assume that for each $j$ there is $E_{j}\in Ext(F_{j},M)$ and
\begin{equation}\label{e32}
A:=\sup||E_{j}||<\infty\ .
\end{equation}
Then there exists $E\in Ext(F,M)$ with
\begin{equation}\label{e33'}
||E||\leq A\ .
\end{equation}
\end{Proposition}
{\bf Proof.}
Let a linear operator $L_{j}:Lip_{0}(F)\rightarrow Lip_{0}(F_{j})$ be given
by
$$
(L_{j}f)(m_{l}^{j}):=f(m_{l})\ ,\ \ 1\leq l\leq n\ , \ \ \ {\rm and}\ \ \
(L_{j}f)(m^{*})=0\ .
$$
Then its norm is bounded by
$\sup_{l'\neq l''}\{d(m_{l'},m_{l''})/d(m_{l'}^{j},m_{l''}^{j})\}$ and
therefore
\begin{equation}\label{e33}
\limsup_{j\to\infty}| |L_{j}||=1\ \ \ {\rm and}\ \ \ ||L_{j}||\leq
2\ ,\ \ j\geq j_{0}\ ,
\end{equation}
for some $j_{0}$.

\noindent Introduce now a linear operator
$\widetilde E_{j}:Lip_{0}(F)\rightarrow Lip_{0}(M)$ by
$$
\widetilde E_{j}:=E_{j}L_{j}\ .
$$
Applying (\ref{e33}) to have
\begin{equation}\label{e34}
\limsup_{j\to\infty}||\widetilde E_{j}||\leq A\ \ \ {\rm and}\ \ \
||\widetilde E_{j}||\leq 2A\ ,\ \ j\geq j_{0}\ .
\end{equation}
Let $R_{j}:K(M)\rightarrow K(F)^{**}$ be given by
$$
R_{j}:=\widetilde E_{j}^{*}\kappa_{M}\ ,
$$
and a vector-valued function $G_{j}:M\rightarrow K(F)^{**}$ be defined by
$$
G_{j}(m):=R_{j}(J(m))\ ,\ \ m\in M\ .
$$
Arguing as in Lemma \ref{l23} we then establish that $G_{j}$ is Lipschitz and
\begin{equation}\label{e35}
||G_{j}||_{Lip(M,K(F)^{**})}\leq 2A\ ,\ \ \ j\geq j_{0}\ ,
\end{equation}
and, moreover,
\begin{equation}\label{e36}
G_{j}(m)=\kappa_{M}(J(m))\ ,\ \ m\in F\ .
\end{equation}
In particular, $G_{j}(m^{*})=0$. These allow
to assert that there is a subnet $\{G_{j_{k}}\}$ of the net $\{G_{j}\}$ such
that $G(m):=\lim_{k}G_{j_{k}}(m)$ exists for each
$m\in M$ in the topology of $K(F)^{**}$
(see the argument of Proposition \ref{p21} after Lemma \ref{l23}).
It remains to define $(Eh)(m)$ for $h\in Lip_{0}(F)$ and $m\in M$ by
$$
(Eh)(m):=<\widetilde\kappa_{F}h,\ G(m)>\ ,
$$
cf. (\ref{e216}). The argument at the end of the proof of Proposition
\ref{p21} with $\widetilde E_{j}$ instead of $E_{G}$ and (\ref{e34})
instead of (\ref{e22}) leads to existence of $E\in Ext(F,M)$ with
$$
||E||\leq\limsup_{j\to\infty}||\widetilde E_{j}||\leq A\ .\ \ \ \ \ \Box
$$

Now we will finalize the proof of Theorem \ref{te11}. Since the inequality
\begin{equation}\label{e37}
\sup_{F\in S}\lambda(F)=\lambda(S)\leq\lambda(M)
\end{equation}
is clear, we have to prove the following assertion.

Given a finite $F\subset (M,m^{*})$ and $\epsilon>0$, there is
$E\in Ext(F,M)$ such that
\begin{equation}\label{e38}
||E||\leq\lambda(S)+\epsilon\ .
\end{equation}
Since $\lambda(M)=\sup_{F}\lambda(F)$, this will prove the converse to
(\ref{e37}).

To establish this assertion one considers first the case of $F$ containing
in \penalty-10000
$S_{\infty}:=\cup_{j\geq 0}S_{j}$ where $S_{j}:=\delta^{j}(S)$,
$j=0,1,\dots$. Since $S_{j}$ is increasing, see assumption (a) of Theorem
\ref{te11}, there is $j=j(F)$ such that $F\subset S_{j}$ and so
$$
F_{j}:=\delta^{-j}(F)\subset S\ .
$$
Then there is $E_{j}\in Ext(F_{j},S)$ such that
$$
||E_{j}||\leq\lambda(S)+\epsilon\ .
$$
Introduce now linear operators $D_{j}:Lip(F)\rightarrow Lip(F_{j})$
and $H_{j}:Lip(S)\rightarrow Lip(S_{j})$ by
$$
\begin{array}{c}
(D_{j}f)(m):=f(\delta^{j}(m))\ ,\ \ \ f\in Lip(F)\ ,\ \ m\in F_{j}\ ,
\\
\\
(H_{j}f)(m):=f(\delta^{-j}(m))\ ,\ \ \ f\in Lip(S)\ ,\ \ m\in S_{j}\ .
\end{array}
$$
Then the operator
\begin{equation}\label{e39}
\widetilde E_{j}:=H_{j}E_{j}D_{j}
\end{equation}
clearly belongs to $Ext(F,S_{j})$. Then we have
$$
D_{j}R_{F}=R_{F_{j}}\Delta^{j}\ \ \ {\rm and} \ \ \
H_{j}R_{S}=R_{S_{j}}\Delta^{-j}\ ,
$$
where $R_{K}:f\mapsto f|_{K}$
stands for the restriction operator to $K\subset M$, and
$\Delta$ is given by  (\ref{dil}).\\
Since each $f\in Lip_{0}(K)$ has a preserving norm extension to
$Lip_{0}(M)$, the above identities imply that
$$
\begin{array}{c}
||D_{j}||=||D_{j}R_{F}||=
||R_{F_{j}}\Delta^{j}||\leq ||\Delta||^{j}\ \ \ {\rm and}\\
\\
||H_{j}||=||H_{j}R_{S}||=||R_{S_{j}}\Delta^{-j}||\leq
||\Delta^{-1}||^{j}\ .
\end{array}
$$
This, in turn, leads to the estimate
\begin{equation}\label{e310}
||\widetilde E_{j}||\leq ||H_{j}||\cdot ||E_{j}||\cdot ||D_{j}||\leq
||E_{j}||(||\Delta||\cdot ||\Delta^{-1}||)^{j}\leq\lambda(S)+\epsilon\ .
\end{equation}
So we have constructed a sequence $\{\widetilde E_{j}\}$ of operators whose
norms are bounded by (\ref{e310}). Moreover, $\{S_{j}\}$ is increasing and
its union is dense in $M$, see assumption (b) of Theorem \ref{te11}.
Hence we are under the conditions of Proposition \ref{p31} that guarantees
existence of an operator $E\in Ext(F,M)$ with $||E||\leq\lambda(S)+\epsilon$.
Therefore the inequality (\ref{e38}) is done for such $F$.

Let now $F:=\{m^{*},m_{1},\dots, m_{n}\}$ be arbitrary. Since
$S_{\infty}$ is dense in $M$, one can find a sequence
$F_{j}:=\{m^{*},m_{1}^{j},\dots, m_{n}^{j}\}$, $j\in\N$, of subsets of
$S_{\infty}$ such that
$$
\lim_{j\to\infty}m_{l}^{j}=m_{l}\ ,\ \ \ 1\leq l\leq n\ .
$$
It had just proved that for each $j$ there is $E_{j}\in Ext(F_{j},M)$
such that
\begin{equation}\label{e311}
||E_{j}||\leq\lambda(S)+\epsilon\ ,\ \ j\in\N\ .
\end{equation}
Applying to this setting Proposition \ref{p32} we conclude that there is
$E\in Ext(F,M)$ satisfying (\ref{e38}). Together with (\ref{e37}) this
completes the proof of Theorem \ref{te11}.\ \ \ \ \ $\Box$
\begin{R}\label{rnew}
{\rm (a) Using the compactness argument of the proof of Proposition
\ref{p21} one can also prove that for every $S\subset M$ there is an
extension operator \penalty-10000 $E_{min}\in Ext(S,M)$ such that}
\begin{equation}\label{enew}
||E_{min}||=\inf\{||E||\ :\ E\in Ext(S,M)\}
\end{equation}
{\rm (b) The same argument combining with some additional consideration
allow to establish the following fact.

The set function $S\mapsto\lambda(S)$ defined on closed subsets of $M$ is
lower semicontinuous in the Hausdorff metric.}
\end{R}
\sect{Proof of Theorem \ref{te14a}} {\bf Trees with all edges of
length one.} In accordance with Corollary \ref{c12} one has to
find for every pair $S\subset S'$ of finite sets in
$M_{p}:=\oplus_{p}\{{\cal T}_{i}\}_{1\leq i\leq n}$, 
an operator $E\in Ext(S,S')$
whose norm is bounded by
\begin{equation}\label{eq71}
||E||\leq cn
\end{equation}
with $c$ independent of $S$ and $S'$; recall that $p=1$ or $\infty$.

To accomplish this we first find a subset
$S''=\prod_{i=1}^{n}S_{i}$, $S_{i}\subset {\cal T}_{i}$, such that
$S'\subset S''$. Further, every finite subset of a tree can be
isometrically embedded into an infinite rooted tree ${\cal T}_{k}$
with vertices of degree $k+2$ for some $k\in\N$ (and all edges of
length one in our case). Note that the {\em degree of vertex} $v$,
written $deg\ \! v$, {\em is the number of its children plus 1}.
Taking some $k$ such that every set $S_{i}$ is an isometric
part of ${\cal T}_{k}$ we therefore can derive 
(\ref{eq71}) from a similar inequality with $S'$ substituted for
$\oplus_{p}\{{\cal T}_{k}\}_{1\leq i\leq n}$ (= the direct $p$-sum of
$n$ copies of ${\cal T}_{k}$) . Hence we have reduced the
required result to the following assertion:
\begin{equation}\label{eq72}
\lambda(\oplus_{p}\{{\cal T}_{k}\}_{1\leq i\leq n})\leq cn
\end{equation}
with $c$ independent of $n$ and $k$.

The proof of this inequality is divided into two parts. We first
prove that each ${\cal T}_{k}$ can be quasi-isometrically embedded
into the hyperbolic plane $\H^{2}$, see Proposition \ref{p161}
below. From here we derive that
\begin{equation}\label{164}
\lambda(\oplus_{p}\{{\cal T}_{k}\}_{1\leq i\leq n})\leq
256\ \!\lambda(\oplus_{p}\{\H ^{2}\}_{1\leq i\leq n}).
\end{equation}
Then we estimate the right-hand side applying inequalities (\ref{e2.13'}) and
(\ref{eq212}) of Theorem \ref{c118} for the case $M_{i}=\H^{2}$
for all $1\leq i\leq n$. Combining these we prove the right-hand side
inequality of Theorem
\ref{te14a} for trees with edges of length one.

We begin with establishing the desired quasi-isometric embedding
of ${\cal T}_{k}$ into $\H^{2}$. In the formulation of this
result, $d$ and $\rho$ are, respectively, the path metrics on
${\cal T}_{k}$ and $\H ^{2}=\{x\in\Re^{2}\ :\ x_{2}>0\}$.
\begin{Proposition}\label{p161}
For every $k\geq 2$ there is an embedding $I:{\cal T}_{k}\to \H ^{2}$ such
that for all $m_{1},m_{2}\in {\cal T}_{k}$
\begin{equation}\label{164a}
A\ \!d(m_{1},m_{2})\leq\rho(I(m_{1}),I(m_{2}))\leq B\ \!d(m_{1},m_{2})
\end{equation}
with constants $0<A<B$ independent of $m_{1},m_{2}$ and satisfying
\begin{equation}\label{165}
BA^{-1}\leq 256\ .
\end{equation}
\end{Proposition}
{\bf Proof.}
It will be more appropriate to work with another metric on $\H ^{2}$ given
by
\begin{equation}\label{166}
\rho_{0}(x,y):=\max_{i=1,2}\ \!
\log\left(1+\frac{|x_{i}-y_{i}|}{\min(x_{2},y_{2})}
\right)\ .
\end{equation}
The following result establishes an
equivalence of this to the hyperbolic metric $\rho$ for pairs of
points far enough from each other.
\begin{Lm}\label{l162}
\begin{itemize}
\item[(a)] $\rho\leq 4\rho_{0}$;
\item[(b)] If $|x-y|\geq\frac{1}{2}\min(x_{2},y_{2})$, then
$$
\rho(x,y)\geq\frac{1}{8}\rho_{0}(x,y)\ .
$$
\end{itemize}
Here and below $|x|$ is the Euclidean norm of $x\in\Re_{+}^{2}$.
\end{Lm}
{\bf Proof.}
For definiteness assume that
\begin{equation}\label{167}
\min(x_{2},y_{2})=y_{2}\ .
\end{equation}
Use for a while the complex form of $\H ^{2}$ with the underlying set
$\{z\in\Co\ :\ Im\ \! z>0\}$. Then the metric $\rho$ is given by
$$
\rho(z_{1},z_{2})=\log\frac{1+\left|
\frac{z_{1}-z_{2}}{z_{1}-\overline{z}_{2}}\right|}{1-\left|
\frac{z_{1}-z_{2}}{z_{1}-\overline{z}_{2}}\right|}\ .
$$
Identifying $z=x_{1}+ix_{2}$, with $(x_{1},x_{2})\in \H ^{2}$ we rewrite this
as
\begin{equation}\label{168}
\rho(x,y)=\log\ \!\frac{(|x-y^{+}|+|x-y|)^{2}}{|x-y^{+}|^{2}-|x-y|^{2}}
\end{equation}
where $y^{+}:=(y_{1},-y_{2})$ is the reflexion of $y$ in the $x_{1}$-axis.
Since the denominator in (\ref{168}) equals $4x_{2}y_{2}$ and
$$
|x-y^{+}|+|x-y|\leq 2|x-y|+|y^{+}-y|=2(|x-y|+y_{2})
$$
we derive from (\ref{168}) and (\ref{167}) the inequality
$$
\rho(x,y)\leq 2\log\frac{y_{2}+|x-y|}{y_{2}}\leq
4\log\left(1+\frac{\max_{i=1,2}|x_{i}-y_{i}|}{y_{2}}\right)\ .
$$
By (\ref{166}) and (\ref{167}) this implies the required result formulated
in (a).

In case (b) we use an equivalent formula for $\rho$ given by
$$
\cosh\rho(x,y)=1+\frac{|x-y|^{2}}{2x_{2}y_{2}}\ .
$$
Since $\cosh t\leq e^{t}$ for $t\geq 0$, this yields
\begin{equation}\label{169}
\rho(x,y)\geq\log\left(1+\frac{|x-y|^{2}}{2x_{2}y_{2}}\right)\ .
\end{equation}
Consider two possible cases
\begin{equation}\label{1610}
y_{2}\leq x_{2}\leq 2y_{2}\ ;
\end{equation}
\begin{equation}\label{1611}
2y_{2}<x_{2}\ .
\end{equation}
In the first case we use (\ref{167}) and the assumption
$\frac{2|x-y|}{y_{2}}\geq 1$ to derive from (\ref{169})
$$
\begin{array}{c}
\displaystyle
\rho(x,y)\geq\log\left(1+\left(\frac{|x-y|}{2y_{2}}\right)^{2}\right)\geq
\frac{1}{8}\log\left(1+\frac{|x-y|}{y_{2}}\right)\geq\\
\\
\displaystyle
\frac{1}{8}
\log\left(1+\frac{\max_{i=1,2}|x_{i}-y_{i}|}{y_{2}}\right):=\frac{1}{8}\ \!
\rho_{0}(x,y)\ .
\end{array}
$$

In the second case we have from (\ref{1611})
$|x-y|\geq |x_{2}-y_{2}|\geq\frac{1}{2}x_{2}$. Inserting this in (\ref{169})
we obtain the required result
$$
\rho(x,y)\geq\log\left(1+\frac{|x-y|}{4y_{2}}\right)\geq
\frac{1}{4}\log\left(1+\frac{\max_{i=1,2}|x_{i}-y_{i}|}{y_{2}}\right):=
\frac{1}{4}\ \!\rho_{0}(x,y)\ .\ \ \ \ \ \Box
$$

We now begin to construct the required embedding $I:{\cal T}_{k}\to \H ^{2}$.
To this end we first introduce coordinates for the set of vertices
${\cal V}_{k}$ of the ${\cal T}_{k}=(R,{\cal V}_{k},{\cal E}_{k})$ where
$R$ stands for the
root. Actually, we assign to $v\in {\cal V}_{k}$ the pair of integers
$(j_{v},l_{v})$ determined as follows. The number $l_{v}$ is the {\em level}
of $v$, the length of the unique path from the root to $v$. To define
$j_{v}$ we visualize ${\cal T}_{k}$ using the natural isometric embedding
of ${\cal T}_{k}$ into $\Re^{2}$. Then $j_{v}$ is the number of $v$ in the
ordering of the vertices of the $l_{v}$-th level from the left to the right.
We use in this ordering the set of numbers $0,1,\dots$ and therefore
$$
0\leq j_{v}< (k+1)^{l_{v}}\ ,
$$
since the number of children of each vertex equals $k+1$. We also assign
$(0,0)$ to the root $R$ of ${\cal T}_{k}$.

Using this we relate the coordinates of $v$ and its parent $v^{+}$. To this
end one uses the $(k+1)$-ary digital system to present $j_{v}$ as
\begin{equation}\label{1612}
j_{v}=\sum_{s=1}^{l_{v}}\delta_{s}(v)(k+1)^{s-1}
\end{equation}
where $\delta_{s}(v)\in\{0,1,\dots,k\}$ are the {\em digits}. Then the
coordinates of $v$ and $v^{+}$ are related by
\begin{equation}\label{1613}
l_{v^{+}}=l_{v}-1\ ;
\end{equation}
\begin{equation}\label{1614}
\delta_{s}(v^{+})=\delta_{s+1}(v)+1\ ,\ \ \ s=1,\dots,l_{v^{+}}\ .
\end{equation}
To express the distance between $v,w\in {\cal V}_{k}$ in their coordinates,
we first introduce the notion of the {\em common ancestor} $a(v,w)$ of these
vertices. This is the vertex of the biggest level in the intersection of
the paths joining the root with $v$ and $w$, respectively. Hence there are
sequences $v:=v_{1}, v_{2},\dots, v_{p}:=a(v,w)$ and $w:=w_{1}, w_{2},
\dots,w_{q}:=a(v,w)$ such that $v_{i+1}=v_{i}^{+}$ (the parent of $v_{i}$),
$w_{i+1}=w_{i}^{+}$ and the children $v_{p-1}$ and $w_{q-1}$ of $a(v,w)$ are
distinct. The distance between $v$ and $w$ in the ${\cal T}_{k}$ is
therefore given by
\begin{equation}\label{1615}
d(v,w)=l_{v}+l_{w}-2l_{a(v,w)}\ .
\end{equation}

We now define the required embedding $I:{\cal T}_{k}\to \H ^{2}$ on the set
of vertices ${\cal V}_{k}\subset {\cal T}_{k}$. To this end we assign to each
$v\in {\cal V}_{k}$ a square $Q(v)$ in $\Re_{+}^{2}$ in the following fashion.
For the root $R$ we define $Q(R)$ to be the square in $\Re_{+}^{2}$ whose
center $c(R)$ and lengthside $\mu(R)$ are given by
\begin{equation}\label{1616}
c(R)=(c_{1}(R),c_{2}(R))=(0,1)\ ,\ \ \ \mu(R)=\frac{2(n-1)}{n+1}
\end{equation}
where here and below
\begin{equation}\label{1617}
n:=k^{2}+1\ .
\end{equation}
Define now $Q(v)$ for a child $v$ of $R$. In this case $l_{v}=1$ and
$j_{v}=\delta_{1}(v)$, see (\ref{1612}), and we introduce $Q(v)$ to be a
square in $\Re_{+}^{2}$ whose center $c(v)$ and lengthside $\mu(v)$ are given
by
$$
c_{1}(v):=\frac{1}{2}\mu(v)(2\delta_{1}(v)k-k^{2}),\ \ \ c_{2}(v):=\frac{1}{n}\ ,
\ \ \ \mu(v)=\frac{2(n-1)}{n+1}\cdot\frac{1}{n}\ .
$$
Note that these squares are introduced by the following geometric
construction. Divide the bottom side of $Q(R)$ into $n$ equal intervals
and construct outside $Q(R)$ $n$ squares with these intervals as their sides.
Number those from the left to the right. Then the squares numbered by
$1,k+1,2k+1,\dots,k\cdot k+1:=n$ form the set $\{Q(v)\ :\ l_{v}=1\}$.
Apply now this construction to each $Q(v)$ with $l_{v}=1$ to obtain
$(k+1)^{2}$ squares corresponding to vertices of the level $l_{v}=2$ and
so on. Straightforward evaluation leads to the following formulas related
the coordinates of $c(v)$, the center of $Q(v)$, and its lengthside
$\mu(v)$ to those for the square $Q(v^{+})$ associated with the parent
$v^{+}$ of $v$
\begin{equation}\label{1618}
\begin{array}{c}
\displaystyle
c_{2}(v)=\frac{1}{n}c_{2}(v^{+})\ ,\ \ \ \mu(v)=\frac{1}{n}\mu(v^{+})\ ,\\
\\
\displaystyle
c_{1}(v)=c_{1}(v^{+})+\frac{1}{2}\mu(v)(2\delta_{1}(v)k-k^{2})\ .
\end{array}
\end{equation}
In particular, for the second coordinate of $c(v)$ and for $\mu(v)$ we get
\begin{equation}\label{1619}
c_{2}(v)=\frac{1}{n^{l_{v}}}\ ,\ \ \ \mu(v)=\frac{2(n-1)}{n+1}
\frac{1}{n^{l_{v}}}\ .
\end{equation}

Compare now the metrics $\rho$ and $\rho_{0}$ on the set of the centers
$c(v)$, $v\in {\cal V}_{k}$. Let $v,w$ be distinct vertices of ${\cal V}_{k}$.
Without loss of generality we assume that
\begin{equation}\label{1620}
l_{v}\geq l_{w}\ ,
\end{equation}
so that
\begin{equation}\label{1621}
\min(c_{2}(v),c_{2}(w))=n^{-l_{v}}\ .
\end{equation}
Moreover, if $l_{v}=l_{w}$ then
$$
|c_{1}(v)-c_{1}(w)|\geq |c_{1}(v)-c_{1}(\widehat v)|
$$
where $\widehat v$ has coordinates satisfying $l_{\widehat v}=l_{v}$ and
$|j_{\widehat v}-j_{v}|=1$. If now $l_{v}>l_{w}$, then
$|c_{2}(v)-c_{2}(w)|\geq |c_{2}(v)-c_{2}(v')|$ provided that $v'$ satisfies
$l_{v}-l_{v'}=1$. By (\ref{1618}) and (\ref{1619}) the right-hand side of the
inequality for $\widehat v$ is equal to $\frac{2(n-1)k}{n^{l_{v}}(n+1)}$,
while the right-hand side of the inequality for $v'$ is equal to
$\frac{n-1}{n^{l_{v}}n}$. Together with (\ref{1621}) these inequalities
yield the estimate
$$
\frac{|c(v)-c(w)|}{\min(c_{2}(v),c_{2}(w))}\geq\frac{n-1}{n}\geq\frac{1}{2}\ .
$$
Hence the assumption of Lemma \ref{l162} holds for $x:=c(v)$ and
$y:=c(w)$ and we have
\begin{equation}\label{1622}
\frac{1}{8}\rho_{0}(c(v),c(w))\leq\rho(c(v),c(w))\leq
4\rho_{0}(c(v),c(w))\ .
\end{equation}

We now introduce the required embedding $I:{\cal T}_{k}\to\H^{2}$ beginning
with its definition on the subset
${\cal V}_{k}\subset {\cal T}_{k}$ of vertices; namely, we let
$$
I(v):= c(v)\ ,\ \ \ v\in {\cal V}_{k}\ .
$$
Show that $I|_{{\cal V}_{k}}$ is a bi-Lipschitz equivalence with
the constants satisfying (\ref{165}). By (\ref{1622}) it suffices to work
with the metric space $(\Re_{+}^{2},\rho_{0})$. So we have to compare
$\rho_{0}(c(v),c(w))$ with the distance $d(v,w)$ in the tree ${\cal T}_{k}$.
\begin{Lm}\label{l163}
The $\rho_{0}$-distance between $c(v)$ and $c(v^{+})$ equals $\log n$.
\end{Lm}
{\bf Proof.}
By (\ref{1619})
$$
\log\left(1+\frac{|c_{2}(v)-c_{2}(v^{+})|}{\min(c_{2}(v),c_{2}(v^{+}))}
\right)=\log(1+n-1)=\log n\ .
$$
On the other hand, the similar expression with $c_{2}$ replaced by
the first coordinates in the numerator equals
$\log\left(1+\frac{n-1}{n+1}|2\delta_{1}(v)k-k^{2}|\right)$,
see (\ref{1618}). Since $0\leq\delta_{1}(v)\leq k$, this is at most $\log n$.
Hence
$$
\rho_{0}(c(v),c(v^{+})):=\max_{i=1,2}\log
\left(1+\frac{|c_{i}(v)-c_{i}(v^{+})|}{\min(c_{2}(v),c_{2}(v^{+}))}\right)
=\log n\ .\ \ \ \ \ \Box
$$

Thus, the length of each edge in $I({\cal V}_{k})$ equals
$\log n$. Using this we prove the first of the desired estimates.\\
Let $a(v,w)$ be the common ancestor of $v$ and $w$ and
$v:=v_{1},v_{2},\dots,v_{p}:=a(v,w)$ and
$w:=w_{1},w_{2},\dots,w_{q}:=a(v,w)$ are the corresponding
connecting paths. So $v_{i+1}=v_{i}^{+}$, $w_{i+1}=w_{i}^{+}$ and
$p:=l_{v}-l_{a(v,w)}+1$, $q=l_{w}-l_{a(v,w)}+1$. By the triangle
inequality, Lemma \ref{l163} and (\ref{1615}) we then have
$$
\begin{array}{c}
\displaystyle
\rho_{0}(c(v),c(w))\leq\sum_{i=1}^{p-1}d(v_{i},v_{i+1})+
\sum_{i=1}^{q-1}d(w_{i},w_{i+1})\leq\\
\\
\displaystyle
(l_{v}-l_{a(v,w)})\log n+
(l_{w}-l_{a(v,w)})\log n=d(v,w)\log n\ .
\end{array}
$$
So we get
\begin{equation}\label{1623}
\rho_{0}(I(v),I(w))\leq\log n\cdot d(v,w)
\end{equation}
and it remains to establish the inverse inequality. To this end we
consider two cases. First, suppose that $w=a(v,w)$. Then by
(\ref{1619}), (\ref{1621}) and (\ref{1615}) we have
\begin{equation}\label{1624a}
\rho_{0}(c(v),c(w))\geq\log\left(1+
\frac{|c_{2}(v)-c_{2}(w)|}{\min(c_{2}(v),c_{2}(w))}\right)=
\log n^{l_{v}-l_{w}}=d(v,w)\log n\ .
\end{equation}
Suppose now that $w\neq a(v,w)$, then we use the inequality
\begin{equation}\label{1624}
\rho_{0}(c(v),c(w))\geq\log\left(1+
\frac{|c_{1}(v)-c_{1}(w)|}{\min(c_{2}(v),c_{2}(w))}\right)\geq l_{v}\log n
+\log|c_{1}(v)-c_{1}(w)|\ ,
\end{equation}
see (\ref{1621}). To estimate the second summand one notes that, by
our geometric construction, the orthogonal projection of $Q(v)$ onto the
bottom side of $Q(v^{+})$ lies inside this side. Applying this
consequently to the vertices of the chains $\{v_{i}\}_{1\leq i\leq p}$ and
$\{w_{i}\}_{1\leq i\leq q}$ joining $v$ and $w$ with $a(v,w)$, see the proof
of (\ref{1623}), and taking into account that $v_{i+1}=v_{i}^{+}$ and
$w_{i+1}=w_{i}^{+}$ we conclude that the orthogonal projections of
$c(v):=c(v_{1})$ and $c(w):=c(w_{1})$ onto the bottom side of
$Q(a(v,w))=Q(v_{p-1}^{+})=Q(w_{q-1}^{+})$ lie, respectively, inside the
top sides of the squares $Q(v_{p-1})$ and $Q(w_{q-1})$ adjoint to
$Q(a(v,w))$. Hence
$$
|c_{1}(v)-c_{1}(w)|\geq dist(Q(v_{p-1}),Q(w_{q-1}))
$$
and this distance is at least $\frac{k\mu(a(v,w))}{n}$ by the definition of
the squares involved. By the equality $n:=k^{2}+1$, $k\geq 2$,
and formula (\ref{1619}) we derive from here that
$$
\log|c_{1}(v)-c_{1}(w)|\geq\log\left(\frac{2(n-1)^{3/2}}{n(n+1)}
n^{-l_{a}}\right)\geq -(l_{a}+1/2)\log n\ ,\ \ a:=a(v,w)\ .
$$
Together with (\ref{1624}) this yields
$$
\rho_{0}(c(v),c(w))\geq (l_{v}-l_{a}-1/2)\log n\ .
$$
Now note that (\ref{1620}) and the inequality $l_{v}-l_{a}\geq 1$ imply that
$$
l_{v}-l_{a}-1/2\geq\frac{1}{8}(l_{v}+l_{w}-2l_{a})=\frac{1}{8}d(v,w)\ ,
$$
see (\ref{1615}). Hence
$$
\rho_{0}(c(v),c(w))\geq\frac{1}{8}\log n\ \! d(v,w)\ .
$$
Together with (\ref{1623}) and (\ref{1624a})
this yields the required bi-Lipschitz equivalence
\begin{equation}\label{1625}
\frac{1}{8}\log n\ \! d(v,w)\leq\rho_{0}(I(v),I(w))\leq\log n\ \! d(v,w)\ ,
\ \ \ v,w\in {\cal V}_{k}\ .
\end{equation}

We now extend $I$ to the whole ${\cal T}_{k}$ by defining it on
each edge $[v,v^{+}]\subset {\cal T}_{k}$; recall that
$[v,v^{+}]$ is identified with the unit interval of $\Re$ and therefore
there is a curve $\gamma_{v}:[0,1]\to [v,v^{+}]$ so that
$$
|\gamma_{v}(t)-\gamma_{v}(t')|=|t-t'|\ ,\ \ \ \gamma_{v}(1)=1    \ .
$$
To define the extension we join $c(v)$ and $c(v^{+})$ by the geodesic
segment in $\H ^{2}$ (the subarc of a Euclidean circle or a straight line
intersecting the $x_{1}$-axis orthogonally). Denote this by $[c(v),c(v^{+})]$.
By our geometric construction and properties of the geodesics of the
hyperbolic plane the interiors of any two such segments do not intersect.
Therefore the union of all $[c(v),c(v^{+})]$,
$v\in {\cal V}_{k}\setminus\{R\}$, forms a metric tree whose edges
$[c(v),c(v^{+})]$ are isometric to the closed intervals of $\Re$ of
lengths $\rho_{v}:=\rho(c(v),c(v^{+}))$. Let now
$\widetilde\gamma_{v}:[0,\rho_{v}]\to [c(v),c(v^{+})]$ be the canonical
parameterization of the geodesic $[c(v),c(v^{+})]$ so that
$$
\rho(\widetilde\gamma_{v}(t'),\widetilde\gamma_{v}(t''))=|t'-t''|\ ,\ \ \
\widetilde\gamma_{v}(\rho_{v})=\rho(c(v),c(v^{+})) (:=\rho_{v})\ .
$$
Let us extend the map $I$ to a point
$m=\gamma_{v}(t)\in [v,v^{+}]\subset {\cal T}_{k}$, $0\leq t\leq 1$, by
$$
I(m):=\widetilde\gamma_{v}(\rho_{v}t)\ .
$$
Then for $m,m'\in [v,v^{+}]$ we get
$$
\rho(I(m),I(m'))=\rho_{v}|t-t'|=\rho_{v}d(m,m')\ .
$$
Using now (\ref{1622}) and Lemma \ref{l163} we obtain the estimate
$$
\frac{1}{8}\log n\leq\rho_{v}\leq 4\log n\ .
$$
Together with the previous equality this yields
\begin{equation}\label{1629}
\frac{1}{8}\log n\ \! d(m,m')\leq\rho(I(m),I(m'))\leq 4\log n\ \! d(m,m')\ .
\end{equation}
Let now $m\in [v,v^{+}]$, $m'\in [w,w^{+}]$ and $v\neq w$. Since $\rho$
and $d$ are path metrics, the latter inequality together with
(\ref{1625}) and (\ref{1622}) imply that
$$
\frac{1}{64}\log n\cdot d(m,m')\leq\rho(I(m),I(m'))\leq 4\log n
\cdot d(m,m')\ .
$$
Hence $I$ is a bi-Lipschitz embedding of ${\cal T}_{k}$ into $\H ^{2}$ and
the inequality (\ref{164a}) of Proposition \ref{p161} holds with
$A:=\frac{1}{64}\log n$ and $B:=4\log n$.

This proves Proposition \ref{p161}.\ \ \ \ \ $\Box$

Using now the established embedding (\ref{164a}) we derive the
required inequality (\ref{164}). To complete the proof of the right-hand
side inequality of Theorem
\ref{te14a} for this case it remains to derive inequality
(\ref{eq72}) from Theorem \ref{c118}. To this end we use the
following particular case of the result from [BSh2, Proposition
5.33] with $n=1$.

There is a distance $\rho_{0}$ on $\H^{2}$ and a family of
measures $\{\mu_{x}\}_{x\in\H^{2}}$ such that
\begin{itemize}
\item[(a)]
$(\H^{2},\rho_{0})$ is of pointwise homogeneous type with respect to this
family;
\item[(b)]
$\rho_{0}$ is equivalent to the hyperbolic metric $\rho$;
\item[(c)]
For every ball $B_{R}^{0}(x):=\{y\in\H^{2}\ :\ \rho_{0}(x,y)\leq
R\}$ we have
$$
\mu_{x}(B_{R}^{0}(x))=2R^{2}\ .
$$
\end{itemize}

Apply now Theorem \ref{c118} for $M_{i}=(\H^{2},\rho_{0})$,
$1\leq i\leq n$. The above formulated statement shows that the
conditions of this theorem are true in this case and therefore
inequalities (\ref{e2.13'}) and (\ref{eq212}) yield
$$
\lambda(\oplus_{p}\{(\H^{2},\rho_{0})\}_{1\leq i\leq n})\leq c_{1}n
$$
for some numerical constant $c_{1}$ and $p=1,\infty$.
It remains to note that $\rho_{0}\sim\rho$ and therefore the
metrics of $\oplus_{p}\{\H^{2}\}_{1\leq i\leq n}$ and
$\oplus_{p}\{(\H^{2},\rho_{0})\}_{1\leq i\leq n}$ are equivalent with 
constants independent of $n$.

So inequality (\ref{eq72}) has proved, and the right-hand side
inequality of Theorem \ref{te14a} is
established for this case. \\
{\bf The general case.} Let now ${\cal T}_{i}$ be an arbitrary
metric tree with edges $e$ of lengths $l_{i}(e)>0$, $1\leq i\leq
n$. The argument of the previous subsection shows that in order to prove the
right-hand side inequality of the theorem it suffices to
derive the inequality
\begin{equation}\label{1626}
\lambda(\oplus_{p}\{{\cal T}_{i}\}_{1\leq i\leq n})\leq
256\lambda(\oplus_{p}\{\H ^{2}\}_{1\leq i\leq n})
\end{equation}
for arbitrary {\em finite} rooted metric trees ${\cal T}_{i}$.\\
We first establish this for finite ${\cal T}_{i}$ with edges of
lengths being rational numbers. Let $N$ be the least common
denominator of all these numbers for all $i$. Introduce a new
rooted metric tree ${\cal T}_{i}^{N}$, $1\leq i\leq n$, whose sets
of vertices ${\cal V}_{i}^{N}$ and edges ${\cal E}_{i}^{N}$ are
defined as follows. Let $e$ be an edge of ${\cal T}_{i}$ and
$l_{i}(e):=\frac{M_{i}(e)}{N}$ where $M_{i}(e)$ is a natural
number. Insert in this edge $M_{i}(e)-1$ equally distributed new
vertices; recall that $e$ is regarded as the closed interval of
$\Re$ of length $l_{i}(e)$. In this way we obtain the new rooted
tree ${\cal T}_{i}^{N}$, a triangulation of ${\cal T}_{i}$, that
we endow by the path metric $D_{i}^{N}$ induced by the metric
$Nd_{i}$ (here $d_{i}$ is the path metric on ${\cal T}_{i}$). Note
that every ${\cal T}_{i}^{N}$ has all edges of length 1 and therefore it
can be embedded into the infinite tree ${\cal T}_{k}$ with a
suitable $k$ (the same for all $i$). 
Hence the inequality (\ref{164}) yields \begin{equation}\label{1627}
\lambda(\oplus_{p}\{{\cal T}_{i}^{N}\}_{1\leq i\leq n})\leq
256\ \!\lambda(\oplus_{p}\{\H^{2}\}_{1\leq i\leq n})\ .
\end{equation}
On the other hand, $({\cal T}_{i},Nd_{i})$ is a metric subspace of
${\cal T}_{i}^{N}$, $1\leq i\leq n$, and therefore
$$
\lambda(\oplus_{p}\{{\cal
T}_{i}\}_{1\leq i\leq n})=\lambda(\oplus_{p}\{({\cal
T}_{i},Nd_{i})\}_{1\leq i\leq n})\leq
\lambda(\oplus_{p}\{{\cal T}_{i}^{N}\}_{1\leq i\leq n})\ .
$$
Together with (\ref{1627}) this implies the required result (\ref{1626}).

Consider now the general situation of finite rooted metric trees
${\cal T}_{i}$ with arbitrary lengths of edges. Given $\epsilon>0$
one replaces the metric of ${\cal T}_{i}$, $1\leq i\leq n$, by a
path metric $d_{i,\epsilon}$ which remains to be linear on edges
but such that lengths of edges $l_{i,\epsilon}(e)$ in this metric
are rational numbers satisfying
$$
l_{i}(e)\leq l_{i,\epsilon}(e)\leq (1+\epsilon)l_{i}(e)\ .
$$
Let ${\cal T}_{i,\epsilon}$ be a rooted metric tree with the
underlying set ${\cal T}_{i}$ and the path metric
$d_{i,\epsilon}$. It is already proved that
$$
\lambda(\oplus_{p}\{{\cal T}_{i,\epsilon}\}_{1\leq i\leq n})\leq
256\ \!\lambda(\oplus_{p}\{\H ^{2}\}_{1\leq i\leq n})\ .
$$
On the other hand, the identity map ${\cal T}_{i}\to {\cal
T}_{i,\epsilon}$ is a bi-Lipschitz equivalence with the constant
of equivalence equals $1+\epsilon$. Therefore
$$
\lambda(\oplus_{p}\{{\cal T}_{i}\}_{1\leq i\leq n})\leq
(1+\epsilon)^{2}\lambda(\oplus_{p}\{{\cal T}_{i,\epsilon}\}_{1\leq i\leq n}).
$$
Since $\epsilon$ is arbitrary, the last two inequalities prove
(\ref{1626}) in the general case. 

Thus we have proved the right-hand side inequality of Theorem \ref{te14a}.

To complete the proof of the theorem it remains to prove the lower estimate
$$
\lambda(\oplus_{p}\{{\cal T}_{i}\}_{1\leq i\leq n})\geq c_{0}\sqrt{n}
$$
for $p=1,\infty$ and $c_{0}>0$ independent of $n$. 

To this end choose in every ${\cal T}_{i}$ a path $P_{i}$ incident to the
root of ${\cal T}_{i}$ of length $l_{i}>0$. Then the interval $[0,l_{i}]$
is isometrically embedded into ${\cal T}_{i}$ and the parallelepiped
$\Pi:=\prod_{i=1}^{n}[0,l_{i}]$ equipped  with the $l_{p}^{n}$-metric is
an isometric part of $\oplus_{p}\{{\cal T}_{i}\}_{1\leq i\leq n}$. By
Corollary \ref{c13} with $S:=\Pi$ and $\delta(x):=2(x-c)$, $x\in\Re^{n}$,
where $c$ is the center of $\Pi$, one has
$$
\lambda(\Pi)=\lambda(l_{p}^{n})
$$
and the latter is greater than $\lambda_{conv}(l_{p}^{n})$. This, in turn,
is at least $c_{0}\sqrt{n}$ for $p=1,\infty$ with $c_{0}>0$ independent
of $n$ (see Theorem \ref{te1.11}). 

This completes the proof of the required lower bound (and the theorem).\ \ \
\ \ $\Box$

\sect{Proof of Theorem \ref{te14}}
We begin with the case of $M$ possessing the WTP. Assume that
$M$ has the Lipschitz preserving linear extension property, but
$\lambda(M)=\infty$. The latter implies existence of a sequence of finite
sets $F_{j}$ with $\lambda(F_{j})\geq j$, $j\in\N$; see Corollary \ref{c12}.
This, in turn, leads to the inequalities
\begin{equation}\label{e41}
\inf\{||E||\ :\ E\in Ext(F_{j},M)\}\geq j\ ,\ j\in\N\ .
\end{equation}
Using the WTP of $M$ to choose an appropriate sequence of $C$-isometries
$\sigma_{j}$ such that for $G_{j}:=\sigma_{j}(F_{j})$ the following
holds
\begin{equation}\label{e42}
dist(G_{j},\cup_{i\neq j}G_{i})\geq C\ \! diam\ F_{j}\ ,\ j\in\N\ .
\end{equation}
For every $j\in\N$, fix a point $m_{j}^{*}\in G_{j}$. From (\ref{e42}) we
derive that an operator $N_{j}$ given for every
$f\in Lip(G_{j})$ by
\begin{equation}\label{e43}
(N_{j}f)(m):=
\left\{
\begin{array}{ccc}
f(m)\ ,&&m\in G_{j}\\
f(m_{j}^{*})\ ,&&m\in\cup_{i\neq j}G_{i}
\end{array}
\right.
\end{equation}
belongs to $Ext(G_{j},G_{\infty})$ where $G_{\infty}:=\cup_{i\in\N}G_{i}$,
and, besides,
\begin{equation}\label{e44}
||N_{j}||=1\ .
\end{equation}
In fact, if $f\in Lip(G_{j})$ and $m'\in G_{j}$,
$m''\in G_{\infty}\setminus G_{j}=\cup_{i\neq j}G_{i}$, then
$$
\begin{array}{c}
|(N_{j}f)(m')-(N_{j}f)(m'')|=|f(m')-f(m_{j}^{*})|\leq\\
\\
||f||_{Lip(G_{j})}diam\ G_{j}\leq (C\ \! diam\ F_{j})
||f||_{Lip(G_{j})}\ .
\end{array}
$$
Together with (\ref{e42}) this leads to
$$
|(N_{j}f)(m')-(N_{j}f)(m'')|\leq ||f||_{Lip(G_{j})}d(m',m'')\ .
$$
Since this holds trivially for all other choices of $m',m''$, the equality
(\ref{e44}) is done.

Now by the ${\cal LE}$ of $M$ there is an operator $E\in Ext(G_{\infty},M)$ with
$||E||\leq A$ for some $A>0$. By (\ref{e44}) the operator
$E_{j}:=EN_{j}\in Ext(G_{j},M)$ and $||E_{j}||\leq A$. Then an operator
$\widetilde E_{j}$ given by the formula
$$
(\widetilde E_{j}f)(m):=(E_{j}(f\circ\sigma_{j}^{-1}))(\sigma_{j}(m))\ ,\ \
m\in M\ ,\ f\in Lip(F_{j})\ ,
$$
with the above introduced $C$-isometries $\sigma_{j}$
belongs to $Ext(F_{j},M)$ and its norm is bounded by $C^{2}A$. Comparing with
(\ref{e41}) to get for each $j$
$$
C^{2}A\geq j\ ,
$$
a contradiction.

Let now $M$ be proper. In order to prove that
$\lambda(M)<\infty$ we need
\begin{Lm}\label{l41}
For every $m\in M$ there is an open ball $B_{m}$ centered at $m$ such that
$\lambda(B_{m})<\infty$.
\end{Lm}
{\bf Proof.} Assume that this assertion does not hold for some $m$. Then
there is a sequence of balls $B_{i}:=B_{r_{i}}(m)$, $i\in\N$, centered at
$m$ of radii $r_{i}$ such that $\lim_{i\to\infty}r_{i}=0$ and
$\lim_{i\to\infty}\lambda(B_{i})=\infty$. According to Corollary \ref{c12}
this implies existence of finite subsets
$F_{i}\subset B_{i}$, $i\in\N$, such that
\begin{equation}\label{e45}
\inf\{||E||\ :\ E\in Ext(F_{i},B_{i})\}\to\infty\ ,\ {\rm as}\ i\to\infty\ .
\end{equation}
We may and will assume that $m\in F_{j}$, $j\in\N$.
Otherwise we replace $F_{j}$ by $G_{j}:=F_{j}\cup\{m\}$ and show
that (\ref{e45}) remains true for $G_{i}$ and $B_{i}$, $i\in\N$. In fact,
let $L_{i}$ be an operator given for every $f\in Lip(F_{i})$ by
$$
(L_{i}f)(m'):=\left\{
\begin{array}{ccc}
f(m_{i}),&\rm if&m'=m\\
f(m'),&\rm if&m'\in F_{i}
\end{array}
\right.
$$
where $m_{i}$ is the closest to $m$ point from $F_{i}$. Then
$L_{i}\in Ext(F_{i},G_{i})$ and $||L_{i}||\leq 2$, since
$$
|(L_{i}f)(m)-(L_{i}f)(m')|=|f(m_{i})-f(m')|\leq ||f||_{Lip(F_{i})}
d(m_{i},m')\leq 2||f||_{Lip(F_{i})}d(m,m')\ .
$$
If now (\ref{e45}) does not hold for $\{G_{i}\}$ substituted for
$\{F_{i}\}$, then there is a sequence
$E_{i}\in Ext(G_{i},B_{i})$ such that $\sup_{i}||E_{i}||<\infty$. But
then the same will be true for the norms of $\widetilde E_{i}:=E_{i}L_{i}\in
Ext(F_{i},B_{i})$, $i\in\N$, in contradiction with (\ref{e45}).

The proof will be now finished by the argument of section 4.1. Actually,
choose a subsequence $F_{i_{k}}\subset B_{i_{k}}:=B_{r_{i_{k}}}(m)$,
$k\in\N$, such that
$$
r_{i_{k+1}}<\min\{r_{i_{k}},\
dist(F_{i_{k+1}}\setminus\{m\},\cup_{s< k+1}F_{i_{s}}\setminus\{m\})\}\ .
$$
Without loss of generality we assume that the sequence $\{F_{i}\}$ already
satisfies this condition, i.e.,
\begin{equation}\label{e46}
r_{i+1}<\min\{r_{i},\
dist(F_{i+1}\setminus\{m\},\cup_{s<i+1}F_{s}\setminus\{m\})\}\
\end{equation}
Set $F_{\infty}:=\bigcup_{s\in\N}F_{s}$ and show that
\begin{equation}\label{e47}
Ext(F_{\infty},M)=\emptyset
\end{equation}
which gives the required contradiction to the ${\cal LE}$ of $M$. To prove
this
we choose the point $m$ as a marked point of $M$. Then all $F_{i}$
are subspaces of $(M,m)$ and $f(m)=0$ if $f\in Lip_{0}(F_{i})$, $i\in\N$.
Define now an operator $N_{i}$ by
$$
(N_{i}f)(m'):=\left\{
\begin{array}{ccc}
f(m'),&\rm if&m'\in F_{i}\\
0,&\rm if&m'\in F_{\infty}\setminus F_{i}\ .
\end{array}
\right.
$$
Then for $f\in Lip_{0}(F_{i})$ and $m'\in F_{i}\setminus\{m\}$ and
$m''\in F_{\infty}\setminus F_{i}$ we have
$$
|(N_{i}f)(m')-(N_{i}f)(m'')|=|f(m')-f(m)|\leq ||f||_{Lip_{0}(F_{i})}
d(m',m)\ .
$$
Moreover, $m''\in B_{j}$ for some $j\neq i$. Assume, first, that $j>i$.
Then by (\ref{e46})
$$
\begin{array}{c}
d(m',m)\leq d(m',m'')+d(m'',m)\leq d(m',m'')+r_{j}\leq\\
\\
d(m',m'')+dist(F_{j}\setminus\{m\},F_{i}\setminus\{m\})\leq
2d(m',m'')\ .
\end{array}
$$
If now $j<i$, then by (\ref{e46}) we have
$$
d(m',m)\leq r_{i}<dist(F_{i}\setminus\{m\},F_{j}\setminus\{m\})\leq 
d(m',m'').
$$
Combining these we prove that
$N_{i}\in Ext(F_{i},F_{\infty})$ and $||N_{i}||\leq 2$.

If now (\ref{e47}) is not true, then there is an
operator $E\in Ext(F_{\infty},M)$, and so every operator
$\widetilde E_{i}:=EN_{i}$ belongs to $Ext(F_{i},M)$ and
$||\widetilde E_{i}||\leq 2||E||$, $i\in\N$, a contradiction to
(\ref{e45}). Hence (\ref{e47}) holds and the proof is complete.\ \ \ \ \
$\Box$
\begin{R}\label{r43}
{\rm In this proof properness of $M$ is not used.}
\end{R}

To prove the next important fact on $\lambda$ we need

\begin{Lm}\label{l44}
Let ${\cal U}$ be a finite cover of a compact set $C\subset M$ by open
sets.
Then there is a partition of unity $\{\rho_{U}\}_{U\in {\cal U}}$ on
$C$ subordinate to ${\cal U}$ such that every $\rho_{U}$ is Lipschitz
with a constant depending only on the cover.
\end{Lm}

Let us recall its proof.

Define $d_{U}:M\rightarrow\Re_{+}$ by
$$
d_{U}(m):=dist(m,M\setminus U)\ ,\ \ m\in M\ .
$$
This is supported on $U$ and is Lipschitz with constant 1.
Moreover, $\sum_{U\in {\cal U}}d_{U}>0$ on $C$, as ${\cal U}$ is a cover of
$C$. Putting now
$$
\rho_{U}(m):=\frac{d_{U}(m)}{\sum_{U\in {\cal U}}d_{U}(m)}\ ,\ \
m\in C\cap U\ ,\ \ U\in {\cal U}\ ,
$$
we get the required partition.\ \ \ \ \ $\Box$

The next result implies finiteness of $\lambda(M)$ for compact $M$.
\begin{Lm}\label{l45}
For every compact set $C\subset M$ the constant $\lambda(C)$ is finite.
\end{Lm}
{\bf Proof.}
We have to show that for every $S\subset C$ there is an operator
$E\in Ext(S,C)$ such that
\begin{equation}\label{e49}
\sup_{S}||E||<\infty\ .
\end{equation}
We may and will assume that $S$ and $C$ are subspaces of $(M,m^{*})$ so that
$f(m^{*})=0$ for $f$ belonging to $Lip_{0}(S)$ or $Lip_{0}(C)$. By
compactness of $C$ and Lemma \ref{l41} there is a finite cover
$\{U_{i}\}_{1\leq i\leq n}$ of $C$ by open balls such that
for some constant $A>0$ depending only on $C$ we have
$$
\lambda(U_{i})<A\ ,\ \ 1\leq i\leq n\ .
$$
By the definition of $\lambda$ this implies existence of
$E_{i}\in Ext(S\cap U_{i},U_{i})$ with
\begin{equation}\label{e410}
||E_{i}||\leq A\ ,\ \ 1\leq i\leq n\ .
\end{equation}
For $f\in Lip_{0}(S)$ one sets
\begin{equation}\label{e411}
f_{i}:=\left\{
\begin{array}{ccc}
f|_{S\cap U_{i}}\ ,&\rm if& S\cap U_{i}\neq\emptyset\\
0\ ,&\rm if&S\cap U_{i}=\emptyset\
\end{array}
\right.
\end{equation}
and introduces a function $f_{ij}$ given on $U_{i}\cap U_{j}$ by
\begin{equation}\label{e412}
f_{ij}:=\left\{
\begin{array}{ccc}
E_{i}f_{i}-E_{j}f_{j}\ ,&\rm if& U_{i}\cap U_{j}\neq\emptyset\\
0\ ,&\rm if& U_{i}\cap U_{j}=\emptyset\ ,
\end{array}
\right.
\end{equation}
here $E_{i}f_{i}:=0$, if $f_{i}=0$.
Then (\ref{e410}) implies that
\begin{equation}\label{e413}
||f_{ij}||_{Lip(U_{i}\cap U_{j})}\leq 2A\ ;
\end{equation}
besides, we get
\begin{equation}\label{e414}
f_{ij}=0\ \ \ {\rm on}\ \ \ S\cap U_{i}\cap U_{j}\ .
\end{equation}
At last, introduce the function $g_{i}$ on $C\cap U_{i}$ by
\begin{equation}\label{e415}
g_{i}(m):=\sum_{1\leq j\leq n}\rho_{j}(m)f_{ij}(m)\ ,\ \ \ m\in C\cap U_{i}\ ,
\end{equation}
where $\rho_{j}:=\rho_{U_{j}}$, $1\leq j\leq n$, is the partition of
unity of Lemma \ref{l44}.\\
A straightforward computation leads to the equalities:
\begin{equation}\label{e416}
g_{i}-g_{j}=f_{ij}\ \ \ {\rm on}\ \ \ U_{i}\cap U_{j}\cap C
\end{equation}
and, moreover,
\begin{equation}\label{e417}
g_{i}|_{S\cap U_{i}}=0\ .
\end{equation}
Introduce now an operator $E$ on $f\in Lip_{0}(S)$ by the formula
\begin{equation}\label{e418}
(Ef)(m):=(E_{i}f_{i}-g_{i})(m)\ ,\ \ \ {\rm if}\ \ \ m\in U_{i}\cap C\ ,
\end{equation}
and show that $E$ is an extension operator. In fact, if $m\in S$, then
$m\in S\cap U_{i}$ for some $1\leq i\leq n$ and by (\ref{e417}) and
(\ref{e411}) we get
$$
(Ef)(m)=(E_{i}f_{i})(m)=f(m)\ .
$$

Show now that $E\in Ext(S,C)$ and $||E||$ is bounded by a constant
depending only on $C$. To this end we denote by $\delta=\delta(C)>0$
the {\em Lebesgue number} of the cover ${\cal U}$, see, e.g., [Ke].
So every subset of $\cup_{i=1}^{n}U_{i}$ of diameter at most $\delta$
lies in one of $U_{i}$. Using this we first establish the corresponding
Lipschitz estimate for $m',m''\in (\cup_{i=1}^{n}U_{i})\cap C$ with
\begin{equation}\label{e419}
d(m',m'')\leq\delta\ .
\end{equation}
In this case both $m',m''\in U_{i_{0}}$ for some $i_{0}$. Further,
(\ref{e415})-(\ref{e418}) imply that for $m\in U_{i_{0}}\cap C$
\begin{equation}\label{e420}
(Ef)(m)=\sum_{U_{i}\cap U_{i_{0}}\neq\emptyset}
(\rho_{i}E_{i}f_{i})(m)\ .
\end{equation}
In this sum each $\rho_{i}$ is Lipschitz with a constant $L(C)$ depending
only on $C$ and $0\leq\rho_{i}\leq 1$. In turn, $E_{i}f_{i}$ is
Lipschitz on $U_{i}\cap U_{i_{0}}$ with constant $A ||f||_{Lip_{0}(S)}$,
see (\ref{e410}) (recall that $E_{i}f_{i}=0$ if $S\cap U_{i}=\emptyset$).
If now $m\in U_{i}$ with $S\cap U_{i}\neq\emptyset$, then for arbitrary
$m_{i}\in S\cap U_{i}$
$$
\begin{array}{c}
|(E_{i}f_{i})(m)|\leq |(E_{i}f_{i})(m)-(E_{i}f_{i})(m_{i})|+
|(E_{i}f_{i})(m_{i})|\leq\\
\\
A||f||_{Lip_{0}(S)}d(m,m_{i})+
|f(m_{i})-f(m^{*})|\leq A||f||_{Lip_{0}(S)}(d(m,m_{i})+d(m_{i},m^{*}))\ .
\end{array}
$$
This implies for all $m\in C\cap U_{i}$ the inequality
\begin{equation}\label{e421}
|(E_{i}f_{i})(m)|\leq 2A\ diam\ \! C\ ||f||_{Lip_{0}(S)}\ .
\end{equation}
Together with
(\ref{e420}) this leads to the estimate
\begin{equation}\label{e422}
|(Ef)(m')-(Ef)(m'')|\leq A n(2L(C)\ diam\ \! C+1)||f||_{Lip_{0}(S)}
d(m',m'')\ ,
\end{equation}
provided $m',m''\in U_{i_{0}}\cap C$.
To prove a similar estimate for $d(m',m'')>\delta$, $m',m''\in C$, we note
that the left-hand side in (\ref{e422}) is bounded by
$$
2\sup_{m\in C}|(Ef)(m)|\leq 4A\ diam\ \! C\ ||f||_{Lip_{0}(S)}\ ,
$$
see (\ref{e421}). In turn, the right-hand side of the last
inequality is $\leq 4\delta^{-1}A\ diam\ C\
||f||_{Lip_{0}(S)}d(m',m'')$. Together this implies that $E$
belongs to $Ext(S,C)$ and its norm is bounded by a constant
depending only on $C$.\ \ \ \ \ $\Box$

Our last basic step is a lemma which formulation uses the notation
\begin{equation}\label{e423}
\lambda(S,M):=\inf\{||E||\ :\ E\in Ext(S,M)\}\ ;
\end{equation}
here we set $\lambda(S,M):=0$, if $S=\emptyset$.
\begin{Lm}\label{l46}
Assume that for a sequence of finite subsets
$F_{i}\subset (M,m^{*})$, $i\in\N$,
\begin{equation}\label{e424}
\sup_{i}\lambda(F_{i},M)=\infty\ .
\end{equation}
Then for every closed ball $\overline B$ centered at $m^{*}$
\begin{equation}\label{e425}
\sup_{i}\lambda(F_{i}\setminus\overline B,M)=\infty\ .
\end{equation}
\end{Lm}
{\bf Proof.} Arguing as in the proof of Lemma
\ref{l41}, we can assume that
$m^{*}\in F_{i}$, $i\in\N$.
If now (\ref{e425}) is not true, then for some $A_{1}>0$ there is
a sequence of operators $E_{i}^{1}$, $i\in\N$, such that
\begin{equation}\label{e426}
E_{i}^{1}\in Ext(F_{i}\setminus\overline B,M)\ \ \ {\rm and}\ \ \
||E_{i}^{1}||\leq A_{1}\ ;
\end{equation}
here we set $E_{i}^{1}:=0$, if $F_{i}\setminus\overline
B=\emptyset$. Let $2B$ be the open ball centered at $m^{*}$ and of
twice more radius than that of $B$. Introduce an open cover of $M$
by
\begin{equation}\label{e427}
U_{1}:=M\setminus\overline B\ \ \ {\rm and}\ \ \ U_{2}:=2B\ ,
\end{equation}
and let $\{\rho_{1},\rho_{2}\}$ be the corresponding Lipschitz partition of
unity, cf. Lemma \ref{l44}, given by
$$
\rho_{j}(m):=\frac{d_{U_{j}}(m)}{d_{U_{1}}(m)+d_{U_{2}}(m)}\ ,\ \
m\in M\ ,\ j=1,2\ .
$$
By this definition
\begin{equation}\label{new}
|\rho_{j}(m_{1})-\rho_{j}(m_{2})|\leq
\frac{3d(m_{1},m_{2})}{\max_{k=1,2}\{d_{U_{1}}(m_{k})+d_{U_{2}}(m_{k})\}}\ ,
\ \ \ m_{1},m_{2}\in M\ .
\end{equation}
Set now $H_{i}:=F_{i}\cap 2B$, $i\in\N$. Since these are subsets of the
compact set
$2\overline B$, Lemma \ref{l45} gives
$$
\sup_{i}\lambda(H_{i},2B)\leq\lambda(2\overline B)<\infty\ .
$$
This, in turn, implies existence of operators $E_{i}^{2}$, $i\in\N$,
such that
\begin{equation}\label{e428}
E_{i}^{2}\in Ext(H_{i},2B)\ \ \
{\rm and}\ \ \ ||E_{i}^{2}||\leq A_{2}
\end{equation}
with $A_{2}$ independent of $i$. We now follow the proof of Lemma
\ref{l45} where the set $S$ and the compact set $C\supset S$ are replaced by
$F_{i}$ and the (noncompact) space $M$, respectively, and the cover
(\ref{e427}) is used. Since
$$
H_{i}=F_{i}\cap U_{2}\ \ \ {\rm and}\ \ \
F_{i}\setminus\overline B=F_{i}\cap U_{1}\ ,
$$
we can use in our derivation the operators $E_{i}^{j}$, $j=1,2$, instead of
those in (\ref{e410}). By (\ref{e426}) and (\ref{e428}) inequalities
similar to (\ref{e410}) hold for these operators. Then we set
$f_{j}:=f|_{F_{i}\cap U_{j}}$, and define for $f\in Lip_{0}(F_{i})$ functions
$f_{12}:=-f_{21}$ on $U_{1}\cap U_{2}$, $g_{1}$ on $U_{1}$ and $g_{2}$ on
$U_{2}$ by
$$
f_{12}:=E_{i}^{1}f_{1}-E_{i}^{2}f_{2}\ ,\ \ \  g_{1}:=\rho_{2}f_{12}\ ,\ \ \
g_{2}:=\rho_{1}f_{21}\ .
$$
At last, we introduce the required operator $E_{i}$ on $Lip_{0}(F_{i})$ by
$$
(E_{i}f)(m):=(E_{i}^{j}f_{j})(m)-g_{j}(m)\ ,\ \ \ \ m\in U_{j}\ ,\ j=1,2\ .
$$
As in Lemma \ref{l45} $E_{i}$ is an operator extending functions from
$F_{i}$ to the whole $M$.
To estimate the Lipschitz constant of $E_{i}f$  we extend each
$E_{i}^{j}f_{j}$ outside $U_{j}$ so that the extensions
$\widetilde E_{j}$ satisfy
$$
||\widetilde E_{j}||_{Lip(M)}=||E_{i}^{j}f_{j}||_{Lip(U_{j})}\ .
$$
Now, the definition of $E_{i}$ implies that
$$
E_{i}f:=\rho_{1}\widetilde E_{1}+\rho_{2}\widetilde E_{2}\ .
$$
Assume without loss of generality that $F_{i}\cap U_{1}\neq\emptyset$, and
choose a point $m'\in F_{i}\cap U_{1}$. Then as in the proof of (\ref{e421})
for arbitrary $m\in M$ we obtain
$$
|\widetilde E_{j}(m)|\leq\left\{
\begin{array}{ccc}
A_{1}||f||_{Lip_{0}(F_{i})}(d(m,m')+d(m',m^{*})),&{\rm if}&j=1\\
\\
A_{2}||f||_{Lip_{0}(F_{i})}d(m,m^{*}),&{\rm if}&j=2
\end{array}
\right.
$$
This implies for all $m$ the inequality
\begin{equation}\label{new1}
|\widetilde E_{j}(m)|\leq\ A(d(m,m^{*})+d(m',m^{*}))||f||_{Lip_{0}(F_{i})}
\end{equation}
with $A:=2\max(A_{1},A_{2})$.\\
Together with (\ref{new}) this leads to the estimate
\begin{equation}\label{new2}
\begin{array}{c}
\displaystyle |(E_{i}f)(m_{1})-(E_{i}f)(m_{2})|\leq\\
\\
\displaystyle
\left(\frac{3(d(m_{1},m^{*})+d(m',m^{*}))}
{\max_{k=1,2}\{d_{U_{1}}(m_{k})+d_{U_{2}}(m_{k})\}}
+1\right)2A\ \!||f||_{Lip_{0}(F_{i})}d(m_{1},m_{2})\ .
\end{array}
\end{equation}
Since $\max_{k=1,2}\{d_{U_{1}}(m_{k})+d_{U_{2}}(m_{k})\}\geq R$, the radius
of $B$, and
$$
\lim_{d(m_{1},m^{*})\to\infty}\frac{d(m_{1},m^{*})+d(m',m^{*})}
{d_{U_{1}}(m_{1})+d_{U_{2}}(m_{1})}=1\ ,
$$
(\ref{new2}) implies that $E_{i}\in Ext(F_{i},M)$ and its norm is bounded
by a constant independent of $i$.
By definition (\ref{e423}) this gets
$$
\sup_{i}\lambda(F_{i},M)<\infty
$$
in a contradiction with (\ref{e424}).\ \ \ \ \ $\Box$

Now we will complete the proof of Theorem \ref{te14}. Recall that it
has already proved for compact $M$, see Lemma \ref{l45}. So it remains to
consider the case of a proper $M$ with
\begin{equation}\label{e429}
diam\ \! M=\infty\ .
\end{equation}
In this case we will show that
\begin{equation}\label{e430}
\sup_{F}\lambda(F,M)<\infty\ ,
\end{equation}
where $F$ is running through all finite subsets $F\subset (M,m^{*})$.
Since $\sup_{F}\lambda(F)$ is bounded by the supremum in (\ref{e430}), this
leads to finiteness of $\lambda(M)$, see Corollary \ref{c12}, and proves the
result.

Let, to the contrary, (\ref{e430}) does not hold. Then there is a sequence
of finite subsets $F_{i}\subset (M,m^{*})$, $i\in\N$, that satisfies the
assumption of Lemma \ref{l46}, see (\ref{e424}). We use this to construct
a sequence $G_{i}$, $i\in\N$, such that
\begin{equation}\label{e431}
\lambda(G_{i},M)\geq i-1\ \ \ {\rm and}\ \ \ dist(G_{i+1},G_{i})\geq
dist(m^{*},G_{i})\ ,\ i\in\N\ .
\end{equation}
As soon as it is done we set $G_{\infty}:=\cup_{i\in\N}G_{i}$, and use
with minimal changes the argument of Lemma \ref{l41} to show that
$$
Ext(G_{\infty},M)=\infty\ .
$$
Since this contradicts to the ${\cal LE}$ of $M$, the result will be done.

To construct the required $\{G_{i}\}$,  set $G_{1}:=F_{1}$ and assume that
the first $j$ terms of this sequence have already defined. Choose the
closed ball $\overline B$ such that
$$
dist(G_{j},M\setminus\overline B)\geq dist(m^{*},G_{j})\ ;
$$
it exists because of (\ref{e429}). Then apply Lemma \ref{l46} to find
$F_{i(j)}$  such that
$$
\lambda(F_{i(j)}\setminus\overline B,M)\geq j-1\ .
$$
Setting $G_{j+1}:=F_{i(j)}\setminus\overline B$ we obtain the next term
satisfying the condition (\ref{e431}).

The proof is complete.\ \ \ \ \ $\Box$
\sect{Proof of Theorem \ref{te1.11}}
We first prove that
\begin{equation}\label{e51}
\lambda_{conv}(l_{2}^{n})=1\ .
\end{equation}
Since the lower bound $1$ is evident, we have to prove that
\begin{equation}\label{e52}
\lambda_{conv}(l_{2}^{n})\leq 1\ .
\end{equation}

Let $C\subset (l_{2}^{n},0)$ be a closed convex set containing 0, and
$p_{C}(x)$ be the (unique) closest to $x$ point from $C$.
Then, see, e.g., [BL, Sect.3.2], the {\em metric projection} $p_{C}$ is
Lipschitz and
\begin{equation}\label{e53}
||p_{C}(x)-p_{C}(y)||_{2}\leq ||x-y||_{2}\ .
\end{equation}
Using this we introduce a linear operator $E$ given on $Lip_{0}(C)$ by
\begin{equation}\label{e54}
(Ef)(x):=(f\circ p_{C})(x)\ ,\ \ \ x\in l_{2}^{n}\ .
\end{equation}
Since $p_{C}$ is identity on $C$ and $p_{C}(0)=0$ as $0\in C$, this
operator belongs to $Ext(C,l_{2}^{n})$. Moreover, by (\ref{e53})
$$
||(Ef)(x)-(Ef)(y)||_{2}\leq ||f||_{Lip_{0}(C)}||x-y||_{2}\ ,
$$
i.e., $||E||\leq 1$, and (\ref{e52}) is established.

Using now the inequalities
\begin{equation}\label{e55}
||x||_{p}\leq ||x||_{2}\leq n^{\frac{1}{2}-\frac{1}{p}}\ ||x||_{p}\ \ \
{\rm for}\ \ \ 2\leq p\leq\infty\ ,
\end{equation}
and
\begin{equation}\label{e56}
n^{\frac{1}{2}-\frac{1}{p}}\ ||x||_{p}\leq ||x||_{2}\leq ||x||_{p}
\ \ \ {\rm for}\ \ \ 1\leq p\leq 2
\end{equation}
we derive from (\ref{e52}) the required upper bound:
\begin{equation}\label{e57}
\lambda(l_{p}^{n})\leq n^{\left|\frac{1}{2}-\frac{1}{p}\right|}\ .
\end{equation}

In order to prove the lower estimate we need the following
result where Banach spaces are regarded as punctured metric
spaces with $m^{*}=0$.
\begin{Proposition}\label{p51}
Let $Y$ be a linear subspace of a finite dimensional Banach space $X$, and
an operator $E$ belongs to $Ext(Y,X)$. Then there is a linear projection
$P$ from $X$ onto $Y$ such that
\begin{equation}\label{e58}
||P||\leq ||E||\ .
\end{equation}
\end{Proposition}
{\bf Proof.} We use an argument similar to that in
[P, Remarks to $\S 2$]. First, we introduce an operator
$S:Lip_{0}(Y)\rightarrow Lip_{0}(X)$ given at $z\in X$ by
\begin{equation}\label{e59}
(Sf)(z):=\int_{X}\left\{\int_{Y}[(Ef)(x+y+z)-(Ef)(x+y)]dy\right\}dx\ .
\end{equation}
Here $\int_{A}\dots da$ is a translation invariant mean on the space
$l_{\infty}(A)$ of all bounded functions on an abelian group $A$, see, e.g.,
[HR]. Since the function within $[\ \ ]$ is bounded for every fixed $z$
(recall that $Ef\in Lip(X)$), this operator is well-defined.
Moreover, as $\int_{A}da=1$ we get
\begin{equation}\label{e510}
||Sf||_{Lip_{0}(X)}\leq ||E||\cdot ||f||_{Lip_{0}(Y)}\ .
\end{equation}
By translation invariance of $dx$ we then derive from (\ref{e59}) that
$$
(Sf)(z_{1}+z_{2})=(Sf)(z_{1})+(Sf)(z_{2})\ ,\ \ \ z_{1},z_{2}\in X\ .
$$
Together with (\ref{e510}) and the equality
\begin{equation}\label{e511}
||f||_{Lip_{0}(X)}=||f||_{X^{*}}\ ,\ \ f\in X^{*}\ ,
\end{equation}
this shows that $Sf$ belongs to $X^{*}$ and therefore $S$ maps
$Lip_{0}(Y)$ linearly and continuously in $X^{*}$. Further, $Y^{*}$ is a
linear subset of $Lip_{0}(Y)$ whose norm coincides with that induced
from $Lip_{0}(Y)$. Therefore the restriction
$$
T:=S|_{Y^{*}}
$$
is a linear bounded operator from $Y^{*}$ to $X^{*}$. Show that $T$
satisfies
\begin{equation}\label{e512}
(Tf)(z)=f(z)\ ,\ \ z\in Y\ ,
\end{equation}
i.e., $T$ is an extension from $Y^{*}$. To this end write
$$
\begin{array}{rl}
\displaystyle
(Tf)(z)=\int_{X}\left\{\int_{Y}[(Ef)(x+y+z)-(Ef)(y+z)]dy\right\}dx\ +\\
\\
\displaystyle
\int_{X}\left\{\int_{Y}[(Ef)(y+z)-(Ef)(x+y)]dy\right\}dx\ .
\end{array}
$$
Since $z\in Y$ and $dy$ is translation invariant with respect to
translations by elements of $Y$, we can omit $z$ in the first summand.
Moreover, $(Ef)(y)=f(y)$ for $f\in Y^{*}\subset Lip_{0}(Y)$. Thus, the
right-hand side is equal to
$$
\int_{X}\left\{\int_{Y}[(Ef)(x+y)-f(y)+(Ef)(y+z)-(Ef)(x+y)]dy\right\}dx\ .
$$
Since $(Ef)(y+z)=f(y+z)=f(y)+f(z)$, this integral equals
$$
f(z)\int_{X}dx\int_{Y}dy=f(z)\ ,
$$
and (\ref{e512}) is done.

Consider now the conjugate to $T$ operator $T^{*}$ acting from $X^{**}=X$
to $Y^{**}=Y$. Since $T$ is a linear extension operator from $Y^{*}$,
its conjugate is a projection onto $Y$. At last, (\ref{e510}) and
(\ref{e511}) give for the norm of this projection the
required estimate
$$
||T^{*}||=||T||\leq ||E||\ .
$$

The proposition is proved.\ \ \ \ \ $\Box$
\begin{Proposition}\label{p52}
Let $X$ be either $l_{1}^{n}$ or $l_{\infty}^{n}$. Then there is a subspace
$Y\subset X$ such that $dim\ \!Y=[n/2]$ and its projection constant
\footnote{$\pi(Y,X):=\inf||P||$ where $P$ runs through all linear
projections from $X$ onto $Y$.} satisfies
\begin{equation}\label{e513}
\pi(X,Y)\geq c_{0}\sqrt{n}
\end{equation}
with $c_{0}$ independent of $n$.
\end{Proposition}

{\bf Proof.} The inequality follows from Theorem 1.2 of the paper [S]
by Sobczyk with the optimal $c_{0}$ greater than $1/4$.\ \ \ \ \ $\Box$

We now complete the proof of Theorem \ref{te1.11}. Applying
Propositions \ref{p51} and \ref{p52} we get for an arbitrary
$E\in Ext(Y,l_{1}^{n})$ the inequality
\begin{equation}\label{e517}
||E||\geq c_{0}\sqrt{n}
\end{equation}
with $c_{0}>0$ independent of $n$. A similar estimate is valid for
$E\in Ext(Y^{\perp},l_{\infty}^{n})$, as well. Hence for $p=1,\infty$
\begin{equation}\label{e518}
\lambda_{conv}(l_{p}^{n})\geq c_{0}\sqrt{n}\ .
\end{equation}
Using this estimate for $p=1$ and applying a similar to
(\ref{e56}) inequality comparing $||x||_{1}$ and $||x||_{p}$ we get for
$1\leq p\leq 2$ the estimate
$$
\lambda(l_{p}^{n})\geq c_{0}n^{\frac{1}{p}-\frac{1}{2}}\ .
$$
Then using (\ref{e518}) for $p=\infty$ and a similar to
(\ref{e55}) inequality comparing $||x||_{p}$ and $||x||_{\infty}$ we get for
$2\leq p\leq\infty$ the inequality
$$
\lambda(l_{p}^{n})\geq c_{0}n^{\frac{1}{2}-\frac{1}{p}}\ .
$$

The proof of the theorem is complete.\ \ \ \ \ $\Box$
\sect{Proof of Theorem \ref{te110}}
{\bf A metric graph without ${\cal LE}$.}
A construction presented here may be of independent interest.
To formulate the result
we recall several notions of Graph Theory, see, e.g., [R].

Let $\Gamma=({\cal V}, {\cal E})$ be a graph with the sets of vertices
${\cal V}$ and edges ${\cal E}$. We consider below only {\em simple}
(i.e., without loops and double edges) and {\em connected} graphs. The
latter means that for every two vertices $v'$, $v''$ there is a
{\em path}\footnote{i.e., an alternative sequence $\{v_{0},e_{1},v_{1},\dots,
e_{n},v_{n}\}$ with pairwise distinct edges $e_{i}$ such that $e_{i}$ joins
$v_{i-1}$ and $v_{i}$. Vertices $v_{0}$ and $v_{n}$ are the {\em head} and
the {\em tail} of this path.} whose head and tail are $v'$ and $v''$,
respectively. The {\em distance} $d_{\Gamma}(v',v'')$ between two vertices
$v',v''$ is the {\em length} (number of edges) of a shortest path between
them. To introduce the {\em metric graph} ${\cal M}_{\Gamma}$
associated with $\Gamma$ we regard every $e\in {\cal E}$ as the unit
interval of $\Re$ and equip the 1-dimensional CW complex obtained in this
way with the {\em path} ({\em length}) metric generated by $d_{\Gamma}$.
Thus the restriction of this metric to ${\cal V}$ coincides with
$d_{\Gamma}$ and every edge is isometric to $[0,1]\subset\Re$, see, e.g.,
[BH] for further details.
\begin{Proposition}\label{pr61}
There exists a graph $\Gamma=({\cal V}, {\cal E})$ and a subset
$S\subset {\cal V}$ such that
\begin{itemize}
\item[{\rm (a)}]
the degree of every vertex\footnote{i.e., the number of edges incident to
$v$. This is denoted by $deg\ \! v$} $v\in {\cal V}$ is at most 3;
\item[{\rm (b)}]
$$
Ext(S,{\cal M}_{\Gamma})=\emptyset \ .
$$
\end{itemize}
\end{Proposition}
\begin{R}\label{rem61}
{\rm In fact, we will prove that}
\begin{equation}\label{eq61}
Ext(S,{\cal V})=\emptyset\ .
\end{equation}
{\rm Here and above $S$ and ${\cal V}$ are regarded as metric subspaces of
${\cal M}_{\Gamma}$. Since $\lambda(S,{\cal M}_{\Gamma})\geq
\lambda(S,{\cal V})$, equality (\ref{eq61}) implies statement (b).}
\end{R}
{\bf Proof.}
Our argument is based on the following result.

Let $\Z_{p}^{n}$ denote $\Z^{n}$ regarding as a metric subspace of
$l_{p}^{n}$.
\begin{Lm}\label{l61}
There is an absolute  constant $c_{0}>0$ such that for every $n$
\begin{equation}\label{e62}
\lambda(\Z_{1}^{n})\geq c_{0}\sqrt{n}\ .
\end{equation}
\end{Lm}
{\bf Proof.}
We apply Theorem \ref{te11} with $M:=l_{1}^{n}$, $S:=\Z_{1}^{n}$ and
the dilation $\delta:x\mapsto \frac{1}{2}x$. Then
$\delta(S)=\frac{1}{2}\Z^{n}\supset\Z^{n}$ and the Lipschitz constants
of $\delta$ and $\delta^{-1}$ equal $\frac{1}{2}$ and $2$, respectively.
Consequently, the assumptions of Theorem \ref{te11} hold for this case
and therefore
\begin{equation}\label{e63}
\lambda(l_{1}^{n})=\lambda(\Z_{1}^{n})\ .
\end{equation}
Since $\lambda(l_{1}^{n})\geq\lambda_{conv}(l_{1}^{n})$, it remains to
apply to (\ref{e63}) the lower estimate of Theorem \ref{te1.11} with
$p=1$.\ \ \ \ \ $\Box$
\begin{R}\label{re62}
{\rm The same argument gives}
\begin{equation}\label{e64}
\lambda(l_{p}^{n})=\lambda(\Z_{p}^{n})\geq c_{0}
n^{\left|\frac{1}{p}-\frac{1}{2}\right|}\ .
\end{equation}
\end{R}

Let now $\Z_{1}^{n}(l)$ denote the discrete cube of the lengthside
$l\in\N$, i.e.
$$
\Z_{1}^{n}(l):=\Z_{1}^{n}\cap [-l,l]^{n}\ .
$$
\begin{Lm}\label{l63}
For every $n\in\N$ there is an integer $l=l(n)$ and a subset
$S_{n}\subset \Z_{1}^{n}(l)$ such that
\begin{equation}\label{e65}
\lambda(S_{n},\Z_{1}^{n}(l))\geq c_{1}\sqrt{n}
\end{equation}
with $c_{1}>0$ independent of $n$.
\end{Lm}
{\bf Proof.}
By Corollary \ref{c12}
$$
\lambda(\Z_{1}^{n})=\sup_{F}\lambda(F)
$$
where $F$ runs through all finite subsets $F\subset\Z^{n}$. On the other
hand
$$
\lambda(\Z_{1}^{n})\geq\sup_{l\in\N}\lambda(\Z_{1}^{n}(l))\ .
$$
At last, Corollary \ref{c12} gives
\begin{equation}\label{e66}
\lambda(\Z_{1}^{n}(l))=\sup_{F\subset\Z_{1}^{n}(l)}\lambda(F)\ .
\end{equation}
These three relations imply that
$$
\lambda(\Z_{1}^{n})=\sup_{l\in\N}\lambda(\Z_{1}^{n}(l))\ .
$$
Together with (\ref{e62}) this gives for some $l=l(n)$
$$
\lambda(\Z_{1}^{n}(l(n)))> \frac{c_{0}}{2}\sqrt{n}\ .
$$
Applying now (\ref{e66}) with $l:=l(n)$ we then find
$S_{n}\subset\Z_{1}^{n}(l(n))$ such that for $l=l(n)$
$$
\lambda(S_{n},\Z_{1}^{n}(l)):=
\inf\{||E||\ :\ E\in Ext(S_{n},\Z_{1}^{n}(l))\}\geq \frac{c_{0}}{2}\sqrt{n}\ .
$$
The result is done.\ \ \ \ \ $\Box$

Let now $G_{n}:=(\Z^{n},{\cal E}^{n})$ be a graph whose set of edges is
given by
$$
{\cal E}^{n}:=\{(x,y)\ :\ x,y\in\Z^{n},\ ||x-y||_{l_{1}^{n}}=1\}\ .
$$
Let $\Gamma_{n}:=(V_{n},E_{n})$ be a subgraph of $G_{n}$ whose set of
vertices is
$$
V_{n}:=\Z^{n}\cap [-l(n),l(n)]^{n}
$$
where $l(n)$ is defined in Lemma \ref{l63}. So the set $S_{n}$
from (\ref{e65}) contains in $V_{n}$.
The metric graph
${\cal M}_{\Gamma_{n}}$ is then a (metric) subspace of the space
$l_{1}^{n}$, but it also can and will be regarded below as a
subspace of $l_{2}^{n}$ with the path metric induced by the Euclidean
metric.
\begin{Lm}\label{l64}
There is a finite connected graph $\widehat\Gamma_{n}:=(\widehat V_{n},
\widehat E_{n})$ and a subset $\widehat S_{n}\subset\widehat V_{n}$ such that
\begin{itemize}
\item[{\rm (a)}]
for every vertex $v\in\widehat V_{n}$
\begin{equation}\label{eq67}
deg\ \!v\leq 3\ ;
\end{equation}
\item[{\rm (b)}]
the underlying set of the metric graph ${\cal M}_{\widehat\Gamma_{n}}$
is a subset of the $n$-dimensional Euclidean space and its metric is
the path metric generated by the Euclidean one;
\item[{\rm (c)}]
there is an absolute constant $c>0$ such that for every $n$
\begin{equation}\label{eq68}
\lambda(\widehat S_{n},\widehat V_{n})\geq c\sqrt{n}\ ;
\end{equation}
here $\widehat S_{n}$ and $\widehat V_{n}$ are regarded as subspaces of
${\cal M}_{\widehat\Gamma_{n}}$.
\end{itemize}
\end{Lm}
{\bf Proof.} Let $\epsilon:=\frac{1}{q\sqrt{2}}$ for some natural $q\geq 2$.
For a vertex $v\in V_{n}\subset l_{2}^{n}$ of $\Gamma_{n}$,
let $S(v)$ stand for the $(n-1)$-dimensional Euclidean
sphere centered at $v$ and of radius $\frac{\epsilon}{1+2\epsilon}$.
Then $S(v)$ intersects $N(v)$ ($n\leq N(v)\leq 2n$) edges of
${\cal M}_{\Gamma_{n}}$ at some points denoted by
$p_{i}(v)$, $i=0,\dots,N(v)-1$. The ordering is choosing in such a way that
any interval $conv\{p_{i}(v),p_{i+1}(v)\}$ does not belong to
${\cal M}_{\Gamma_{n}}$ (here and below
$p_{N(v)}(v)$ is identified with $p_{0}(v)$).
Let us introduce a new graph with the set of
vertices $\{p_{i}(v)\ :\ i=0,\dots,N(v)-1,\  v\in V_{n}\}$ and
the set of edges defined as follows. This set contains edges determined by
all pairs $(p_{i}(v),p_{i+1}(v))$ with $0\leq i\leq N(v)-1$ and $v\in V_{n}$
and, moreover, all
edges formed by all pairs $(p_{i}(v'),p_{j}(v''))$ where $v',v''$ are
the head and the tail of an edge $e\in E_{n}$, and $i\neq j$
satisfy the condition
\begin{equation}\label{con0}
conv\{p_{i}(v'),p_{j}(v'')\}\subset e\ ;
\end{equation}
here $e$ is regarded as a subset (interval) of ${\cal M}_{\Gamma_{n}}$.
In this way we obtain
a new graph (and an associated metric space with the path metric
induced by the Euclidean one) whose vertices has degree at most 3.

The lengths of edges $(p_{i}(v),p_{i+1}(v))$ of this graph are
$\frac{1}{q(1+2\epsilon)}$ while the lengths of edges
$(p_{i}(v'),p_{i}(v''))$ satisfying (\ref{con0}) are
$\frac{1}{1+2\epsilon}$. Then we add new vertices (and edges) by inserting
into every edge satisfying (\ref{con0}) the $(q-1)$ equally distributed
new vertices. (Note that every new vertex obtained in this way has degree 2.)
At last, by dilation (with respect to $0\in l_{2}^{n}$)
with factor $q(1+2\epsilon)$ we obtain a new graph
$\widehat\Gamma_{n}:=(\widehat E_{n},\widehat V_{n})$ whose edges are of
length one, and such that
$$
deg\ \! v\leq 3\ ,\ \ \ v\in\widehat V_{n}\ ,\ \ \ {\rm and}\ \ \ \exists
\ v_{n}\in \widehat V_{n}\ : \ \
deg\ \! v_{n}=2\ .
$$
Moreover, the metric graph ${\cal M}_{\widehat\Gamma_{n}}$ is a (metric)
subspace of $l_{2}^{n}$ equipped with the path metric induced by that of
$l_{2}^{n}$.

It remains to introduce the required subset
$\widehat S_{n}\subset\widehat V_{n}$. To this end we define a map
$i:V_{n}\rightarrow\widehat V_{n}$ by
$$
i(v):=q(1+2\epsilon)\cdot p_{0}(v)\ ,\ \ \ v\in V_{n}\ .
$$
Recall that $p_{0}(v)$ is a point of the sphere $S(v)\subset l_{2}^{n}$.
Since our construction depends on $\epsilon$ continuously,
and $\Gamma_{n}$ is a finite graph, we clearly have for a sufficiently
small $\epsilon$
\begin{equation}\label{e69}
(q/2)(1+2\epsilon)\cdot d(v',v'')\leq\widehat d(i(v'),i(v''))\leq
2q(1+2\epsilon)\cdot d(v',v'')\ ,\ \ \ v',v''\in V_{n}\ .
\end{equation}
Here $d,\widehat d$ are the metrics of ${\cal M}_{\Gamma_{n}}$ and
${\cal M}_{\widehat\Gamma_{n}}$, respectively. Note now that the constant
$\lambda(S_{n},V_{n})$ does not change if we replace the metric $d$ by
$q(1+2\epsilon)\cdot d$. Therefore
(\ref{e69}) and (\ref{e65}) imply the estimate
$$
\lambda(i(S_{n}),i(V_{n}))\geq \frac{1}{4}\lambda(S_{n},V_{n})\geq
\frac{c_{1}}{4}\sqrt{n}\ .
$$
We set $\widehat S_{n}:=i(S_{n})\cup\{v_{n}\}$ where
$v_{n}\in\widehat V_{n}$ satisfies $deg\ \! v_{n}=2$.
Noting that $i(V_{n})\subset\widehat V_{n}$ and there is an
$L\in Ext(i(S_{n}),\widehat S_{n})$ such that $||L||\leq 2$ (cf.
the proof of Lemma \ref{l41}), we get from here
$$
\lambda(\widehat S_{n},\widehat V_{n})\geq\frac{1}{2}
\lambda(i(S_{n}),i(V_{n}))\geq\frac{c_{1}}{8}\sqrt{n}\ .
$$
This proves (\ref{eq68}) and the lemma.\ \ \ \ \ $\Box$

Let now $\widehat\Gamma_{n}=(\widehat V_{n},\widehat E_{n})$ and
$\widehat S_{n}\subset\widehat V_{n}$ be as in Lemma \ref{l64}. Then
${\cal M}_{\widehat\Gamma_{n}}$ is a subset of the space $l_{2}^{n}$
equipped with the path metric generated by the Euclidean metric. We now
identify $l_{2}^{n}$ with its isometric copy $P_{n}$, an $n$-dimensional
plane of the Hilbert space $l_{2}(\N)$ orthogonal to the line
$\{x\in l_{2}(\N)\ :\ x_{i}=0\ {\rm for}\ i>1\}$ and intersecting this line 
at the point
$v_{n}:=(n,0,\dots)$. Then ${\cal M}_{\widehat\Gamma_{n}}$ is a subset of
$P_{n}\subset l_{2}(\N)$ equipped with the path metric generated by the
metric of $l_{2}(\N)$. Using an appropriate translation we also may and
will assume that
$$
v_{n}\subset\widehat S_{n}\ \ \ {\rm and}\ \ \ deg\ \! v_{n}=2\ .
$$
Note that
$$
dist(P_{n},P_{n+1})=||v_{n}-v_{n+1}||=1
$$
and therefore the sets ${\cal M}_{\widehat\Gamma_{n}}$ are pairwise
disjoint.

Introduce now the set of vertices and edges of the required graph
$\Gamma=({\cal V},{\cal E})$ by letting
$$
{\cal V}:=\bigcup_{n\in\N}\widehat V_{n}
$$
and, moreover,
$$
{\cal E}:=\bigcup_{n\in\N}(\widehat E_{n}\cup e_{n})
$$
where $e_{n}$ denotes the new edge joining $v_{n}$ with $v_{n+1}$.

This definition and Lemma \ref{l64} imply that
$$
deg\ \! v\leq 3,\ \ \ v\in {\cal V}\ ,
$$
and so assertion (a) of Proposition \ref{pr61} holds.

Set now
$$
S:=\bigcup_{n\in\N}\widehat S_{n}\ .
$$
We claim that this $S$ satisfies (\ref{eq61}) and assertion (b) of the
proposition. If, to the contrary, there is an operator $E\in Ext(S,{\cal V})$,
we choose $n$ so that
\begin{equation}\label{e610}
c\sqrt{n}>||E||
\end{equation}
with the constant $c$ from (\ref{eq68}). Introduce for this $n$ an operator
$T_{n}$ given on \penalty-10000 $f\in Lip(\widehat S_{n})$ by
$$
(T_{n}f)(v):=
\left\{
\begin{array}{ccc}
f(v),&{\rm if}&v\in \widehat S_{n}\\
f(v_{n}),&{\rm if}&v\in S\setminus\widehat S_{n}\ .
\end{array}
\right.
$$
Show that $T_{n}$ maps $Lip(\widehat S_{n})$ into $Lip(S)$ and its norm is
$1$. To accomplish this we have to show that for $v'\in\widehat S_{n}$
and $v''\in S\setminus\widehat S_{n}$
$$
|(T_{n}f)(v')-(T_{n}f)(v'')|\leq ||f||_{Lip(\widehat S_{n})}d(v',v'')\ .
$$
But the left-hand side here is
$$
|f(v')-f(v_{n})|\leq ||f||_{Lip(\widehat S_{n})}d(v',v_{n})
$$
and $d(v',v_{n})\leq d(v',v'')$ by the definition of $S$ and the metric $d$
of ${\cal M}_{\Gamma}$. This implies the required statement for $T_{n}$.

Finally, introduce the restriction operator
$R_{n}:Lip({\cal V})\rightarrow Lip(\widehat V_{n})$ by
$$
R_{n}f=f|_{\widehat V_{n}}
$$
and set $E_{n}:=R_{n}ET_{n}$.
Then $E_{n}\in Ext(\widehat S_{n},\widehat V_{n})$ and its norm is bounded by
$||E||$. This immediately implies that
$$
\lambda(\widehat S_{n},\widehat V_{n})\leq ||E||
$$
in contradiction with (\ref{eq68}) and our choice of $n$, see (\ref{e610}).

So we establish (\ref{eq61}) and complete the proof of the
proposition.\ \ \ \ \ $\Box$\\
{\bf Two-dimensional metric space without ${\cal LE}$.}
In order to complete the proof of Theorem \ref{te110} we have
to construct a connected two-dimensional metric space $M$ of bounded
geometry so that
$$
Ext(S,M)=\emptyset
$$
for some its subspace $S$. In fact, $M$ will be a Riemannian manifold
(with the geodesic (inner) metric). This will be done by
sewing surfaces of three types along the metric graph ${\cal M}_{\Gamma}$
of the previous part. At the first stage we introduce an open cover of
${\cal M}_{\Gamma}$ by balls and a related coordinate system and
partition of unity. To simplify evaluation we replace the metric of
${\cal M}_{\Gamma}$ by $\widetilde d_{\Gamma}:=4d_{\Gamma}$. So every
edge $e\subset {\cal M}_{\Gamma}$ is a closed interval of length $4$.
(Note that the abstract graph $\Gamma=({\cal V}, {\cal E})$ remains
unchanged.) The required cover $\{B(v)\}_{v\in {\cal V}}$ is given by
\begin{equation}\label{e611}
B(v):=\{m\in {\cal M}_{\Gamma}\ :\ \widetilde d(m,v)< 3\}\ .
\end{equation}
So $B(v)$ is the union of at most three intervals of length $3$ each of which
has a form $e\cap B(v)$ where every $e$ belongs to the set of edges
${\cal E}(v)$ incident to $v$. We numerate these intervals by numbers
from the set $\omega\subset\{1,2,3\}$ where $\omega=\{1\}$, $\{1,2\}$
or $\{1,2,3\}$, if $deg\ \!v=1$, $2$ or $3$, respectively. This set of
indices will be denoted by $\omega(v)$ and $i(e,v)$ (briefly, $i(e)$)
will stand for the number of $e\cap B(v)$ in this numeration.

We then introduce a {\em coordinate system}
$\psi_{v}: B(v)\rightarrow\Re^{3}$, $v\in {\cal V}$, of ${\cal M}_{\Gamma}$
as follows. Let $\{b_{1},b_{2},b_{3}\}$ be the standard basis in
$\Re^{3} (:=l_{2}^{3})$. We define $\psi_{v}$ as the isometry sending
$v$ to $0$ and each interval $e\cap B(v)$, $e\in {\cal E}(v)$, to the
interval $\{tb_{i}\ :\ 0\leq t< 3\}$ of the $x_{i}$-axis with $i:=i(e)$.

Now we introduce the desired {\em partition} of unity
$\{\rho_{v}\}_{v\in {\cal V}}$ subordinate to the cover
$\{B(v)\}_{v\in {\cal V}}$. To this end one first considers a function
$\widetilde\rho_{v}:\psi_{v}(B(v))\rightarrow [0,1]$ with support strictly
inside of its domain such that $\widetilde\rho_{v}=1$ in the neighbourhood
$\cup_{e\in {\cal E}(v)}\{tb_{i(e)}\ :\ 0\leq t\leq 1\}$ of $0$ and is
$C^{\infty}$-smooth outside $0$. This function gives rise to a function
$\widehat\rho_{v}:{\cal M}_{\Gamma}\rightarrow [0,1]$ equals
$\widetilde\rho_{v}\circ\psi_{v}$ on $B(v)$ and $0$ outside. It is important
to note that there exist only three types of the functions
$\widetilde\rho_{v}$ corresponding to the types of the balls $B(v)$. Finally
we determine the required partition of unity by setting
\begin{equation}\label{e612}
\rho_{v}:=\widehat\rho_{v}/\sum_{v}\widehat\rho_{v}\ ,\ \ \ v\in {\cal V}\ .
\end{equation}

At the second stage we introduce the building blocks of our construction,
$C^{\infty}$-smooth surfaces $\Sigma_{\{1\}}$, $\Sigma_{\{1,2\}}$ and
$\Sigma_{\{1,2,3\}}$ embedded in $\Re^{3}$. We begin with a
$C^{\infty}$-function $f:[-1,3)\rightarrow [0,1]$ given by
$$
f(t):=\left\{
\begin{array}{ccc}
\sqrt{1-t^{2}},&{\rm if}&-1\leq t\leq 1-\epsilon:=\frac{3}{4}\\
\\
\sqrt{1-10\epsilon^{2}},&{\rm if}&1\leq t< 3
\end{array}
\right.
$$
In the remaining interval $[1-\epsilon,1] (:=[3/4,1])$ $f$ is an arbitrary
decreasing function smoothly joining the given endpoint values. Then we
introduce $\Sigma_{\{1\}}$ as a surface of revolution
\begin{equation}\label{e613}
\Sigma_{\{1\}}:=\{(t,f(t)\cos\theta ,f(t)\sin\theta)\in \Re^{3}\ :\
-1\leq t< 3\ ,\ 0\leq\theta<2\pi\}\ ,
\end{equation}
the result of rotating the graph of $f$ about the $x_{1}$-axis. By the
definition of $f$ this surface is the union of the unit sphere
$S^{2}\subset\Re^{3}$ with the spherical hole $S(b_{1})$ centered at $b_{1}$
and of the curvilinear (near the bottom) cylinder $T(b_{1})$ attached to
the circle $\partial S(b_{1})$ (of radius $\sqrt{1-(1-\epsilon)^{2}}$).
In turn, $T(b_{1})$ is the union of the curvilinear cylinder and that of
circular. The latter, denoted by $\widehat T(b_{1})$, is of height $2$.

Similarly $\Sigma_{\{1,2\}}$ and $\Sigma_{\{1,2,3\}}$ are the unions of the
unit sphere $S^{2}$ with the holes $S(b_{i})$ and of the cylinders
$T(b_{i})$ attached to $\partial S(b_{i})$; here $i=1,2$ or
$i=1,2,3$, respectively. Note that each $T(b_{i})$ with $i\neq 1$ is obtained
from $T(b_{1})$ by a fixed turn around the $x_{j}$-axis, $j\neq 1,i$.
This determines the isometry
\begin{equation}\label{e614}
J_{i}:\widehat T(b_{i})\rightarrow\widehat T(b_{1})\ ,\ \ \ i=1,2,3\ ,
\end{equation}
where $J_{1}$ stands for the identity map.

Using these blocks and the previous notations for $B(v)$ we now assign
to every $v\in {\cal V}$ the smooth surface
\begin{equation}\label{e615}
\Sigma(v):=\Sigma_{\omega(v)}\subset\Re^{3}\ .
\end{equation}
We denote by $S(v)\subset\Sigma_{\omega(v)}$ the corresponding sphere $S^{2}$
with holes $\{S(b_{i})\}$, $i\in\omega(v)$, and by $T(e)$,
$e\in {\cal E}(v)$, the corresponding curvilinear cylinder ($=T(b_{i(e)})$).
The circular part of the latter is denoted by $\widehat T(e)$ and the
corresponding isometry of $\widehat T(e)$ onto $\widehat T(b_{1})$ is
denoted by
$J_{e}$ ($=J_{i(e)}$)\footnote{Since $e$ belongs to two different sets, say,
${\cal E}(v)$ and ${\cal E}(v')$, we will also write $T(e,v)$,
$\widehat T(e,v)$ and $J_{e,v}$ to distinguish them from the corresponding
objects determined by $e$ as an element of ${\cal E}(v')$.}. Finally,
we equip $\Sigma(v)$ with the Riemannian metric induced by the canonical
Riemannian structure of $\Re^{3}$, and denote the corresponding geodesic
metric by $d_{v}$.

According to our construction there exists for every $\Sigma(v)$ a continuous
surjection $p_{v}:\Sigma(v)\rightarrow\psi_{v}(B(v))$ such that the
restriction of $p_{v}$ to every cylinder $\widehat T(e,v)$,
$e\in {\cal E}(v)$, is the orthogonal projection onto its axis
$I_{e}:=\{tb_{i(e,v)}\ :\ 1< t< 3\}$. Using this and the polar coordinate
$\theta$ from (\ref{e613}) we then equip each $x\in\widehat T(e,v)$ with
the cylindrical coordinates:
$$
r(x):=\psi_{v}^{-1}(p_{v}(x))\ ,\ \ \ \theta(x):=\theta(J_{e,v}(x))\ .
$$

Now we define the required smooth surface $M$ as the quotient of the
disjoint union $\sqcup_{v\in {\cal V}}\Sigma(v)$ by the equivalence relation:
$$
\begin{array}{c}
\displaystyle
x\sim y\ ,\ \ {\rm if}\ \ \ x\in\widehat T(e,v_{0})\ ,\
y\in\widehat T(e,v_{1})\ \ \ {\rm for\ some}\ \ \
e\in {\cal E}(v_{0})\cap {\cal E}(v_{1})\ \ \ {\rm and}\\
\\
\displaystyle
(r(x),\theta(x))=(r(y),\theta(y))\ .
\end{array}
$$

Let $\pi:\sqcup_{v\in {\cal V}}\Sigma(v)\rightarrow M$ be the quotient
projection. Then $\{\pi(\Sigma(v))\}_{v\in {\cal V}}$ is an open cover of
$M$. Using the partition of unity (\ref{e612}) we now introduce a partition
of unity subordinate to this cover as follows. Define a function
$\widehat\phi_{v}:\Sigma(v)\rightarrow [0,1]$ as a pullback of the
function $\rho_{v}:B(v)\rightarrow [0,1]$ given by
\begin{equation}\label{e616}
\widehat\phi_{v}:=\rho_{v}(\psi_{v}^{-1}(p_{v}(x)))\ ,\ \ \ x\in\Sigma(v).
\end{equation}
By the definitions of all functions used here, the function
$\widehat\phi_{v}$ is $C^{\infty}$-smooth in every $\widehat T(e,v)$ and is
equal to 1 outside $\cup_{e\in {\cal E}(v)}\widehat T(e,v)$. In particular,
$\widehat\phi_{v}$ is $C^{\infty}$-smooth and its support is strictly inside
$\Sigma(v)$. Since $\pi|_{\Sigma(v)}$ is a smooth embedding, the function
$\phi_{v}:M\rightarrow [0,1]$ equals $\widehat\phi_{v}\circ\pi$ on
$\pi(\Sigma(v))$ and $0$ outside is $C^{\infty}$-smooth. By (\ref{e616})
the family $\{\phi_{v}\}_{v\in {\cal V}}$ forms the required partition of
unity subordinate to the cover $\{\pi(\Sigma(v))\}_{v\in {\cal V}}$.

Using this we now determine a Riemannian metric tensor $R$ of $M$ by
$$
R:=\sum_{v\in{\cal V}}\phi_{v}\cdot(\pi^{-1})^{*}(R_{v})
$$
where $R_{v}$ stands for the metric tensor of $\Sigma(v)$. If now $d$ is
the geodesic (inner) metric of $M$ determined by $R$, then the metric space
$(M,d)$ is clearly of bounded geometry, because in this construction we used
the objects of only three different types.

It remains to find a subspace $\widetilde S$ of $(M,d)$ such that
\begin{equation}\label{e617}
Ext(M,\widetilde S)=\emptyset\ .
\end{equation}
To this end we first consider two the hole spheres $\pi(S(v_{i}))\subset M$,
$i=1,2$, such that $v_{1}$ and $v_{2}$ are joined by an edge. Let
$m_{i}\in\pi(S(v_{i})$ be arbitrary points, $i=1,2$. Then by the
definition of the metric $d$ and by a compactness argument
\begin{equation}\label{e618}
0<c\leq d(m_{1},m_{2})\leq C
\end{equation}
where $c,C$ are independent of $m_{i}$ and $v_{i}$. On the other hand in the
space ${\cal M}_{\Gamma}$
\begin{equation}\label{e619}
d_{\Gamma}(v_{1},v_{2})=1
\end{equation}
for this choice of $v_{i}$.

Let now $m_{a}\in\pi(S(v_{a}))\subset M$, $a\in\{A,B\}$, be arbitrary points
and $v_{A},v_{B}$ are distinct and
may not necessarily be joined by an edge. Let
$\{v_{i}\}_{i=1}^{n}$ be a path in the graph $\Gamma$ joining $v_{A}$ and
$v_{B}$ (here $v_{1}:=v_{A}$ and $v_{n}:=v_{B}$) such that
$$
d_{\Gamma}(v_{A},v_{B})=\sum_{i=1}^{n-1}d_{\Gamma}(v_{i},v_{i+1})\ .
$$
Together with (\ref{e618}) and (\ref{e619}) this implies that
\begin{equation}\label{e620}
cd_{\Gamma}(v_{A},v_{B})\leq\sum_{i=1}^{n-1}d(m_{i},m_{i+1})\leq
Cd_{\Gamma}(v_{A},v_{B})\ .
\end{equation}
On the other hand, the definitions of $M$ and $d$ get
\begin{equation}\label{e620a}
\widetilde c\cdot\sum_{i=1}^{n-1}d(m_{i},m_{i+1})\leq d(m_{A},m_{B})\leq
\sum_{i=1}^{n-1}d(m_{i},m_{i+1})\
\end{equation}
with some $\widetilde c>0$ independent of $m_{i}$'s.

Introduce now a map $T:{\cal V}\rightarrow M$ sending a point
$v\in {\cal V}$ to an arbitrary point $T(v)\in\pi(S(v))$. Because
of (\ref{e620}) and (\ref{e620a}) $T$ is a quasi-isometric
embedding of ${\cal V}\subset {\cal M}_{\Gamma}$ into $M$. We then
define the required subset $\widetilde S$ as the image under $T$
of the set $S\subset {\cal V}$ for which $Ext(S,{\cal
V})=\emptyset$, see (\ref{eq61}). Then we have for $\widetilde
S:=T(S)$
$$
Ext(\widetilde S,M)=\emptyset\ .
$$

The proof of Theorem \ref{te110} is complete.\ \ \ \ \ $\Box$
\begin{R}\label{nash}
{\rm Using the Nash embedding theorem one can realize the
Riemannian manifold $M$ as a $C^{\infty}$-surface in an open ball
of $\Re^{3}$ (with the Riemannian quadratic form induced from the
canonical Riemannian structure of $\Re^{3}$).}
\end{R}
\sect{Proofs of Theorem \ref{te17} and Corollaries \ref{2.12''} and 
\ref{c111b}} 
{\bf Proof of Theorem \ref{te17}.} Let $\Gamma$ be the $R$-lattice and
${\cal B}:=\{B_{R}(\gamma)\}_{\gamma\in\Gamma}$. By the definition
of an $R$-lattice,  ${\cal B}$ and $\frac{1}{2}{\cal
B}:=\{B_{R/2}(\gamma)\}_{\gamma\in\Gamma}$ are covers of $M$.
\begin{Lm}\label{l71}
Multiplicity of ${\cal B}$ is bounded by a constant $\mu$ depending only on
$c_{\Gamma}$ and $N=N_{M}$.
\end{Lm}

For convenience of the reader we outline the proof of this well-known
fact.

Let $B_{R}(\gamma_{i})$, $1\leq i\leq k$, contain a point $m$. Then
all $\gamma_{i}$ are in the ball $B_{R}(m)$. Since
$d(\gamma_{i},\gamma_{j})>2cR$, $i\neq j$, see Definition \ref{def2}, any
cover of $B_{R}(m)$ by balls of radius $cR$ separates $\gamma_{i}$,
i.e., distinct $\gamma_{i}$ lie in the distinct balls. Hence cardinality
of such a cover is at least $k$.
On the other hand the doubling condition
implies that there is a cover of $B_{R}(m)$ by balls of radius $cR$ and
cardinality $N^{s}$ where $s:=[\log_{2}\frac{1}{c_{\Gamma}}]+1$.
So multiplicity of ${\cal B}$ is bounded by $N^{s}$.\ \ \ \ \ $\Box$
\begin{Lm}\label{l72}
There is a partition of unity $\{\rho_{\gamma}\}_{\gamma\in\Gamma}$
subordinate to ${\cal B}$ such that
\begin{equation}\label{e71}
K:=\sup_{\gamma}||\rho_{\gamma}||_{Lip(M)}<\infty
\end{equation}
where $K$ depends only on $c_{\Gamma}$, $N=N_{M}$ and $R=R_{M}$.
\end{Lm}
{\bf Proof.}
Set
$$
B_{\gamma}:=B_{R}(\gamma)\ \ \ {\rm and}\ \ \ ^{c}B_{\gamma}:=
M\setminus B_{\gamma}
$$
and define
$$
d_{\gamma}(m):=dist(m,\ \! ^{c}B_{\gamma})\ ,\ \ \ m\in M\ .
$$
It is clear that
\begin{equation}\label{e72}
supp\ d_{\gamma}\subset B_{\gamma}\ \ \  {\rm and}\ \ \
||d_{\gamma}||_{Lip(M)}\leq 1\ .
\end{equation}
Let now $\phi:\Re_{+}\rightarrow [0,1]$ be continuous, equal one on $[0,R/2]$,
zero on $[R,\infty)$ and linear on $[R/2,R]$. Introduce the function
\begin{equation}\label{e73}
s:=\sum_{\gamma}\phi\circ d_{\gamma}\ .
\end{equation}
By Lemma \ref{l71} only at most $\mu$ terms here are nonzero at every point.
Therefore
$$
||s||_{Lip(M)}\leq 2\mu ||\phi||_{Lip(\Re_{+})}
\sup_{\gamma}||d_{\gamma}||_{Lip(M)}
$$
and by (\ref{e72}) and the definition of $\phi$ we get
\begin{equation}\label{e74}
||s||_{Lip(M)}\leq 4\mu/R\ .
\end{equation}
On the other hand, every $m\in M$ is contained in some ball
$B_{R/2}(\gamma)$ of the cover $\frac{1}{2}{\cal B}$. For this $\gamma$
$$
(\phi\circ d_{\gamma})(m)\geq\phi(R/2)=1
$$
and therefore
\begin{equation}\label{e75}
s\geq 1\ .
\end{equation}
Introduce now the required partition by
$$
\rho_{\gamma}:=\frac{\phi\circ d_{\gamma}}{s}\ ,\ \ \ \gamma\in\Gamma\ .
$$
Then $\{\rho_{\gamma}\}$ is clearly a partition of unity subordinate to
${\cal B}$. Moreover, we have
$$
|\rho_{\gamma}(m)-\rho_{\gamma}(m')|\leq
\frac{|\phi(d_{\gamma}(m))-\phi(d_{\gamma}(m'))|}{s(m)}+
\frac{\phi(d_{\gamma}(m'))}{s(m)\cdot s(m')}\cdot |s(m)-s(m')|
$$
and application of (\ref{e75}), (\ref{e74}) and (\ref{e72}) leads to
the desired inequality
$$
||\rho_{\gamma}||_{Lip(M)}\leq\frac{2}{R}(2\mu+1)\ .\ \ \ \ \ \Box
$$
\begin{Lm}\label{l73}
$$
Ext(\Gamma,M)\neq\emptyset\ .
$$
\end{Lm}
{\bf Proof.}
By the assumption (\ref{e16}) of the theorem,
for every $\gamma\in\Gamma$ there
is a linear operator $E_{\gamma}\in Ext(\Gamma\cap B_{\gamma},B_{\gamma})$
such that
\begin{equation}\label{e76}
||E_{\gamma}||\leq\lambda_{R}\ ,\ \ \ \gamma\in\Gamma\ .
\end{equation}
Using this we introduce the required linear operator by
\begin{equation}\label{e77}
Ef:=\sum_{\gamma\in\Gamma}(E_{\gamma}f_{\gamma})\rho_{\gamma}\ ,\ \ \
f\in Lip(\Gamma)\ ,
\end{equation}
where $\{\rho_{\gamma}\}$ is the partition of unity from Lemma \ref{l72} and
$f_{\gamma}:=f|_{\Gamma\cap B_{\gamma}}$; here we assume that
$E_{\gamma}f_{\gamma}$ is zero outside of $B_{\gamma}$.
We have to show that
\begin{equation}\label{e78}
Ef|_{\Gamma}=f|_{\Gamma}
\end{equation}
and estimate $||Ef||_{Lip(M)}$.

Given $\hat\gamma\in\Gamma$ we can write
$$
(Ef)(\hat\gamma)=\sum_{B_{\gamma}\ni\hat\gamma}
(E_{\gamma}f_{\gamma})(\hat\gamma)\rho_{\gamma}(\hat\gamma)\ .
$$
Since $E_{\gamma}$ is an extension from $B_{\gamma}\cap\Gamma$ we get
$$
(E_{\gamma}f_{\gamma})(\hat\gamma)=f_{\gamma}(\hat\gamma)=
f(\hat\gamma)\ .
$$
Moreover,
$\displaystyle\sum_{B_{\gamma}\ni\hat\gamma}\rho_{\gamma}(\hat\gamma)=1$,
and (\ref{e78}) is done.

To estimate the Lipschitz constant of $Ef$, we extend $E_{\gamma}f_{\gamma}$
outside of $B_{\gamma}$ so that the (non-linear)
extension $F_{\gamma}$ satisfies
\begin{equation}\label{e79}
||F_{\gamma}||_{Lip(M)}=||E_{\gamma}f_{\gamma}||_{Lip(B_{\gamma})}\ .
\end{equation}
Since $\rho_{\gamma}F_{\gamma}=\rho_{\gamma}E_{\gamma}f_{\gamma}$, we have
\begin{equation}\label{e710}
Ef=\sum_{\gamma}F_{\gamma}\rho_{\gamma}\ .
\end{equation}
Given $\hat\gamma\in\Gamma$ introduce a function $G_{\hat\gamma}$ by
\begin{equation}\label{e711}
G_{\hat\gamma}:=\sum_{\gamma}(F_{\gamma}-F_{\hat\gamma})\rho_{\gamma}:=
\sum_{\gamma}F_{\gamma\hat\gamma}\rho_{\gamma}\ .
\end{equation}
Then we can write for every $\hat\gamma$
\begin{equation}\label{e712}
Ef=F_{\hat\gamma}+G_{\hat\gamma}\ .
\end{equation}
It follows from (\ref{e76}) and (\ref{e79}) that
\begin{equation}\label{e713}
||F_{\hat\gamma}||_{Lip(M)}\leq\lambda_{R}||f||_{Lip(\Gamma)}\ .
\end{equation}
Prove now that
\begin{equation}\label{e714}
|F_{\gamma\hat\gamma}(m)|\leq 4R\lambda_{R}||f||_{Lip(\Gamma)}\ ,\ \ \
m\in B_{\gamma}\cap B_{\hat\gamma}\ .
\end{equation}
In fact, we have for these $m$
$$
\begin{array}{c}
|F_{\gamma\hat\gamma}(m)|=|(E_{\gamma}f_{\gamma}-E_{\hat\gamma}f_{\hat\gamma})
(m)|\leq |f(\gamma)-f(\hat\gamma)|+|(E_{\gamma}f_{\gamma})(m)-
(E_{\gamma}f_{\gamma})(\gamma)|+\\
\\
|(E_{\hat\gamma}f_{\hat\gamma})(m)-
(E_{\hat\gamma}f_{\hat\gamma})(\hat\gamma)|\ .
\end{array}
$$
Using now (\ref{e76}) to estimate the right-hand side we get
$$
|F_{\gamma\hat\gamma}(m)|\leq\lambda_{R}||f||_{Lip(\Gamma)}
(d(\gamma,\hat\gamma)+d(m,\gamma)+d(m,\hat\gamma))\leq 4R\lambda_{R}
||f||_{Lip(\Gamma)}\ .
$$
We apply this to estimate
$$
\Delta G_{\hat\gamma}:=G_{\hat\gamma}(m)-G_{\hat\gamma}(m')
$$
provided that $m,m'\in B_{\hat\gamma}$. We get
$$
|\Delta G_{\hat\gamma}|\leq\sum_{B_{\gamma}\cap B_{\hat\gamma}\ni m}
|\Delta\rho_{\gamma}|\cdot |F_{\gamma\hat\gamma}(m)|+
\sum_{B_{\gamma}\cap B_{\hat\gamma}\ni m'}\rho_{\gamma}(m')\cdot
|\Delta F_{\gamma\hat\gamma}|
$$
(here $\Delta\rho_{\gamma}$ and $\Delta F_{\gamma\hat\gamma}$ are defined
similarly to $\Delta G_{\hat\gamma}$).
The first sum is estimated by (\ref{e714}), (\ref{e71}) and
Lemma \ref{l71}, while the second sum is at most
$2\lambda_{R}||f||_{Lip(\Gamma)}d(m,m')$ by (\ref{e79}) and
(\ref{e76}). This leads to the estimate
$$
|\Delta G_{\hat\gamma}|\leq (8RK\mu+2)\cdot\lambda_{R}
\cdot ||f||_{Lip(\Gamma)}d(m,m')\ ,\ \ \ m,m'\in B_{\hat\gamma}\ .
$$
Together with (\ref{e713}) this gives for these
$m,m'$:
\begin{equation}\label{e715}
|(Ef)(m)-(Ef)(m')|\leq C||f||_{Lip(\Gamma)}d(m,m')\ .
\end{equation}
Here and below $C$ denotes a constant depending only on the basic parameters,
that may change from line to line.

It remains to prove (\ref{e715}) for $m,m'$ belonging to distinct balls
$B_{\gamma}$. Let $m\in B_{\gamma}$ and $m'$ be a point of some
$B_{R/2}(\hat\gamma)$ from the cover $\frac{1}{2}{\cal B}$. Then
$m\in B_{\gamma}\setminus B_{\hat\gamma}$ and so
\begin{equation}\label{e716}
d(m,m')\geq R/2\ .
\end{equation}
Using now (\ref{e712}) we have
$$
|(Ef)(m)-(Ef)(m')|\leq |F_{\gamma}(m)-F_{\hat\gamma}(m')|+
|G_{\gamma}(m)-G_{\hat\gamma}(m')|:=I_{1}+I_{2}\ .
$$
By the definition of $F_{\gamma}$, we then get
$$
I_{1}\leq |f(\gamma)-f(\hat\gamma)|+
|(E_{\gamma}f_{\gamma})(m)-(E_{\gamma}f_{\gamma})(\gamma)|+
|(E_{\hat\gamma}f_{\hat\gamma})(m')-
(E_{\hat\gamma}f_{\hat\gamma})(\hat\gamma)|\ .
$$
Together with (\ref{e76}) this leads to the estimate
$$
\begin{array}{c}
I_{1}\leq\lambda_{R}||f||_{Lip(\Gamma)}(d(\gamma,\hat\gamma)+
d(m,\gamma)+d(m',\hat\gamma))\leq\\
\\
 2\lambda_{R}||f||_{Lip(\Gamma)}
(d(m,m')+d(m,\gamma)+d(m',\hat\gamma))\ .
\end{array}
$$
Since $d(m,\gamma)+d(m',\hat\gamma)\leq 2R\leq 4d(m,m')$ by (\ref{e716}),
we therefore have
\begin{equation}\label{e717}
I_{1}\leq C||f||_{Lip(\Gamma)}d(m,m')\ .
\end{equation}
To estimate $I_{2}$, note that for $m\in B_{\gamma}$
$$
G_{\gamma}(m)=\sum_{B_{\gamma'}\cap B_{\gamma}\ni m}
(\rho_{\gamma'}F_{\gamma'\gamma})(m)\ .
$$
In combination with (\ref{e714}) and (\ref{e716}) this gives
$$
|G_{\gamma}(m)|\leq 4\lambda_{R}R||f||_{Lip(\Gamma)}\leq
C||f||_{Lip(\Gamma)}d(m,m')\ .
$$
The same argument estimates $G_{\hat\gamma}(m')$ for $m'\in B_{\hat\gamma}$.
Hence
$$
I_{2}\leq |G_{\gamma}(m)|+|G_{\hat\gamma}(m')|\leq C||f||_{Lip(\Gamma)}d(m,m')
\ .
$$
Together with (\ref{e717}) and (\ref{e715}) this leads to the inequality
$$
||Ef||_{Lip(M)}\leq C||f||_{Lip(\Gamma)}\ .
$$
Hence $E$ is an operator from $Ext(\Gamma,M)$.\ \ \ \ \ $\Box$

We now in a position to complete the proof of Theorem \ref{te17}.
According to Theorem \ref{te11} we have to show that
\begin{equation}\label{e718}
\sup_{F}\lambda(F)<\infty
\end{equation}
where $F$ runs through all finite point subspaces of $M$. To this end consider
a \penalty-10000 ``$\Gamma$-envelope'' of such $F$ given by
$$
\widehat F:=\{\gamma\in\Gamma\ :\ B_{\gamma}\cap F\neq\emptyset\}\ .
$$
Then $\{B_{\gamma}\ :\ \gamma\in\widehat F\}\subset {\cal B}$ is an open
cover of $F$. By assumption (\ref{e16}) of the theorem for every
$\gamma\in\widehat F$ there is an operator
$E_{\gamma}\in Ext(F\cap B_{\gamma},B_{\gamma})$ such that
$$
||E_{\gamma}||\leq\lambda_{R}\ .
$$
Introduce now a linear operator $T$ given on
$f\in Lip(F)$ by
$$
(Tf)(\gamma):=(E_{\gamma}f_{\gamma})(\gamma)\ ,\ \ \ \ \gamma\in\widehat F
$$
where $f_{\gamma}:=f|_{B_{\gamma}\cap F}$.
Show that
\begin{equation}\label{e719}
T: Lip(F)\rightarrow Lip(\widehat F)\ \ \ {\rm and}\ \ \
||T||\leq\lambda_{R}(2/c_{\Gamma}+1)\ .
\end{equation}
Actually, let $\gamma_{i}\in\widehat F$ and $m_{i}\in B_{\gamma_{i}}\cap F$,
$i=1,2$. Then $(E_{\gamma_{i}}f_{\gamma_{i}})(m_{i})=f(m_{i})$ and
$$
\begin{array}{c}
|(Tf)(\gamma_{1})-(Tf)(\gamma_{2})|\leq\sum_{i=1,2}
|(E_{\gamma_{i}}f_{\gamma_{i}})(\gamma_{i})-
(E_{\gamma_{i}}f_{\gamma_{i}})(m_{i})|+|f(m_{1})-f(m_{2})|\leq\\
\\
\lambda_{R}||f||_{Lip(F)}(d(\gamma_{1},m_{1})+d(\gamma_{2},m_{2})+
d(m_{1},m_{2}))\ .
\end{array}
$$
The sum in the brackets does not exceed $2R+d(m_{1},m_{2})\leq
4R+d(\gamma_{1},\gamma_{2})$. Moreover, by the definition of an $R$-lattice,
$d(\gamma_{1},\gamma_{2})\geq 2c_{\Gamma}R$. Combining these estimates to have
$$
|(Tf)(\gamma_{1})-(Tf)(\gamma_{2})|\leq\lambda_{R}||f||_{Lip(F)}
(2/c_{\Gamma}+1)d(\gamma_{1},\gamma_{2})\ .
$$

This establishes (\ref{e719}).

Now the assumption (\ref{e16}) of the theorem implies that there is an
operator $L$ from $Ext(\widehat F,\Gamma)$ whose norm is bounded by
$\lambda_{\Gamma}$. Composing $T$ and $L$ with the operator
$E\in Ext(\Gamma,M)$ of Lemma \ref{l73} we obtain the operator
\begin{equation}\label{e720}
\widetilde E:=ELT: Lip(F)\rightarrow Lip(M)
\end{equation}
whose norm is bounded by a constant depending only on $\lambda(\Gamma)$,
$\lambda_{R}$, $c_{\Gamma}$, $R$ and $N$. This definition also implies that
\begin{equation}\label{e721}
(\widetilde Ef)(\gamma):=(E_{\gamma}f_{\gamma})(\gamma)\ ,\ \ \
\gamma\in\widehat F\ .
\end{equation}
Unfortunately, $\widetilde E$ is not extension from $F$ and we modify it to
obtain the required extension operator. To accomplish this we, first,
introduce an operator $\widehat T$ given on $f\in Lip(F)$ by
\begin{equation}\label{e722}
(\widehat Tf)(m):=\left\{
\begin{array}{ccc}
(\widetilde Ef)(m),&{\rm if}&m\in\widehat F\\
f(m),&{\rm if}&m\in F\setminus\widehat F
\end{array}
\right.
\end{equation}
\begin{Lm}\label{l74}
$\widehat T:Lip(F)\rightarrow Lip(F\cup\widehat F)$\ \ and\ \
$||\widehat T||\leq C$.
\end{Lm}
{\bf Proof.}
It clearly suffices to estimate
$$
I:=|(\widehat Tf)(m_{1})-(\widehat Tf)(m_{2})|
$$
for $m_{1}\in\widehat F$ and $m_{2}\in F\setminus\widehat F$. Let, first,
these points belong to a ball $B_{\gamma}$ (hence $m_{1}=\gamma$). Then
(\ref{e721}) and the implication $E_{\gamma}\in
Ext(F\cap B_{\gamma},B_{\gamma})$ imply
$$
\begin{array}{c}
I=|(E_{\gamma}f_{\gamma})(\gamma)-f(m_{2})|=
|(E_{\gamma}f_{\gamma})(\gamma)-(E_{\gamma}f_{\gamma})(m_{2})|\leq\\
\\
\lambda_{R}||f||_{Lip(F)}d(\gamma,m_{2}):=\lambda_{R}||f||_{Lip(F)}
d(m_{1},m_{2})\ .
\end{array}
$$
In the remaining case $m_{2}\in B_{\hat\gamma}\setminus B_{\gamma}$ for some
$\hat\gamma\in\widehat F$ and therefore
\begin{equation}\label{e723}
d(m_{1},m_{2})=d(\gamma,m_{2})\geq R\ .
\end{equation}
Similarly to the previous estimate we now get
$$
\begin{array}{c}
I\leq |(\widetilde Ef)(\gamma)-(\widetilde Ef)(\hat\gamma)|+
|(E_{\hat\gamma}f_{\hat\gamma})(\hat\gamma)-
(E_{\hat\gamma}f_{\hat\gamma})(m_{2})|\leq\\
\\
(||\widetilde E||d(\gamma,\hat\gamma)+
\lambda_{R}d(\hat\gamma,m_{2}))||f||_{Lip(F)}\ .
\end{array}
$$
Moreover, (\ref{e723}) implies
$$
d(\gamma,\hat\gamma)+d(\hat\gamma,m_{2})\leq d(\gamma,m_{2})+
2d(\hat\gamma,m_{2})\leq d(\gamma,m_{2})+2R\leq 3d(\gamma,m_{2})=
3d(m_{1},m_{2})\ .
$$
Together with the previous inequalities this gives the required estimate
of $I$.
\ \ \ \ \ $\Box$

The next operator that will be used in our construction is defined
on \penalty-10000 $f\in Lip(F)$ by
\begin{equation}\label{e724}
(\widehat Sf)(m):=(\widehat Tf)(m)-(\widetilde Ef)(m)\ ,\ \ \ m\in
F\cup\widehat F\ .
\end{equation}
\begin{Lm}\label{l75}
$||\widehat Sf||_{l_{\infty}(F\cup\widehat F)}\leq C||f||_{Lip(F)}$ and,
moreover,
$$
\widehat S:Lip(F)\rightarrow Lip(F\cup\widehat F)\ \ \ {\rm and}\ \ \
||\widehat S||\leq C\ .
$$
\end{Lm}
{\bf Proof.}
The second statement follows straightforwardly from (\ref{e724}).
If, now, $m\in B_{\gamma}\cap (F\cup\widehat F)$, then by the same
definition
$$
(\widehat Sf)(m)=[(\widehat Tf)(m)-(\widehat Tf)(\gamma)]+
[(\widetilde Ef)(\gamma)-(\widetilde Ef)(m)]
$$
which implies that
$$
|(\widehat Sf)(m)|\leq C||f||_{Lip(F)}d(m,\gamma)\leq CR||f||_{Lip(F)}\ .
\ \ \ \ \ \Box
$$

Finally we introduce an operator $\widehat K$ given on
$g\in (Lip\cap l_{\infty})(F\cup\widehat F)$ by
\begin{equation}\label{e725}
\widehat Kg:=\sum_{\gamma}(E_{\gamma}g_{\gamma})\rho_{\gamma}\ ;
\end{equation}
here $\gamma$ runs through the set $\{\gamma\in\Gamma\ :\
(F\cup\widehat F)\cap B_{\gamma}\neq\emptyset\}$, and $\{\rho_{\gamma}\}$
is the partition of unity of Lemma \ref{l72}. Besides,
$E_{\gamma}$ is an operator from
$Ext((F\cup \widehat F)\cap B_{\gamma}, B_{\gamma})$ with
\begin{equation}\label{e726}
||E_{\gamma}||\leq\lambda_{R}\ ,
\end{equation}
and $g_{\gamma}:=g|_{(F\cup\widehat F)\cap B_{\gamma}}$.
\begin{Lm}\label{l76}
$$
||\widehat Kg||_{Lip(M)}\leq C||g||_{(Lip\cap l_{\infty})(F\cup\widehat F)}\ .
$$
\end{Lm}
{\bf Proof.}
As in the proof of Lemma \ref{l73} it is convenient to extend
every $E_{\gamma}g_{\gamma}$ outside $B_{\gamma}$ so that the
extension $F_{\gamma}$ satisfies
\begin{equation}\label{e726a}
||F_{\gamma}||_{Lip(M)}=||E_{\gamma}g_{\gamma}||_{Lip(B_{\gamma})}\ .
\end{equation}
Then we clearly have
\begin{equation}\label{e725b}
\widehat Kg=\sum_{\gamma}F_{\gamma}\rho_{\gamma}\ .
\end{equation}

Now, according to (\ref{e725}) we get for $m\in F\cup\widehat F$
$$
(\widehat Kg)(m)=\sum\rho_{\gamma}(m)g(m)=g(m)\ .
$$
So it remains to estimate the right-hand side of the inequality
$$
\begin{array}{c}
I:=|(\widehat Kg)(m_{1})-(\widehat Kg)(m_{2})|\leq\\
\\
\displaystyle
\sum_{\gamma}|\rho_{\gamma}(m_{1})-\rho_{\gamma}(m_{2})|\cdot
|F_{\gamma}(m_{1})|+\sum_{\gamma}\rho_{\gamma}(m_{2})
|F_{\gamma}(m_{1})-F_{\gamma}(m_{2})|\ .
\end{array}
$$
By (\ref{e726a}) the second sum is at most
$(\sup_{\gamma}||E_{\gamma}||)||g||_{Lip(F\cup\widehat F)}d(m_{1},m_{2})$
and together with (\ref{e726}) this leads to an appropriate bound. In turn,
the first sum is at most
$$
2\mu\cdot K\cdot\max_{\gamma}|(E_{\gamma}g_{\gamma})(m_{1})|
\cdot d(m_{1},m_{2})
$$
where $\mu$ (multiplicity) and $K$ are defined in Lemmae \ref{l71} and
\ref{l72}. To estimate the maximum,  one notes that
$(E_{\gamma}g_{\gamma})(\gamma)=g(\gamma)$ and therefore
$$
\begin{array}{c}
|(E_{\gamma}g_{\gamma})(m_{1})|\leq |(E_{\gamma}g_{\gamma})(m_{1})-
(E_{\gamma}g_{\gamma})(\gamma)|+|g(\gamma)|\leq\\
\\
\lambda_{R}d(m_{1},\gamma)||g||_{Lip(F\cup\widehat F)}+
||g||_{l_{\infty}(F\cup\widehat F)}\leq (R\lambda_{R}+1)
||g||_{(Lip\cap l_{\infty})(F\cup\widehat F)}\ .
\end{array}
$$
Together with the estimate of the second sum this proves the lemma.\ \ \ \ \
$\Box$

We are now ready to define the required operator $\widehat E$ from
$Ext(F,M)$. Actually, we use the above introduced operators $\widetilde E$,
$\widehat K$ and $\widehat S$ and set
$$
\widehat Ef:=\widetilde Ef+\widehat K(\widehat Sf|_{F\cup\widehat F})\ .
$$
If $m\in F$, we then have
$$
(\widehat Ef)(m)=(\widetilde Ef)(m)+(\widehat Sf)(m)=
(\widetilde Ef)(m)+(\widehat Tf)(m)-(\widetilde Ef)(m)=f(m)\ ,
$$
i.e., $\widehat E$ is an extension from $F$. To obtain the necessary
estimate of $||\widehat Ef||_{Lip(M)}$ it suffices by (\ref{e720})
to estimate
$||\widehat K(\widehat Sf|_{F\cup\widehat F})||_{Lip(M)}$.
The latter by Lemmae \ref{l76}, \ref{l75}, \ref{l74} and (\ref{e720})
is bounded by
$$
C||\widehat Sf||_{(Lip\cap l_{\infty})(F\cup\widehat F)}\leq
C||f||_{Lip(F)}\ .
$$
Hence $\widehat E\in Ext(F,M)$ and its norm is bounded as required.

The proof of the theorem is complete.\ \ \ \ \ $\Box$\\
{\bf Proof of Corollary \ref{2.12''}.}
Our initial proof derives this corollary from Theorem \ref{te17} and an
important result by Bonk and Schramm [BoSch]; it will be outlined below.
However, a recently established embedding theorem, see [NPSS],
allows to prove the desired result as an immediate consequence of Theorem
\ref{te14a}. The embedding theorem is formulated as follows.
 
{\em Let $M$ be a $\delta$-hyperbolic space of bounded geometry.
Then there exist constants $N\in\N$ and $C>0$ (depending on $M$) such that
$M$ is $C$-isometric to a subset of the direct sum of $N$ metric trees.}

This and Theorem \ref{te14a} immediately imply Corollary \ref{2.12''}.
\ \ \ \ \ $\Box$
\begin{R}\label{anoproof}
{\rm We outline another proof of Corollary \ref{2.12''} based (as the 
above formulated embedding theorem) on the main result of [BoSch]. This
result asserts:}

Let $M$ be a $\delta$-hyperbolic space of bounded geometry. Then there
exists an integer $n$ such that $M$ is roughly similar to a convex
subset of hyperbolic $n$-space $\H^{n}$.

{\rm Recall that a map $f:(M,d)\to (M',d')$ is a {\em $(K,L)$-rough
similarity}, if for all $m$, $n$ from $M$ it is true that}
$$
K^{-1}d(m,n)-L\leq d'(f(m),f(n))\leq Kd(m,n)+L\ .
$$

{\rm As a consequence of this result we obtain that the finite direct $p$-sum 
$M=\oplus_{p}\{M_{i}\}_{1\leq i\leq k}$ of Gromov-hyperbolic metric spaces 
$M_{i}$ is 
roughly similar to a subset of the direct $p$-sum 
$\oplus_{p}\{\H^{n_{i}}\}_{1\leq i\leq k}$
for some natural $n_{i}$. Now, the required corollary can be easily derived 
from the above result and Theorem \ref{te17} if we observe the following.\\
(1) Restriction $f|_{\Gamma}$ of a rough similar map $f$ to an $R$-lattice 
$\Gamma$ (see Definition 2.9) with $R$ big enough is a
$C$-isometric embedding into $M'$ for an appropriate $C$. If, in addition, 
$\lambda(M')<\infty$, then $\lambda(\Gamma)$ is finite.\\
(2) If $M_{i}$ is a geodesic metric space of $R_{0}$-bounded geometry, 
then every its ball of radius $R$ can be covered by at most $k=k(R,R_{0})$
balls of radius $R_{0}$. Therefore the same is true for 
the finite direct $p$-sum $M=\oplus_{p}\{M_{i}\}_{1\leq i\leq k}$ of such 
spaces.
From here one easily deduce that for the space $M$ the
constants $\lambda_{R}$ defined by (\ref{e16}) are finite for any $R$.\\
(3) According to [BSh2, Proposition 5.33] every $\H^{n_{i}}$ satisfies the 
assumptions of Corollary \ref{c118}. Therefore 
$\lambda(\oplus_{p}\{\H^{n_{i}}\}_{1\leq i\leq k})<\infty$.}
\end{R}
{\bf Proof of Corollary \ref{c111b}.}
We first prove that the condition
\begin{equation}\label{e729}
\lambda(G,d_{A})<\infty
\end{equation}
of Corollary \ref{c111b} is necessary for finiteness of $\lambda(M)$.
Let the latter be true. Then for a $G$-orbit
$G(m):=\{g(m)\ :\ g\in G\}$ we have
\begin{equation}\label{e730}
\lambda(G(m))\leq\lambda(M)<\infty\ .
\end{equation}
Then the \v{S}varc-Milnor lemma, see, e.g., [BH, p.140], states that under
the hypothesis (b) of Corollary \ref{c111b} there is a
constant $C\geq 1$ independent of $m$ so that
\begin{equation}\label{e731}
C^{-1}d_{A}(g,h)\leq d(g(m),h(m))\leq Cd_{A}(g,h)
\end{equation}
for all $g,h\in G$. This, in particular, means that the metric subspace
$G(m)$ is quasi-isometric to the metric space $(G,d_{A})$. Hence
(\ref{e730}) implies the required inequality (\ref{e729}).

To prove sufficiency of the condition (\ref{e729}) for finiteness of
$\lambda(M)$, we choose a point $m_{0}$ of the generating compact set $K_{0}$
from Definition \ref{d111a}, see
(\ref{e17a}), and show that the $G$-orbit $G(m_{0})$ is an $R$-lattice for
some $R>0$. Let $B_{R_{0}}(m_{0})$ be a ball
containing $K_{0}$. Then we have by (\ref{e17a}) for $\Gamma:=G(m_{0})$
$$
\bigcup_{m\in\Gamma}B_{R_{0}}(m)=G(B_{R_{0}}(m_{0}))\supset G(K_{0})=M\ .
$$
Hence the family of balls $\{B_{R_{0}}(m)\ :\ m\in\Gamma\}$
covers $M$. Moreover, (\ref{e731}) implies that for $m:=g(m_{0})$,
$m':=h(m_{0})$ with $g\neq h$
$$
d(m,m')\geq C^{-1}d_{A}(g,h)\geq C^{-1}\ ,
$$
that is to say, the family $\{B_{cR_{0}}(m)\ :\ m\in\Gamma\}$
with $c=c_{\Gamma}:=\frac{1}{2CR_{0}}$ consists of pairwise disjoint balls.
So $\Gamma$ is an $R$-lattice, $R:=2R_{0}$, satisfying, by (\ref{e729}) and
(\ref{e731}), the condition
$$
\lambda(\Gamma)<\infty\ .
$$

We now will apply Theorem \ref{te17} with that $R$-lattice $\Gamma$ to
derive finiteness of $\lambda(M)$.
To this end we have to establish validity of the assumptions of the theorem
with this $R$.

First we prove that $M$ belongs to the class of doubling metric spaces
${\cal D}(R,N)$ for some $N=N(R,M)$. In other words, we show that every ball
$B_{r}(m)$ with $r<R$ can be covered by at most $N$ balls of radius
$r/2$. Indeed, by the hypothesis (a) of the corollary,
$M\in {\cal G}_{n}(\widetilde R,\widetilde C)$ for certain
$\widetilde R,\widetilde C$ and $n$. This implies that
$M\in {\cal D}(\widetilde R/2,N)$ for some $N=N(\widetilde C,n)$ and shows
that the required statement is true for $R\leq\widetilde R/2$. Suppose now
that
\begin{equation}\label{e732a}
\widetilde R/2\leq r<R\ .
\end{equation}
Note that it suffices to consider balls with $m\in K_{0}$. In fact,
$G(K_{0})=M$ and therefore $g_{0}(m)\in K_{0}$ for some isometry
$g_{0}\in G$. Hence $g_{0}(B_{r}(m))=B_{r}(g_{0}(m))$ and we can work with
$B_{r}(m)$ for $m\in K_{0}$. Let us fix a point $m_{0}\in K_{0}$ and set
$R':=R+diam\ \! K_{0}$. Then
\begin{equation}\label{e732}
B_{r}(m)\subset B_{R'}(m_{0})\ ,\ \ \ m\in K_{0}\ ,
\end{equation}
and it remains to show that $B_{R'}(m_{0})$ can be covered by a finite
number, say $N$, of (open) balls of radius $r/2$ with $N$ independent of
$r$. We use the following
\begin{Lm}\label{l733}
Suppose that $G$ acts properly, freely and cocompactly on a path space $M$
by isometries. Then every bounded closed set $S\subset M$ is compact.
\end{Lm}
{\bf Proof.} For every $m\in S$ there is a finite number of
isometries $g_{im}\in G$, $i=1,\dots, k_{m}$, such that
$g_{im}(m)\in K_{0}$. Here $K_{0}$ is the generating compact of
Definition \ref{d111a} (c). Let $H:=\{g_{im}^{-1}\in G\ :\ 1\leq
i\leq k_{m},\ m\in S\}$. Then $S\subset H(K_{0})$, and, by that
definition, $diam\ \! H(K_{0})<\infty$. For a fixed $m_{0}\in
K_{0}$ let us consider the orbit $H(m_{0})$. Show that $H(m_{0})$
consists of a finite number of points. Otherwise there is a
sequence of points $m_{i}=h_{i}(m_{0})\in H(m_{0})$ such that
$d_{A}(h_{i},1)\to\infty$ as $i\to\infty$. This and inequality
(\ref{e731}) imply $d(m_{i},m_{0})\to\infty$ in $M$ as
$i\to\infty$ and this contradicts to the condition $diam\ \!
H(K_{0})<\infty$. From finiteness of $H(m_{0})$ we also obtain
that $H$ is finite. Thus $S$ is covered by a finite number of
compact sets, and, since $S$ is closed, it is compact. \ \ \ \ \
$\Box$

According to this lemma the closure $\overline{B_{R'}(m_{0})}$ is
compact. Thus $B_{R'}(m_{0})$ can be covered by a finite number
$N$ of open balls of radius $\widetilde R/4$. This, (\ref{e732a})
and (\ref{e732}) show that $B_{r}(m)$ can be covered by $N$ open
balls of radius $r/2$ as it is required.

To establish the second condition of Theorem \ref{te17}, finiteness of
$$
\lambda_{R}:=\sup\{\lambda(B_{R}(m))\ :\ m\in M\}\ ,
$$
we first use the previous argument and (\ref{e732}) which immediately get
$$
\lambda_{R}\leq\lambda(B_{R'}(m_{0}))\ .
$$
Show that the right-hand side is bounded. Since $M\in {\cal
G}_{n}(\widetilde R,\widetilde C)$, for every $m\in M$,
$$
\lambda(B_{\widetilde R}(m))<\infty\ .
$$
This, compactness of $\overline{B_{R'}(m_{0})}$ and the argument used in the
proof of Lemma \ref{l45} lead to the required inequality
$$
\lambda(B_{R'}(m_{0}))<\infty\ .
$$

The proof of Corollary \ref{c111b} is completed.\ \ \ \ \ $\Box$

\sect{Proofs of Theorem \ref{te114} and its Corollaries} {\bf
Proof of Theorem \ref{te114}; Part I.} Given a metric space 
$(M,d)$ of pointwise homogeneous type 
of Definition \ref{d218} and a
subspace $S\subset M$ we will construct an operator $E\in
Ext(S,M)$ whose norm is bounded by a constant depending only on
the constants $C$ of (\ref{eq28}), and
\begin{equation}\label{e141}
D(l):=\sup_{m\in
M}\sup_{R>0}\frac{\mu_{m}(B_{lR}(m))}{\mu_{m}(B_{R}(m))}
\end{equation}
with $l>1$ that will be specified later. By the uniform doubling
condition (\ref{eq27}), this is finite and depends only on the
constant in (\ref{eq27}). Our construction is similar to that of
[BSh2, pp.535-540]; unfortunately, the latter used a Borel
measurable selection of the multivalued function $m\mapsto\{m'\in
M\ : d(m,S)\leq d(m,m')\leq 2d(m,S)\}$ for some specific spaces
$(M,d)$ (including e.g. the $n$-dimensional hyperbolic space).
Generally speaking, such a selection may not exist even in the
case of the Euclidean plane\footnote{see, in particular, the
corresponding counterexample in [N].}. Fortunately, Theorem
\ref{te11} allows to restrict our consideration to the case of
finite point subspaces $S$ in which case the corresponding
measurable selection trivially exists.

So we consider a {\em finite point} metric subspace $S\subset M$ and construct
in this case an operator $E\in Ext(S,M)$ with
\begin{equation}\label{e142}
||E||\leq K=K(D(l),C)<\infty\ .
\end{equation}
By Theorem \ref{te11} $\lambda(M)$ will be then bounded by the same
constant $K$. 

For this purpose let us arrange $S$ in a sequence
$s_{1},\dots, s_{l}$ and introduce functions $d:M\rightarrow\Re_{+}$ and
$p:M\rightarrow S$ by the conditions
\begin{equation}\label{e143}
d(m):=\min\{d(m,m')\ :\ m'\in S\}
\end{equation}
and
\begin{equation}\label{e144}
p(m):=s_{i}
\end{equation}
where $i$ is the minimal number for which
$s_{i}\in \{m'\in S\ :\ d(m)=d(m',m)\}$.
\begin{Lm}\label{l141}
(a) For every $m_{1},m_{2}\in M$
\begin{equation}\label{e145}
|d(m_{1})-d(m_{2})|\leq d(m_{1},m_{2})\ .
\end{equation}
(b) If $f:S\rightarrow\Re$ is an arbitrary function, then $f\circ p$ is
Borel measurable.
\end{Lm}
{\bf Proof.} (a) follows directly from (\ref{e143}). To check (b) note that
for each $1\leq i\leq l$ the set $p^{-1}(\{s_{1},\dots,s_{i}\})\subset M$ 
is closed. This implies the required result.\ \ \ \ \ $\Box$

To introduce the required extension operator $E$ we also use the
average with respect to the Borel measure $\mu_{m}$ of Definition
\ref{d218}, letting for a Borel measurable function
$g:M\rightarrow\Re$
\begin{equation}\label{e146}
I(g;m,R):=\frac{1}{\mu_{m}(B_{R}(m))}\cdot\int_{B_{R}(m)}g\
\!d\mu_{m}\ ,\ \ \ m\in M\ , R>0\ .
\end{equation}
Finally, we define $E$ on functions $f\in Lip(S)$ by
\begin{equation}\label{e147}
(Ef)(m):=\left\{
\begin{array}{ccc}
f(m),&{\rm if}&m\in S\\
I(f\circ p\ \!;m,d(m)),&{\rm if}&m\in M\setminus S
\end{array}
\right.
\end{equation}
We now have to show that for every $m_{1},m_{2}$
\begin{equation}\label{e148}
|(Ef)(m_{1})-(Ef)(m_{2})|\leq K||f||_{Lip(S)}d(m_{1},m_{2})
\end{equation}
where $K=K(D(l),C)$ will be specified later.\\
It suffices to consider only two cases:
\begin{itemize}
\item[(a)]
$m_{1}\in S$\ and\ $m_{2}\not\in S$;
\item[(b)]
$m_{1}, m_{2}\not\in S$.
\end{itemize}
We assume without loss of generality that
\begin{equation}\label{e149}
||f||_{Lip(S)}=1
\end{equation}
and simplify computations by introducing the following notations:
\begin{equation}\label{e1410}
R_{i}:=d(m_{i})\ ,\ \mu_{i}:=\mu_{m_{i}}\ ,\ B_{ij}:=B_{R_{j}}(m_{i})\ ,
\ v_{ij}:=\mu_{i}(B_{ij})\ ,\ \ 1\leq i,j\leq 2\ .
\end{equation}
We assume also for definiteness that
\begin{equation}\label{e1411}
0< R_{1}\leq R_{2}\ .
\end{equation}
By Lemma \ref{l141} we then have
\begin{equation}\label{e1412}
0\leq R_{2}-R_{1}\leq d(m_{1},m_{2})\ .
\end{equation}
In what follows we will prove the required result under the next additional
assumption on the family $\{\mu_{m}\}$:

{\em There is a constant $A>0$ such that for all $0< R_{1}\leq
R_{2}$ and $m\in M$}
\begin{equation}\label{e1413}
\mu_{m}(B_{R_{2}}(m))-\mu_{m}(B_{R_{1}}(m))\leq\frac{A\mu_{m}(B_{R_{2}}(m))}{R_{2}}
(R_{2}-R_{1})\ .
\end{equation}

This restriction will be removed at the final stage of the proof.

Under the notations and the assumptions introduced the following
is true
\begin{equation}\label{e1414}
v_{i2}-v_{i1}\leq\frac{Av_{i2}}{R_{2}}(R_{2}-R_{1})\ ,
\end{equation}
\begin{equation}\label{e1416}
|\mu_{1}-\mu_{2}|(B_{ij})\leq\frac{Cv_{ij}}{R_{j}}d(m_{1},m_{2})\
,
\end{equation}
see (\ref{e1413}) and (\ref{eq28}).

To estimate the difference in (\ref{e148}) for $m_{1}\in S$ we need
\begin{Lm}\label{l142}
It is true that
\begin{equation}\label{e1417}
\max\{|\widetilde f(m)|\ :\ m\in B_{i2}\}\leq 4R_{2}+(i-1)d(m_{1},m_{2})\ ;
\end{equation}
here $i=1,2$ and
\begin{equation}\label{e1418}
\widetilde f(m):=(f\circ p)(m)-(f\circ p)(m_{1})\ .
\end{equation}
\end{Lm}
{\bf Proof.}
Let first $i=1$ and $m\in B_{12}$. By (\ref{e1418}), (\ref{e149}), and
the triangle inequality
$$
|\widetilde f(m)|\leq d(p(m),p(m_{1}))\leq d(m)+d(m,m_{1})+d(m_{1})\ .
$$
But $d(m,m_{1})\leq R_{2}$, since $m\in B_{12}$. Besides,
$d(m)\leq d(m,m_{1})+d(m_{1})\leq 2R_{2}$ by Lemma \ref{l141}. Taking
these together to get
$$
|\widetilde f(m)|\leq 4R_{2}\ ,\ m\in B_{12}\ .
$$
Let now $i=2$ and $m\in B_{22}$. As before, the triangle inequality gives
$$
|\widetilde f(m)|\leq d(m,m_{1})+d(m,m_{2})+d(m_{1})+d(m_{2})\ .
$$
Since $d(m,m_{2})\leq R_{2}$ and $d(m,m_{1})\leq d(m,m_{2})+d(m_{1},m_{2})$,
we therefore have
$$
|\widetilde f(m)|\leq d(m_{1},m_{2})+4R_{2}\ ,\ \ \ m\in B_{22}\ .
\ \ \ \ \ \Box
$$

Prove now (\ref{e148}) for $m_{1}\in S$ and $m_{2}\not\in S$. We clearly
have under the notation (\ref{e1410})
$$
|(Ef)(m_{2})-(Ef)(m_{1})|=
\frac{1}{v_{22}}\left|\int_{B_{22}}\widetilde f(m)d\mu_{2}\right|\leq
\max_{B_{22}}|\widetilde f|\ .
$$
Applying (\ref{e1417}) with $i=2$ we then bound this difference by
$4R_{2}+d(m_{1},m_{2})$. But $m_{1}\in S$ and so
$$
R_{2}=d(m_{2})\leq d(m_{1},m_{2})\ ,
$$
and therefore (\ref{e148}) holds in this case with $K=5$.

The remaining case $m_{1},m_{2}\not\in S$ requires some additional
auxiliary results. To their formulations we first write
\begin{equation}\label{e1419}
(Ef)(m_{1})-(Ef)(m_{2}):=D_{1}+D_{2}
\end{equation}
where
\begin{equation}\label{e1420}
\begin{array}{c}
D_{1}:=I(\widetilde f; m_{1},R_{1})-I(\widetilde f; m_{1},R_{2})\\
\\
D_{2}:=I(\widetilde f; m_{1},R_{2})-I(\widetilde f; m_{2},R_{2})\ ,
\end{array}
\end{equation}
see (\ref{e147}) and (\ref{e1418}).
\begin{Lm}\label{l143}
It is true that
$$
|D_{1}|\leq 8Ad(m_{1},m_{2})\ .
$$
Recall that $A$ is the constant in (\ref{e1413}).
\end{Lm}
{\bf Proof.}
By (\ref{e1420}), (\ref{e1418}) and (\ref{e1410}),
$$
D_{1}=\frac{1}{v_{11}}\int_{B_{11}}\widetilde f d\mu_{1}-
\frac{1}{v_{12}}\int_{B_{12}}\widetilde f d\mu_{1}=
\left(\frac{1}{v_{11}}-\frac{1}{v_{12}}\right)\int_{B_{11}}
\widetilde f d\mu_{1}-\frac{1}{v_{12}}\int_{B_{12}\setminus B_{11}}
\widetilde f d\mu_{1}\ .
$$
This immediately implies that
$$
|D_{1}|\leq 2\cdot \frac{v_{12}-v_{11}}{v_{12}}\cdot
\max_{B_{12}}|\widetilde f|\ .
$$
Applying now (\ref{e1414}) and (\ref{e1412}), and then Lemma
\ref{l142} with $i=1$ we get the desired estimate.\ \ \ \ \ $\Box$

To obtain a similar estimate for $D_{2}$ we will use the following
two facts.
\begin{Lm}\label{l144}
Assume that for a given $l>1$
\begin{equation}\label{e1421}
d(m_{1},m_{2})\leq (l-1)R_{2}\ .
\end{equation}
Let for definiteness
\begin{equation}\label{eq1421}
v_{22}\leq v_{12}\ .
\end{equation}
Then it is true that
\begin{equation}\label{e1422}
\mu_{2}(B_{12}\Delta B_{22})\leq 2(A+C)D(l)\frac{v_{12}}{R_{2}}
d(m_{1},m_{2})
\end{equation}
(here $\Delta$ denotes symmetric difference of sets).
\end{Lm}
{\bf Proof.}
Set
$$
R:=R_{2}+d(m_{1},m_{2})\ .
$$
Then $B_{12}\cup B_{22}\subset B_{R}(m_{1})\cap B_{R}(m_{2})$ and
\begin{equation}\label{eq1223}
\mu_{2}(B_{12}\Delta B_{22})\leq (\mu_{2}(B_{R}(m_{1}))-\mu_{2}(B_{12}))
+(\mu_{2}(B_{R}(m_{2}))-\mu_{2}(B_{22}))\ .
\end{equation}
The first summand on the right-hand side is at most
$$
|\mu_{2}-\mu_{1}|(B_{R}(m_{1}))+|\mu_{2}-\mu_{1}|(B_{R_{2}}(m_{1}))+
(\mu_{1}(B_{R}(m_{1}))-\mu_{1}(B_{R_{2}}(m_{1}))\ .
$$
Estimating the first two summands by (\ref{eq28}) and the third by
(\ref{e1413}) we bound this sum by
$$
C\left(\frac{\mu_{1}(B_{R}(m_{1}))}{R}+\frac{\mu_{1}(B_{R_{2}}(m_{1}))}{R_{2}}\right)+
A\frac{\mu_{1}(B_{R}(m_{1}))}{R}\ \! (R-R_{2})\ .
$$
Besides, $R_{2}\leq R\leq lR_{2}$ and $R-R_{2}:=d(m_{1},m_{2})$,
see (\ref{e1421}); taking into account (\ref{e141}) and the
notations (\ref{e1410}) we therefore have
$$
\mu_{2}(B_{R}(m_{1}))-\mu_{2}(B_{12})\leq
[C(D(l)+1)+AD(l)]\frac{v_{12}}{R_{2}}\ \! d(m_{1},m_{2})\ .
$$
Similarly, by (\ref{e1413}) and (\ref{eq1421})
$$
\begin{array}{c}
\displaystyle \mu_{2}(B_{R}(m_{2}))-\mu_{2}(B_{22})\leq
A\frac{\mu_{2}(B_{R}(m_{2}))}{R}\ \! (R-R_{2})\leq\\
\\
\displaystyle AD(l)\frac{v_{22}}{R_{2}}d(m_{1},m_{2})\leq
AD(l)\frac{v_{12}}{R_{2}}d(m_{1},m_{2})\ .
\end{array}
$$
Combining the last two estimates with (\ref{eq1223}) we get the
result. \ \ \ \ \ $\Box$
\begin{Lm}\label{l145}
Under the assumptions of the previous lemma it is true that
\begin{equation}\label{e1423}
v_{12}-v_{22}\leq 3(A+C)D(l)\frac{v_{12}}{R_{2}}d(m_{1},m_{2})\ .
\end{equation}
\end{Lm}
{\bf Proof.} By (\ref{e1410}) the left-hand side is bounded by
$$
|\mu_{1}(B_{12})-\mu_{2}(B_{12})|+\mu_{2}(B_{12}\Delta B_{22})\ .
$$
Estimating these summands by (\ref{e1416}) and (\ref{e1422}) we
get the result. \ \ \ \ \ $\Box$

We now estimate $D_{2}$ from (\ref{e1420}) beginning with
\begin{Lm}\label{l146}
Under the conditions of Lemma \ref{l144} it is true that
$$
|D_{2}|\leq K(l)d(m_{1},m_{2})
$$
where
\begin{equation}\label{eq1225}
K(l):=6(A+C)D(l)(l+3)\ .
\end{equation}
\end{Lm}
{\bf Proof.} By the definition of $D_{2}$ and our notations, see
(\ref{e1420}), (\ref{e1417}) and (\ref{e1410}),
$$
\begin{array}{c}
\displaystyle
|D_{2}|:=\left|\frac{1}{v_{12}}\int_{B_{12}}\widetilde f d\mu_{1}-
\frac{1}{v_{22}}\int_{B_{22}}\widetilde f d\mu_{2}\right|\leq
\frac{1}{v_{12}}\int_{B_{12}}|\widetilde f|\ \! d|\mu_{1}-\mu_{2}|\ +
\\
\\
\displaystyle
\frac{1}{v_{12}}\int_{B_{12}\Delta B_{22}}|\widetilde f|\ \! d\mu_{2}
+\left|\frac{1}{v_{12}}-\frac{1}{v_{22}}\right|\int_{B_{22}}|\widetilde f|
\ \! d\mu_{2}:=J_{1}+J_{2}+J_{3}\ .
\end{array}
$$
By (\ref{e1416}), (\ref{e1417}) with $i=1$ and (\ref{e1421})
$$
J_{1}\leq\frac{1}{v_{12}}|\mu_{1}-\mu_{2}|(B_{12}) \sup_{B_{12}}
|\widetilde f|\leq
\frac{C}{R_{2}}d(m_{1},m_{2})(d(m_{1},m_{2})+4R_{2})\leq
C(l+3)d(m_{1},m_{2})\ .
$$
In turn, by (\ref{e1422}), (\ref{e1411}) and (\ref{e1417})
$$
\begin{array}{c}
\displaystyle J_{2}\leq\frac{1}{v_{12}}\mu_{2}(B_{12}\Delta
B_{22})\sup_{B_{12}\Delta B_{22}}|\widetilde f|\leq
\frac{2(A+C)D(l)}{R_{2}}
d(m_{1},m_{2})(d(m_{1},m_{2})+4R_{2})\leq\\
\\
\displaystyle 2(A+C)D(l)(l+3)d(m_{1},m_{2})\ .
\end{array}
$$
Finally, (\ref{e1423}), (\ref{e1417}) and (\ref{e1421}) yield
$$
J_{3}\leq 3(A+C)D(l)(l+3)d(m_{1},m_{2})\ .
$$
Combining we get the required estimate.\ \ \ \ \ $\Box$

It remains to consider the case of $m_{1},m_{2}\in M$ satisfying
the inequality
$$
d(m_{1},m_{2})>(l-1)R_{2}
$$
converse to (\ref{e1421}). Now the definition (\ref{e1420}) of
$D_{2}$ and (\ref{e1417}) imply that
$$
|D_{2}|\leq 2\sup_{B_{12}\cup B_{22}}|\widetilde f|\leq
2(4R_{2}+d(m_{1},m_{2}))\leq 2\left(\frac{4}{l-1}+1\right)d(m_{1},m_{2})\ .
$$
Combining this with the inequalities of Lemmae \ref{l143} and
\ref{l146} and equality (\ref{e1419}) we obtain the required
estimate of the Lipschitz norm of the extension operator $E$
defined by (\ref{e147}). Actually, we have proved that
\begin{equation}\label{e1424}
||E||\leq  8A+ \max\left(\frac{2(l+3)}{l-1},K(l)\right)
\end{equation}
where $K(l)$ is the constant in (\ref{eq1225}). Hence Theorem
\ref{te114} has proved under the additional assumption
(\ref{e1413}) with $\lambda(M)$ estimated by (\ref{e1424}).
\begin{R}\label{r147}
{\rm (a) Let $M$ be a metric space of homogeneous type with respect to a
doubling measure $\mu$. Taking $\mu_{m}:=\mu$ for all $m\in M$ and
noting that (\ref{e1416}) is now trivially held we improve the
estimate (\ref{e1424}) as follows:
\begin{equation}\label{eq1227}
\lambda(M)\leq 8A+\max\left(\frac{2(l+3)}{l-1},K_{\mu}(l)\right)
\end{equation}
where $K_{\mu}(l):=2A(l+3)D_{\mu}(l)$. Here $D_{\mu}(l)$ is the
dilation function for $\mu$, see (\ref{e141}). In fact,
(\ref{eq1223}) is now bounded by
$\frac{2AD(l)}{R_{2}}d(m_{1},m_{2})$ and this constant appears
in (\ref{e1424}). Besides, $J_{1}=J_{3}=0$ in this case.\\
(b) If, on the other hand, for some $a,n>0$ and all $B_{R}(m)$
\begin{equation}\label{eq1228}
\mu_{m}(B_{R}(m))=aR^{n}\ ,
\end{equation}
the estimate (\ref{eq1227}) can be sharpen. Note that in this case
condition (\ref{e1413}) clearly holds with $A=n$. Hence Theorem
\ref{te114} had already proved in this case. Besides, in the proof
of Lemma \ref{l146}, $J_{1}$ is now bounded by
$C(l+3)d(m_{1},m_{2})$, and $J_{2}$ and $J_{3}$ by
$n(l+3)l^{n-1}d(m_{1},m_{2})$ and $nl^{n-1}(l+3)d(m_{1},m_{2})$,
respectively. Collecting these we get in this case
\begin{equation}\label{eq1229}
\lambda(M)\leq 8n+\max\left(\frac{2(l+3)}{l-1}, K_{n}(l)\right)
\end{equation}
where
$$
K_{n}(l):=(l+3)(C+2nl^{n-1})\ .
$$
(c) Finally, for the case of the doubling measure $\mu$ of part
(a) satisfying condition (\ref{eq1228}) the constant $C$ in
(\ref{eq1229}) disappears and we get the estimate (\ref{eq1229})
with $K_{n}(l)=2n(l+3)l^{n-1}$.

Let us recall that $l>1$ is arbitrary and we may and will optimize
all these estimates with respect to $l$.}
\end{R}
{\bf Proof of Theorem \ref{te114}; Part II.} We apply now the result of the 
previous part to a metric space  
$(\widehat M,\widehat d)$ of pointwise homogeneous type 
satisfying the following two conditions
\begin{itemize}
\item[(a)] {\em The original metric space $(M,d)$ embeds isometrically
to $(\widehat M,\widehat d)$.}
\end{itemize}

This immediately implies the inequality
\begin{equation}\label{eq1230}
\lambda(M)\leq\lambda(\widehat M)\ .
\end{equation}
\begin{itemize}
\item[(b)] {\em Condition (\ref{e1413}) holds for $(\widehat M,\widehat d)$.}
\end{itemize}

This implies validity of estimate (\ref{e1424}) for the extension operator
$E\in Ext(S,\widehat M)$. Of course, this estimate includes now the
dilation function and the consistency constant for $\widehat M$
which must be evaluated via the corresponding amounts for $M$, see
(\ref{e141}) and (\ref{eq28}). This goal will be achieved by two
auxiliary results, Lemmae \ref{l128} and \ref{l129} presented below.

To their formulation, introduce a metric space $(M_{N},d_{N})$ by
$$
M_{N}:=M\times l_{1}^{N}\ ;
$$
where $l_{1}^{N}$ is the $N$-dimensional vector space defined by the metric
$$
\delta_{1}^{N}(x,y):=||x-y||_{1}=\sum_{i=1}^{N}|x_{i}-y_{i}|\ ,\ \
x,y\in\Re^{N}\ .
$$
In turn, $d_{N}$ is given by
$$
d_{N}(\widetilde m,\widetilde m'):=d(m,m')+\delta_{1}^{N}(x,x')
$$
where here and below
$$
\widetilde m:=(m,x)\ \ \ {\rm with}\ \ \ m\in M\ \ \ {\rm and}\ \ \ 
x\in l_{1}^{N}\ .
$$
At last, we equip $M_{N}$ with a family
$\widetilde{\cal F}:=\{\mu_{\widetilde{m}}\}_{\widetilde m\in M_{N}}$
of positive Borel measures on $M_{N}$ introduced by 
$$
\mu_{\widetilde{m}}:=\mu_{m}\otimes\lambda_{N}
$$
where $\lambda_{N}$ is the Lebesgue measure of $\Re^{N}$ and 
${\cal F}:=\{\mu_{m}\}_{m\in M}$ is the family of doubling measures from the 
definition of $M$. 

Now we estimate the dilation function of family $\widetilde{\cal F}$ using
the corresponding doubling inequality for family ${\cal F}$. Recall
that this inequality, see (\ref{eq27}), implies that
\begin{equation}\label{e1231}
\mu_{m}(B_{2R}(m))\leq D\mu_{m}(B_{R}(m))
\end{equation}
for all $m\in M$ and $R>0$.

Thus, we consider now the function
\begin{equation}\label{e1232}
D_{N}(l):=\sup\left\{\frac{\mu_{\widetilde m}(B_{lR}(\widetilde m))}
{\mu_{\widetilde m}(B_{R}(\widetilde m))}\right\}
\end{equation}
where the supremum is taken over all $\widetilde m\in M_{N}$ and $R>0$.
\begin{Lm}\label{l128}
Assume that $N$ is related to the doubling constant $D$ of (\ref{e1231})
by
\begin{equation}\label{e1233}
N\geq [\log_{2} D]+5\ .
\end{equation}
Then it is true that
$$D_{N}(1+1/N)\leq \frac{6}{5}e^{4}\ .
$$
\end{Lm}
{\bf Proof.}
Note that open ball $B_{R}(\widetilde m)$ of $M_{N}$ is the set
$$
\{(m',y)\in M\times l_{1}^{N}\ :\ d(m',m)+||x-y||_{1}<R\}
$$
(recall that $\widetilde m=(m,x)$). Therefore application of Fubini's theorem
yields
\begin{equation}\label{e1234}
\mu_{\widetilde m}(B_{R}(\widetilde m))=\gamma_{N}\int_{B_{R}(m)}
(R-d(m,m'))^{N}d\mu_{m}(m')\ ;
\end{equation}
here $B_{R}(m)$ is a ball of $M$ and $\gamma_{N}$ is the volume of the 
unit $l_{1}^{N}$-ball. 

Estimate this measure with $R$ replaced by
$$
R_{N}:=(1+1/N)R\ .
$$
Split the integral in (\ref{e1234}) into those over $B_{3R/4}(m)$ and
over the remaining part $B_{R_{N}}(m)\setminus B_{3R/4}(m)$. Denote these
integrals by $I_{1}$ and $I_{2}$. For $I_{2}$ we get from (\ref{e1234})
$$
I_{2}\leq \gamma_{N}(R_{N}-3R/4)^{N}\int_{B_{R_{N}}(m)}d\mu_{m}(m')=
\gamma_{N}
\left(\frac{1}{4}+\frac{1}{N}\right)^{N}R^{N}\mu_{m}(B_{R_{N}}(m))\ .
$$
Using the doubling constant for ${\cal F}=\{\mu_{m}\}$, see (\ref{e1231}),
we further have
$$
\mu_{m}(B_{R_{N}}(m))\leq D\mu_{m}(B_{R_{N}/2}(m))\ .
$$
Moreover, by (\ref{e1233}), $D<2^{[\log_{2}D]+1}\leq\frac{1}{16}2^{N}$.
Combining all these inequalities we obtain
\begin{equation}\label{e1235}
I_{2}\leq\gamma_{N}\frac{1}{16}2^{-N}\left(1+\frac{4}{N}\right)^{N}R^{N}
\mu_{m}(B_{R_{N}/2}(m))\ .
\end{equation}

To estimate $I_{1}$ we present its integrand (which equals to that in
(\ref{e1234}) with $R$ replaced by $R_{N}$) in the following way.
$$
\left(1+\frac{1}{N}\right)^{N}(R-d(m,m'))^{N}\left(1+
\frac{d(m,m')}{(N+1)(R-d(m,m'))}\right)^{N}\ .
$$
Since $d(m,m')\leq 3R/4$ for $m'\in B_{3R/4}(m)$, the last factor is 
at most $\left(1+\frac{3}{N+1}\right)^{N}$. Hence, we have
$$
I_{1}\leq\gamma_{N}
\left(1+\frac{1}{N}\right)^{N}\left(1+\frac{3}{N+1}\right)^{N}
\int_{B_{3R/4}(m)}(R-d(m,m'))^{N}d\mu_{m}(m')\ .
$$
Using then (\ref{e1234}) we, finally, obtain
\begin{equation}\label{e1236}
I_{1}\leq e^{4}\mu_{\widetilde m}(B_{R}(\widetilde m))\ .
\end{equation}

To estimate $D_{N}(l)$ with $l=1+1/N$ it remains to bound fractions
$$
\widetilde I_{k}:=\frac{I_{k}}{\mu_{\widetilde m}(B_{R}(\widetilde m))}\ ,
\ \ \ k=1,2\ .
$$
For $k=2$ estimate the denominator from below as follows. Since 
$R_{N}<2R$, we bound $\mu_{\widetilde m}(B_{R}(\widetilde m))$
from below by
$$
\begin{array}{c}
\displaystyle
\gamma_{N}\int_{B_{R_{N}/2}(m)}(R-d(m,m'))^{N}d\mu_{m}(m')\geq
\gamma_{N} 2^{-N}
\left(1-\frac{1}{N}\right)^{N}R^{N}\int_{B_{R_{N}/2}(m)}d\mu_{m}(m')=\\
\\
\displaystyle
\gamma_{N}2^{-N}\left(1-\frac{1}{N}\right)^{N}R^{N}\mu_{m}(B_{R_{N}/2}(m))\ 
.
\end{array}
$$
Combining this with (\ref{e1235}) to get
$$
\widetilde I_{2}\leq\frac{1}{16}\left(1-\frac{1}{N}\right)^{-N}
\left(1+\frac{4}{N}\right)^{N}\ .
$$
Since $\left(1-\frac{1}{N}\right)^{-N}\leq\left(1-\frac{1}{5}\right)^{-5}$ 
as $N\geq 5$, we finally obtain
$$
\widetilde I_{2}\leq\frac{1}{5}e^{4}\ .
$$
For $\widetilde I_{1}$ using (\ref{e1236}) one immediately has
$$
\widetilde I_{1}\leq e^{4}\ .
$$
Hence, by the definition of $D_{N}$, see (\ref{e1232}), we have
$$
D_{N}(1+1/N)\leq\sup_{\widetilde m, R}(\widetilde I_{1}+\widetilde I_{2})<
\frac{6}{5}e^{4}\ .\ \ \ \ \ \Box
$$

Our next auxiliary result evaluate the consistency constant $C_{N}$ for family
$\widetilde{\cal F}=\{\mu_{\widetilde m}\}$ using that for
${\cal F}:=\{\mu_{m}\}$. Recall that the latter constant stands in the
inequality, see (\ref{eq28}),
\begin{equation}\label{e1237}
|\mu_{m_{1}}-\mu_{m_{2}}|(B_{R}(m_{i}))\leq
\frac{C\mu_{m_{i}}(B_{R}(m_{i}))}{R}d(m_{1},m_{2})
\end{equation}
where $m_{1},m_{2}$ are arbitrary points of $M$, $R>0$ and $i=1,2$.
\begin{Lm}\label{l129}
$$
C_{N}\leq\left(1+\frac{4e}{3}\right)NC\ .
$$
\end{Lm}
{\bf Proof.}
Using the Fubini theorem, rewrite (\ref{e1234}) in the form
\begin{equation}\label{e1238}
\mu_{\widetilde m}(B_{R}(\widetilde m))=\beta_{N}\int_{0}^{R}
\mu_{m}(B_{s}(m))(R-s)^{N-1}ds\ 
\end{equation}
where $\beta_{N}$ is the volume of the unit sphere in $l_{1}^{N}$.
Then for $i=1,2$ we have
$$
|\mu_{\widetilde m_{1}}-\mu_{\widetilde m_{2}}|(B_{R}(\widetilde m_{i}))
\leq\beta_{N}\int_{0}^{R}|\mu_{m_{1}}-\mu_{m_{2}}|(B_{s}(m_{i}))\cdot
(R-s)^{N-1}ds\ .
$$
Divide here the interval of integration into subintervals $[0,R/N]$ and
$[R/N,R]$ and denote the corresponding integrals over these intervals by
$I_{1}$ and $I_{2}$. It suffices to find appropriate upper bounds for
$I_{k}$. Replacing $B_{s}(m_{i})$ in $I_{1}$ by bigger ball $B_{s+R/N}(m_{i})$
and applying (\ref{e1237}) we obtain
$$
I_{1}\leq 
C\left(\beta_{N}\int_{0}^{R/N}\frac{\mu_{m_{i}}(B_{s+R/N}(m_{i}))}
{s+R/N}(R-s)^{N-1}ds\right)d(m_{1},m_{2})\ .
$$
Changing $s$ by $t=s+R/N$ we bound the expression in the brackets by
$$
\left(\beta_{N}\int_{R/N}^{2R/N}\mu_{m_{i}}(B_{t}(m_{i})) (R-t)^{N-1}dt
\right)
\max_{R/N\leq t\leq 2R/N}\frac{(R+R/N-t)^{N-1}}{t(R-t)^{N-1}}\ .
$$
Since $[R/N,2R/N]\subset [0,R]$ and the maximum $<\frac{N}{R}
\left(1+\frac{1}{N-2}\right)^{N-1}<\frac{4e}{3}\frac{N}{R}$, as $N\geq 5$, 
this and (\ref{e1238}) yield
$$
I_{1}\leq\frac{4e}{3}CN\frac{\mu_{\widetilde m_{i}}(B_{R}(\widetilde m_{i}))}
{R}d(m_{1},m_{2})\ .
$$
For the second term we get from (\ref{e1237})
$$
I_{2}\leq 
C\left(\beta_{N}\int_{R/N}^{R}\frac{\mu_{m_{i}}(B_{s}(m_{i}))}{s}
(R-s)^{N-1}ds\right)d(m_{1},m_{2})
$$
and by (\ref{e1238}) the factor in the brackets is at most 
$\mu_{\widetilde m_{i}}(B_{R}(\widetilde m_{i}))\cdot\frac{N}{R}$.
Hence, we have
$$
I_{2}\leq CN\frac{\mu_{\widetilde m_{i}}(B_{R}(\widetilde m_{i}))}{R}
d(m_{1},m_{2})\ .
$$
Further note that $d(m_{1},m_{2})\leq 
d_{N}(\widetilde m_{1},\widetilde m_{2})$. Hence, we obtain finally
the inequality
$$
|\mu_{\widetilde m_{1}}-\mu_{\widetilde m_{2}}|(B_{R}(\widetilde m_{i}))\leq
\left(1+\frac{4e}{3}\right)NC\frac{\mu_{\widetilde m_{i}}
(B_{R}(\widetilde m_{i}))}{R}d(\widetilde m_{1},\widetilde m_{2})
$$
whence $C_{N}\leq\left(1+\frac{4e}{3}\right)NC$.\ \ \ \ \ $\Box$

At the last step we introduce the desired metric space 
$(\widehat M,\widehat d)$ simply setting
$$
\widehat M:=M_{N}\times\Re\ (=M\times l_{1}^{N+1})\ \ \ {\rm and}\ \ \
\widehat d:=d_{N+1}
$$ 
with $N:=[\log_{2}D]+5$. Moreover, we introduce the family of
measures $\widehat {\cal F}:=\{\mu_{\widehat m}\}_{\widehat m\in\widehat M}$
by
$$
\mu_{\widehat m}:=\mu_{\widetilde m}\otimes\lambda\ ;
$$
here and below $\widehat m:=(\widetilde m,x)$ where $\widetilde m\in M_{N}$
and $x\in\Re$, and $\lambda$ is the Lebesgue measure on $\Re$.

By Lemma \ref{l128} the dilation function $\widehat D=D_{N+1}$ of family
$\widehat{\cal F}$ at point $1+\frac{1}{N+1}$ is estimated as
\begin{equation}\label{e1239}
\widehat D\left(1+\frac{1}{N+1}\right)\leq\frac{6}{5}e^{4}\ .
\end{equation}
Moreover, by Lemma \ref{l129} the consistency constant 
$\widehat C=C_{N+1}$ of $\widehat{\cal F}$ satisfies
\begin{equation}\label{e1240}
\widehat C\leq\left(1+\frac{4e}{3}\right)(N+1)C\ .
\end{equation}

Show now that family $\widehat{\cal F}=\{\mu_{\widehat m}\}$ satisfies
condition (\ref{e1413}) with constant $\widehat A$ satisfying
\begin{equation}\label{e1241}
\widehat A\leq \frac{6}{5}e^{4}(N+1)\ .
\end{equation}
In fact, $\mu_{\widehat m}=\mu_{\widetilde m}\otimes\lambda$ and by the
Fubini theorem we have for $0<R_{1}<R_{2}$
$$
\mu_{\widehat m}(B_{R_{2}}(\widehat m))-\mu_{\widehat m}(B_{R_{1}}
(\widehat m))
=2\int_{R_{1}}^{R_{2}}\mu_{\widetilde m}(B_{s}(\widetilde m))ds\leq
\frac{2R_{2}\mu_{\widetilde m}(B_{R_{2}}(\widetilde m))}{R_{2}}(R_{2}-R_{1})\ .
$$
Prove that for arbitrary $l>1$ and $R>0$
\begin{equation}\label{e1242}
R\mu_{\widetilde m}(B_{R}(\widetilde m))\leq\frac{lD_{N}(l)}{2(l-1)}
\mu_{\widehat m}(B_{R}(\widehat m))\ .
\end{equation}
Together with the previous inequality this will yield
$$
\mu_{\widehat m}(B_{R_{2}}(\widehat m))-
\mu_{\widehat m}(B_{R_{1}}(\widehat m))\leq\frac{lD_{N}(l)}{l-1}\cdot
\frac{\mu_{\widehat m}(B_{R_{2}}(\widehat m))}{R_{2}}(R_{2}-R_{1})\ ,
$$
that is, inequality (\ref{e1413}) for family $\{\mu_{\widehat m}\}$ will be
proved with
\begin{equation}\label{eq1243}
\widehat A\leq\frac{lD_{N}(l)}{l-1}\ .
\end{equation}
Finally choose here $l=1+\frac{1}{N}$ and use Lemma \ref{l128}. This gives
the required inequality (\ref{e1241}).

Hence, it remains to establish (\ref{e1242}). By the definition of
$D_{N}(l)$, see (\ref{e141}), we have for $l>1$
$$
\mu_{\widehat m}(B_{lR}(\widehat m))=2l\int_{0}^{R}
\mu_{\widetilde m}(B_{ls}(\widetilde m))ds\leq lD_{N}(l)\mu_{\widehat m}
(B_{R}(\widehat m))\ .
$$
On the other hand, replacing $[0,R]$ by $[l^{-1}R,R]$ we also have
$$
\mu_{\widehat m}(B_{lR}(\widehat m))\geq
2l\mu_{\widetilde m}(B_{R}(\widetilde m))
(R-l^{-1}R)= 2(l-1)R\mu_{\widetilde m}(B_{R}(\widetilde m))\ .
$$
Combining the last two inequalities to get (\ref{e1242}).
\begin{R}\label{rem210}
{\rm For the proofs of corollaries it is useful to single out the next two
inequalities
\begin{equation}\label{eq1244}
\widehat D(l)\leq lD_{N}(l)\ \ \ \ {\rm and}\ \ \ \
\widehat C\leq\frac{C_{N}}{l-1}D_{N}(l)\ .
\end{equation}
The first of them follows from the inequality next to (\ref{eq1243}).
To prove the second one, write for $i=1,2$
$$
|\mu_{\widehat m_{1}}-\mu_{\widehat m_{2}}|(B_{R}(\widehat m_{i}))\leq 2
\int_{0}^{R}|\mu_{\widetilde m_{1}}-\mu_{\widetilde m_{2}}|
(B_{s}(\widetilde m_{i}))\ \!ds\leq 2C_{N}
\mu_{\widetilde m_{i}}(B_{R}(\widetilde m_{i}))
d(\widetilde m_{1},\widetilde m_{2})\ .
$$
Combining this with inequality (\ref{e1242}) we obtain the second inequality
in (\ref{eq1244}). We will use inequalities (\ref{eq1243}) and (\ref{eq1244})
for $N=0$, i.e., for $C_{N}$ equals the consistency
constant $C$ for $(M,d)$ and $D_{N}(l)=D(l)$.}
\end{R}

Now use the main result of Part I for the case of 
$S\subset M\subset\widehat M$. We conclude from here that there exists an
extension operator $\widehat E\in Ext(S,\widehat M)$ with norm satisfying
the inequality
$$
||\widehat E||\leq 8\widehat A+\max\left(\frac{2(l+3)}{l-1},K(l)\right)
$$
where
$$
K(l)=6(\widehat A+\widehat C)\widehat D(l)(l+3)\ ,
$$
see (\ref{eq1225}) and (\ref{e1424}).

Choose here $l:=1+\frac{1}{N+1}$ and apply inequalities (\ref{e1239})-
(\ref{e1241}) with $N=[\log_{2} D]+5$. This yields
\begin{equation}\label{e1243}
||\widehat E||\leq a_{0}(C+a_{1})(\log_{2}D+6)
\end{equation}
with some $a_{0}\ (<7575)$ and $a_{1}\ (<15)$. 
Then the restriction of 
$\widehat Ef$ to
$M$ gives the required extension operator from $Ext(S,M)$ with the norm 
bounded by the  right-hand side of (\ref{e1243}).

The proof of Theorem \ref{te114} is complete.\ \ \ \ \ $\Box$\\
{\bf Proof of Corollary \ref{co222}.} According to (\ref{eq29})
the dilation function for $\{\mu_{m}\}$ satisfies
\begin{equation}\label{eq1236}
D(l)\leq al^{n}\ ,\ \ \ 1\leq l<\infty\ ,
\end{equation}
with $a\geq 1$ and $n\geq 0$. To derive the require estimate of
$\lambda(M)$ we first use inequality (\ref{e1424}) for the space
$(\widehat M,\widehat d)$ where $\widehat M:=M\times\Re$ and $\widehat m$,
$\widehat d$ and $\{\mu_{\widehat m}\}$ are defined as in the above proof,
i.e., $\widehat m:=(m,x)$ with $m\in M$ and $x\in\Re$, \
$\widehat d=d_{1}$ and
$\mu_{\widehat m}:=\mu_{\widetilde m}\otimes\lambda_{1}$. Then
\begin{equation}\label{eq1237}
\lambda(M)\leq 8\widehat A+\max\left(\frac{2(l+3)}{l-1},\widehat
K(l)\right)
\end{equation}
where
$$
\widehat K(l)=6(\widehat A+\widehat C)(l+3)\widehat D(l)
$$
and the quantities with the hat are estimated in (\ref{eq1243}) and 
(\ref{eq1244}) with $N=0$ (see Remark \ref{rem210}).
In particular, one has $\widehat D(l)\leq al^{n+1}$, 
$\widehat A\leq \frac{al^{n+1}}{l-1}$ and $\widehat
C\leq\frac{al^{n+1}}{l-1}C$. Taking in the last two inequalities
$l=1+\frac{1}{n+1}$ we have 
$$
\widehat A\leq ea (n+1)\ ,\ \ \ \widehat C\leq ea (n+1) C\ .
$$
Inserting this in (\ref{eq1237}) we get
$$
\lambda(M)\leq 8ea (n+1)+\max\left(\frac{2(l+3)}{l-1},\widehat
K(l)\right)
$$
where now
$$
\widehat K(l)\leq 6e a (n+1) (C+1)(l+3)\widehat D(l)\ ,
$$
and $\widehat D(l)\leq al^{n+1}$. Since $l:=1+\frac{1}{n+1}$, we
straightforwardly obtain the inequality
$$
\lambda(M)\leq 225 a^{2}(C+1)(n+1)\ .\ \ \ \ \ \Box
$$
{\bf Proof of Corollary \ref{c116}.} We follow the same
argument using now the inequality (\ref{eq1229}) and then choosing
for $n\geq 1$ the value $l=1+\frac{2}{3n}$. Then a straightforward
computation yields for this choice of $l$
$$
\lambda(M)\leq
8n+\max\left(\frac{2(l+3)}{l-1},(l+3)(C+2nl^{n-1})\right)\leq
24(n+C)\ .
$$
For $n<1$ we simply choose $l=2$.\ \ \ \ \ $\Box$\\
{\bf Proof of Theorem \ref{c118}.} {\bf (a)} Let, first, $p=\infty$.
Since the metric in
$M:=\oplus_{\infty}\{(M_{i},d_{i})\}_{1\leq i\leq N}$ is given by 
$d(m,m'):=\max_{1\leq i\leq
N}d_{i}(m_{i},m_{i}')$, the ball $B_{R}(m)$ of $M$ is the product
of balls $B_{R}(m_{i})$ of $M_{i}$, $1\leq i\leq N$. Therefore for
a family of doubling measures $\{\mu_{m}\}_{m\in M}$ given by the
tensor product
\begin{equation}\label{eq1238}
\mu_{m}:=\bigotimes_{i=1}^{N}\mu_{m_{i}}^{i}\ ,\ \ \
m=(m_{1},\dots, m_{N})\ ,
\end{equation}
we get
\begin{equation}\label{eq1239}
\mu_{m}(B_{R}(m))=\prod_{i=1}^{N}\mu_{m_{i}}^{i}(B_{R}(m_{i}))\ .
\end{equation}
Hence for the dilation function (\ref{e141}) of the family
$\{\mu_{m}\}_{m\in M}$ we get
\begin{equation}\label{eq1240}
D(l)=\prod_{i=1}^{N} D_{i}(l)
\end{equation}
where $D_{i}$ is the dilation function of $\{\mu_{m}^{i}\}_{m\in
M_{i}}$. In particular, $\{\mu_{m}\}_{m\in M}$ satisfies the
uniform doubling condition (\ref{eq27}) with $D:=D_{1}\cdots D_{N}$.

Check that the condition (\ref{eq28}) holds for this family with
the constant
\begin{equation}\label{eq1241}
\widetilde C_{\infty}:=\left(\prod_{i=1}^{N}
K_{i}\right)\sum_{i=1}^{N}C_{i}\ .
\end{equation}

In fact, the identity
\begin{equation}\label{eq1241'}
\mu_{m}-\mu_{\widetilde
m}=\sum_{i=1}^{N}(\otimes_{j=1}^{i-1}\mu_{\widetilde
m_{j}}^{j})\otimes (\mu_{m_{i}}^{i}-\mu_{\widetilde
m_{i}}^{i})\otimes(\otimes_{j=i+1}^{N}\mu_{m_{j}}^{j})
\end{equation}
together with (\ref{eq1239}), and (\ref{eq28})  and
$K_{j}$-uniformity of $\{\mu_{m}^{j}\}_{m\in M_{j}}$ implies that
for $\widehat m=m $ or $\widetilde m$
$$
|\mu_{m}-\mu_{\widetilde m}|(B_{R}(\widehat
m))\leq\sum_{i=1}^{N}\left(\prod_{j\neq i}
K_{j}\right)C_{i}\frac{\mu_{m}(B_{R}(m))}{R}d_{i}(m_{i},\widetilde
m_{i})\leq\widetilde C_{\infty}\frac{\mu_{m}(B_{R}(m))}{R}d(m,\widetilde
m).
$$

Thus $\oplus_{\infty}\{(M_{i},d_{i})\}_{1\leq i\leq N}$ is of
pointwise homogeneous type with respect
to the family (\ref{eq1238}) with the optimal constants bounded by $D$ and
$\widetilde C_{\infty}$ (and so we have the required estimate for 
$\lambda(M)$ in this case).

Let now $\mu_{m}^{i}(B_{R}(m))=\gamma_{i}R^{n_{i}}$ for some
$\gamma_{i},n_{i}>0$ and all $m\in M_{i}$ and $R>0$, $1\leq i\leq
N$. In this case $\{\mu_{m}^{i}\}_{m\in M_{i}}$ is clearly
$K_{i}$-uniform with $K_{i}=1$. Moreover, by (\ref{eq1239})
$$
\mu_{m}(B_{R}(m))=\gamma R^{n}\ ,\ \ \ n:=\sum_{i=1}^{N} n_{i}\ .
$$
Hence $M$ equipped with the family (\ref{eq1238}) satisfies the
conditions of Corollary \ref{c116} with this $n$ and
$C=\sum_{i=1}^{N} C_{i}$, see (\ref{eq1241}). Applying this
corollary we get
$$
\lambda(\oplus_{\infty}\{M_{i}\}_{1\leq i\leq N})\leq 
24\sum_{i=1}^{N}(n_{i}+C_{i}).
$$
{\bf (b)} Let now $1\leq p<\infty$. In this case we cannot estimate 
the optimal constants $C$ and $D$ for the space
\begin{equation}\label{eq1253}
(M,d):=\oplus_{p}\{(M_{i},d_{i})\}_{1\leq i\leq N}
\end{equation}
directly. To overcome this difficulty we use the argument of Theorem
\ref{te114} and isometrically embed this space into the space
$$
(\widehat M,\widehat d):=(M,d)\oplus_{1}l_{1}^{a}
$$
with a suitable $a$. Hence, a point $\widehat m\in \widehat M$ is an
$(N+a)$-tuple
$$
\widehat m:=(m,x):=(m_{1},\dots, m_{N},x_{1},\dots, x_{a})
$$
with $m\in\prod_{i=1}^{N}M_{i}$ and $x\in\Re^{a}$. Moreover, the metric
$\widehat d$ is given by
$$
\widehat d(\widehat m,\widehat m'):=\left(\sum_{i=1}^{N}
d_{i}(m_{i},m_{i}')^{p}\right)^{1/p}+\sum_{i=1}^{a}|x_{i}-x_{i}'|.
$$
Endow $\widehat M$ with a family of measures given by the tensor product
$$
\mu_{\widehat m}:=\mu_{m}\otimes\Lambda_{a} ,\ \ \ \widehat m\in\widehat M,
$$
where $\Lambda_{a}$ is the Lebesgue measure on $\Re^{a}$ and
$\mu_{m}:=\otimes_{i=1}^{N}\mu_{m_{i}}^{i}$. \\
We will show that $\lambda(\widehat M)$ is bounded as required in Theorem
\ref{c118}. This immediately gets the desired estimate for $\lambda(M)$
and completes the proof of the theorem.

To accomplish this we need
\begin{Lm}\label{le1211}
The optimal uniform doubling constant $D$ of the family $\{\mu_{m}\}_{m\in M}$
satisfies
$$
D\leq\prod_{i=1}^{N}D_{i}.
$$
Recall that $D_{i}$ is the optimal uniform doubling constant of
$\{\mu_{m_{i}}^{i}\}_{m_{i}\in M_{i}}$.
\end{Lm}
{\bf Proof} (induction on $N$). For the $\mu_{m}$-measure of the ball
$$
B_{2R}(m):=\{m'\in M\ :\ \sum_{i=1}^{N}d_{i}(m_{i},m_{i}')^{p}\leq (2R)^{p}\}
$$
we get by the Fubini theorem:
$$
\mu_{m}(B_{2R}(m))=\int_{d^{1}<(2R)^{p}}d\mu^{1}(m')\int_{d_{1}<
(2R)^{p}-d^{1}}d\mu_{1}(m_{1}').
$$
Here we set for simplicity:
$$
d^{1}:=\sum_{i=2}^{N}d_{i}(m_{i},m_{i}')^{p},\ \ \ d_{1}:=d(m_{1},m_{1}')^{p},
\ \ \
\mu^{1}:=\bigotimes_{i=2}^{N}\mu_{m_{i}}^{i},\ \ \ \mu_{1}:=\mu_{m_{1}}^{1}.
$$
The second integral is the $\mu_{1}$-measure of the ball $B_{2\rho}(m_{1})$
where $\rho:=\sqrt[p]{R^{p}-2^{-p}d^{1}}$ which is bounded by 
$D_{1}\mu_{1}(B_{\rho}(m_{1}))$. This and the Fubini theorem imply
$$
\begin{array}{c}
\displaystyle
\mu_{m}(B_{2R}(m))\leq D_{1}\int_{2^{-p}d^{1}<R^{p}}d\mu^{1}(m')
\int_{d_{1}<R^{p}-2^{-p}d^{1}}d\mu_{1}(m_{1}')=\\
\\
\displaystyle
D_{1}\int_{d_{1}<R^{p}}
d\mu_{1}(m_{1}')\int_{d^{1}<(2R)^{p}-2^{p}d_{1}}d\mu^{1}(m').
\end{array}
$$
By the induction hypothesis the inner integral in the right-hand side
is bounded by
$$
\left(\prod_{i=2}^{N}D_{i}\right)\mu^{1}(B_{\sqrt[p]{R^{p}-d_{1}}}
(m_{2},\dots, m_{N}))=
\prod_{i=2}^{N}D_{i}\int_{d^{1}<R^{p}-d_{1}}d\mu^{1}(m').
$$
Combining this with the previous inequality to get the required result:
$$
\mu_{m}(B_{2R}(m))\leq\left(\prod_{i=1}^{N}D_{i}\right)\mu_{m}(B_{R}(m)).
\ \ \ \ \ \Box
$$

Using Lemma \ref{le1211} we estimate now the dilation function
$D_{a}(s)$ of the family $\{\mu_{\widehat m}\}$. Recall that for $s>1$
\begin{equation}\label{eq1254}
D_{a}(s):=\sup_{\widehat m\in\widehat M}\left\{
\frac{\mu_{\widehat m}(B_{sR}(\widehat m))}{\mu_{\widehat m}
(B_{R}(\widehat m))}\right\}
\end{equation}
To this end we simply apply to this setting Lemma \ref{l128} with $D$
replaced by $\prod_{i=1}^{N}D_{i}$ and $N$ by $a$. This gets
\begin{Lm}\label{le1212}
If $a\geq [\log_{2}\prod_{i=1}^{N}D_{i}]+5$, then
$$
D_{a}(1+1/a)\leq\frac{6}{5}e^{4}.\ \ \ \ \ \Box
$$
\end{Lm}

Now we estimate the consistency constant for the family 
$\{\mu_{\widehat m}\}_{\widehat m\in\widehat M}$, see Definition \ref{d218}.
To this goal we use (\ref{eq1241'}) for $\mu_{\widehat m}-\mu_{\widehat m'}$
and then apply the Fubini theorem to have for $\widehat m'':=\widehat m$ or
$\widehat m'$
\begin{equation}\label{eq1255}
\begin{array}{c}
\displaystyle
|\mu_{\widehat m}-\mu_{\widehat m'}|(B_{R}(\widehat m''))\leq
\\
\\
\displaystyle
\sum_{i=1}^{N}\int_{\delta_{a}<R}d\Lambda_{a}\int_{d^{i}<(R-\delta_{a})^{p}}
d\mu_{i}'d\mu_{i}\int_{d_{i}<(R-\delta_{a})^{p}-d^{i}}
d|\mu_{m_{i}}^{i}-\mu_{m_{i}'}^{i}| .
\end{array}
\end{equation}
Here we use the notations:
$$
\begin{array}{c}
\displaystyle
\delta_{a}:=\sum_{j=1}^{a}|x_{j}-x_{j}''|,\ \ \
d^{i}:=\sum_{j\neq i}d_{j}(m_{j}'',m_{j})^{p},\ \ \ 
d_{i}:=d(m_{i}'',m_{i})^{p},\\
\\
\displaystyle
\mu_{i}':=\bigotimes_{j<i}\mu_{m_{j}'}^{j},\ \ \ \mu_{i}:=\bigotimes_{j>i}
\mu_{m_{j}}^{j}.
\end{array}
$$
Recall that $\widehat m=(m,x)\in M\times\Re^{a}$.

The inner integral in the $i$-th term of the right-hand side of 
(\ref{eq1255}) equals $|\mu_{m_{i}}^{i}-\mu_{m_{i}'}^{i}|
(B_{\rho}(m_{i}''))$ where $\rho:=\sqrt[p]{(R-\delta_{a})^{p}-d^{i}}$.
Replacing here $\rho$ by $\rho_{a}:=\sqrt[p]{(R_{a}-\delta_{a})^{p}-d^{i}}$
with $R_{a}:=(1+\frac{1}{a})R$ and applying the consistency inequality for
$(M_{i},d_{i})$ we then bound this inner integral by
$$
\frac{C_{i}\ \!\mu_{m_{i}''}^{i}(B_{\rho_{a}}(m_{i}''))}{\rho_{a}}
d_{i}(m_{i},m_{i}').
$$
Since $d^{i}\leq (R-\delta_{a})^{p}$, the denominator here is at least
$R_{a}-R=\frac{1}{a}R$. Therefore the inner integral is bounded by 
$$
\frac{a\ \! C_{i}\ \!
d_{i}(m_{i},m_{i}')}{R}\int_{d_{i}<(R_{a}-\delta_{a})^{p}-d^{i}}
d\mu_{m_{i}''}^{i}.
$$
Inserting this in (\ref{eq1255}) and replacing there $R$ by $R_{a}$ we
get
$$
|\mu_{\widehat m}-\mu_{\widehat m'}|(B_{R}(\widehat m''))\leq
\frac{a}{R}\sum_{i=1}^{N}C_{i}d_{i}(m_{i},m_{i}')
\int_{B_{R_{a}}(\widehat m'')}d\Lambda_{a}\ \!d\mu_{i}'
\ \!d\mu_{i}\ \!d\mu_{m_{i}''}^{i}.
$$
To replace in this inequality each $\mu_{m_{j}'}^{j}$ (or
$\mu_{m_{j}}^{j}$) by $\mu_{m_{j}''}^{j}$ we now use $K_{j}$-uniformity of
the family $\{\mu_{m_{j}}^{j}\}_{m_{j}\in M_{j}}$, see Definition
\ref{de224}. Applying this to the right-hand side of the previous inequality
and recalling definition (\ref{eq1254}) we estimate the $i$-th integral there
by
$$
\begin{array}{c}
\displaystyle
\left(\prod_{i=1}^{N}K_{i}\right)\int_{B_{R_{a}}(\widehat m'')}
d\Lambda_{a}d\mu_{m''}=\left(\prod_{i=1}^{N}K_{i}\right)\mu_{\widehat m''}
(B_{R_{a}}(\widehat m''))\leq
\\
\\
\displaystyle
 D_{a}(1+1/a)\left(\prod_{i=1}^{N}K_{i}\right)
\mu_{\widehat m''}(B_{R}(\widehat m'')).
\end{array}
$$
Combining with the previous inequality we get for $\widehat m''=\widehat m$
or $\widehat m'$
$$
|\mu_{\widehat m}-\mu_{\widehat m'}|(B_{R}(\widehat m''))\leq
\frac{aD_{a}(1+1/a)}{R}\left(\prod_{i=1}^{N}K_{i}\right)\left(
\sum_{i=1}^{N}C_{i}d_{i}(m_{i},m_{i}')\right)
\mu_{\widehat m''}(B_{R}(\widehat m'')).
$$
By the H\"{o}lder inequality the sum in the brackets is at most
$$
\left(\sum_{i=1}^{N}C_{i}^{q}\right)^{1/q}\left(\sum_{i=1}^{N}
d_{i}(m_{i},m_{i}')^{p}\right)^{1/p}=:
\left(\sum_{i=1}^{N}C_{i}^{q}\right)^{1/q}d(m,m');
$$ 
here $\frac{1}{p}+\frac{1}{q}=1$. Hence the consistency constant
$\widehat C$ of the family $\{\mu_{\widehat m}\}_{\widehat m\in\widehat M}$
satisfies
\begin{equation}\label{eq1257}
\widehat C\leq a D_{a}(1+1/a)\left(\prod_{i=1}^{N}K_{i}\right)
\left(\sum_{i=1}^{N}
C_{i}^{q}\right)^{1/q}.
\end{equation}

Choose now $a:=[\log_{2}\prod_{i=1}^{N}D_{i}]+5$ and use (\ref{eq1243})
for the space $(\widehat M,\widehat d)$ equipped with the family
$\{\mu_{\widehat m}\}_{\widehat m\in\widehat M}$. Since $\widehat A$ in
(\ref{eq1243}) is bounded by $\frac{sD_{a}(s)}{s-1}$, with 
$s=1+1/a$, we therefore get from Lemma \ref{le1212}
$$
\widehat A\leq\frac{6}{5}e^{4}\left(\log_{2}\left(\prod_{i=1}^{N}D_{i}
\right)+6\right).
$$
Combining Lemma \ref{le1212} with (\ref{eq1257}) and the above inequality
we finally obtain the required result (see (\ref{e1424}))
$$
\lambda(\widehat M)\leq c_{0}(\widetilde C_{p}+1)
\left(\log_{2}\left(\prod_{i=1}^{N}D_{i}\right)+1\right)
$$
with $\widetilde C_{p}:=\left(\sum_{i=1}^{N}
C_{i}^{q}\right)^{1/q}\left(\prod_{i=1}^{N}K_{i}\right)$ and 
$\frac{1}{p}+\frac{1}{q}=1$.
\ \ \ \ \ $\Box$

\sect{Appendix: A Duality Theorem}
Our goal is to prove Theorem A of section 5, that is, we have to find a
Banach space $V$ such that its dual
\begin{equation}\label{ap1}
V^{*}=Lip_{0}(M)
\end{equation}
and all evaluations $\delta_{m}:\phi\to\phi(m)$, $\phi\in Lip_{0}(M)$,
$m\in M$, belong to $V$.

We introduce this as a (closed) subspace of the Banach space $l_{\infty}(B)$
where $B$ is the closed unit ball of $Lip_{0}(M)$. To this end, define a
map
$$
\Phi: M\to l_{\infty}(B)
$$
given for $m\in M$ by
\begin{equation}\label{ap2}
\Phi(m)(b):=b(m)\ ,\ \ \ b\in B\ .
\end{equation}
As all functions of $Lip_{0}(M)$ vanish at a prescribed point $m^{*}$, we
have for $b\in B$
$$
|b(m)-b(m^{*})|\leq d(m,m^{*})\ ,
$$
and $\Phi(M)$ is, actually, a subset of $l_{\infty}(B)$.

We now define the desired Banach space by
\begin{equation}\label{ap3}
V:=\overline{span\ \Phi(M)}\ ,
\end{equation}
the closure in $l_{\infty}(B)$ of the linear span of $\Phi(M)$ endowed by
the norm induced from $l_{\infty}(B)$.

Then introduce the required isometry $I$ of the dual $V^{*}$ to
$V$ onto $Lip_{0}(M)$ as the pullback of the map $\Phi:M\to V$;
that is to say, we let for $l\in V^{*}$
\begin{equation}\label{ap4}
I(l)(m):=l(\Phi(m))\ ,\ \ \ m\in M\ .
\end{equation}
{\bf Assertion 1.} {\em The linear operator $I$ is an injection.}

In fact, if $I(l)=0$ for some $l\in V^{*}$, then $l|_{\Phi(M)}=0$ and,
by (\ref{ap3}), $l=0$.\\
{\bf Assertion 2.} {\em It is true that}
\begin{equation}\label{ap5}
Lip_{0}(M)\subset I(V^{*})\ , \ \ \ {\rm and}\ \ \
\end{equation}
\begin{equation}\label{ap6}
||I||:=\sup\{||I(l)||_{Lip(M)}\ :\ ||l||\leq 1\}\geq 1\ .
\end{equation}

Actually, $\Phi(m^{*})=0$ and therefore each function $I(l)$ vanishes at
$m^{*}$. Let now $b\in B$ and $\pi_{b}:l_{\infty}(B)\to\Re$ be the
canonical projection given by
\begin{equation}\label{ap7}
\pi_{b}(x):=x(b)\ ,\ \ \ x\in l_{\infty}(B)\ .
\end{equation}
Then by (\ref{ap4}) and (\ref{ap2})
\begin{equation}\label{ap7'}
I(\pi_{b}|_{V})(m)=\pi_{b}(\Phi(m))=\Phi(m)(b)=b(m)
\end{equation}
for all $m\in M$. Since $\pi_{b}|_{V}\in V^{*}$ and $b$ is an arbitrary
element of the unit ball in $Lip_{0}(M)$, the embedding (\ref{ap5}) holds.
Besides,
$$
||I||\geq\sup_{b\in B}||I(\pi_{b}|_{V})||_{Lip_{0}(M)}=1\ ,
$$
and (\ref{ap6}) is also true.\\
{\bf Assertion 3.} {\em It is true that}
\begin{equation}\label{ap8}
I(V^{*})\subset Lip_{0}(M)\ ,\ \ \ {\rm and}
\end{equation}
\begin{equation}\label{ap9}
||I||\leq 1\ .
\end{equation}

Let $l\in V^{*}$ and $m_{1}\neq m_{2}\in M$. We have to show that
\begin{equation}\label{ap10}
|I(l)(m_{1})-I(l)(m_{2})|\leq d(m_{1},m_{2})\ \!||l||\ ;
\end{equation}
as, in addition, $I(l)(m^{*})=0$, this will prove the assertion.

To establish (\ref{ap10}), extend $l$ by the Hahn-Banach theorem to
$\widehat l\in l_{\infty}(B)^{*}$. Hence,
\begin{equation}\label{ap11}
\widehat l|_{V}=l\ \ \ {\rm and}\ \ \ ||\widehat l||=||l||\ .
\end{equation}
Now, using the Gelfand transform we identify $l_{\infty}(B)$ with the space
$C(\beta B)$ of continuous functions on the space $\beta B$ of maximal
ideals of the Banach algebra $l_{\infty}(B)$. In fact, $\beta B$ is the
Stone-\v{C}ech compactification of $B$ regarded as a topological space endowed
by the {\em discrete topology}, see, e.g., [Lo]. By the F. Riesz theorem
there exists a bounded (regular) Borel measure $\mu_{\widehat l}$ on
$\beta B$ such that
\begin{equation}\label{ap12}
\widehat l(x):=\int_{\beta B}g(x)d\mu_{\widehat l}\ ,\ \ \ x\in l_{\infty}(B)
\ ;
\end{equation}
here $g(x)\in C(\beta B)$ is the Gelfand transform of $x$; recall that in this
case $g(x)$ is the continuous extension of $x$ from $B$ to $\beta B$.

Let now $x\in V$, and $\{U_{i}\}$ be a finite open cover of $\beta B$ such
that the oscillation of $g(x)$ on each $U_{i}$ is at most $\epsilon$.
Since $B$ is dense in $\beta B$, every $U_{i}$ contains a point $b_{i}\in B$, and therefore
$$
|g(x)(\omega)-g(x)(b_{i})|=|g(x)(\omega)-x(b_{i})|<\epsilon
$$
for every $\omega\in U_{i}$. Hence for such an $x$
$$
\left|\int_{\beta B}g(x)d\mu_{\widehat l}-\sum
x(b_{i})\mu_{\widehat l}(U_{i})\right|<\epsilon\cdot Var\ \!\mu_{\widehat l}=
\epsilon||l||\ ,
$$
see (\ref{ap11}). The sum here can be written as $l_{\epsilon}(x)$ where
$l_{\epsilon}\in V^{*}$ is given by
$$
l_{\epsilon}:=\sum(\pi_{b_{i}}|_{V})\mu_{\widehat l}(U_{i})\ ,
$$
see (\ref{ap7}) and (\ref{ap7'}). So, together with (\ref{ap12})
and (\ref{ap11}) this leads to the estimate
$$
|l(x)-l_{\epsilon}(x)|<\epsilon ||l||\ .
$$
Choose here $x:=\Phi(m_{1})-\Phi(m_{2})$ and use
(\ref{ap4}). This implies that
$$
|(I(l)(m_{1})-I(l)(m_{2}))-
(I(l_{\epsilon})(m_{1})-I(l_{\epsilon})(m_{2}))|<\epsilon ||l||\ .
$$
Besides, by the definition of $l_{\epsilon}$
$$
I(l_{\epsilon})(m)=\sum b_{i}(m)\mu_{\widehat l}(U_{i})\ ,\ \ \ m\in M\ ,
$$
and so $I(l_{\epsilon})\in Lip_{0}(M)$ and
$$
||I(l_{\epsilon})||_{Lip_{0}(M)}\leq Var\ \!\mu_{\widehat l}=||l||\ ,
$$
see (\ref{ap11}) and (\ref{ap12}). Together with the previous inequality
this yields
$$
|I(l)(m_{1})-I(l)(m_{2})|\leq d(m_{1},m_{2})||l||+\epsilon ||l||\ .
$$
Letting $\epsilon$ to 0, we conclude that $I(l)\in Lip_{0}(M)$ and
$||I(l)||_{Lip_{0}(M)}\leq ||l||$. This proves (\ref{ap8}) and (\ref{ap9}).

Now the Assertions 1-3 prove that $I$ is an isometry of $V^{*}$ onto
$Lip_{0}(M)$. Besides, the evaluation functionals
$\delta_{m}:\phi\to\phi(m)$, $\phi\in Lip_{0}(M)$,  can be presented as
$\delta_{m}=\Phi(m)$, see (\ref{ap2}), and therefore belong to $V$.

The proof is complete.\ \ \ \ \ $\Box$
\begin{R}\label{ap}
{\rm Let $U\subset V$ be the closure of the convex hull of the set
\penalty-10000
$\{v_{x,y}:=(x-y)/d(x,y)\ :\ x,y\in\Phi(M),\ x\neq y\}$.
It is easily seen that $U$ is the closed unit ball of $V$.}
\end{R}


\end{document}